\documentclass[12pt]{article}

\usepackage{amssymb}

\usepackage{latexsym,amsmath,amsfonts,amscd}
\usepackage{xfrac}
\usepackage{graphics}
\usepackage{caption}
\usepackage{subcaption}
\usepackage{epstopdf}
\usepackage{color}
\usepackage{changebar}
\usepackage{indentfirst}
\usepackage{verbatim}
\usepackage{xcolor}
\usepackage{hyperref}
\usepackage{amssymb}
\usepackage{lscape}
\usepackage{array,multirow}
\usepackage{booktabs}

\renewcommand{\d}{{\rm d}}			
\newcommand{\CFL}{\mathsf{CFL}}		
\newtheorem{definition}{Definition}
\newtheorem{remark}{Remark}
\newtheorem{proof}{Proof}

\usepackage{authblk}
\usepackage[margin=3.0cm]{geometry}

\title{Multiscale constitutive framework of 1D blood flow modeling: Asymptotic limits and numerical methods}

\author[$\dagger$]{Giulia Bertaglia \footnote{Corresponding author. Email address: \textit{giulia.bertaglia@unife.it}}}
\author[$\ddagger$]{Lorenzo Pareschi}
\affil[$\dagger$]{\small Department of Environmental and Prevention Sciences, University of Ferrara, Italy}
\affil[$\ddagger$]{\small Department of Mathematics and Computer Science, University of Ferrara, Italy}

\begin{document}

\maketitle

\begin{abstract}
In this paper, a multiscale constitutive framework for one-dimensional blood flow modeling is presented and discussed. By analyzing the asymptotic limits of the proposed model, it is shown that different types of blood propagation phenomena in arteries and veins can be described through an appropriate choice of scaling parameters, which are related to distinct characterizations of the fluid-structure interaction mechanism (whether elastic or viscoelastic) that exist between vessel walls and blood flow. In these asymptotic limits, well-known blood flow models from the literature are recovered. Additionally, by analyzing the perturbation of the local elastic equilibrium of the system, a new viscoelastic blood flow model is derived. The proposed approach is highly flexible and suitable for studying the human cardiovascular system, which is composed of vessels with high morphological and mechanical variability. The resulting multiscale hyperbolic model of blood flow is solved using an asymptotic-preserving Implicit-Explicit Runge-Kutta Finite Volume method, which ensures the consistency of the numerical scheme with the different asymptotic limits of the mathematical model without affecting the choice of the time step by restrictions related to the smallness of the scaling parameters. Several numerical tests confirm the validity of the proposed methodology, including a case study investigating the hemodynamics of a thoracic aorta in the presence of a stent.
\end{abstract}


\tableofcontents

\section{Introduction}
The modeling of blood flow has undergone considerable development in recent years thanks to the interest of numerous researchers who have expanded its treatment, focusing on various fundamental aspects and issues, in order to arrive at an increasingly reliable description of the hemodynamics of the circulatory system \cite{formaggia2009,quarteroni2017}. 

Several studies have already shown that, in general, one-dimensional (1D) modeling coupled with lumped-parameter, zero-dimensional (0D) models, derived from full three-dimensional (3D) models by means of simplifying assumptions about flow, structure, and their interaction, is sufficient to obtain realistic and accurate numerical results, particularly when the flow is predominantly unidirectional \cite{xiao2014,pfaller2022,vidotto2019}. Moreover, in contrast to 3D simulations with prohibitively high computational costs, 1D models allow for the investigation of the hemodynamics of the entire main circulatory system \cite{alastruey2011,mynard2015,piccioli2021}.

One of the challenges presented by the modeling of the cardiovascular system is that blood vessels exhibit high morphological and mechanical variability, interacting with the blood flow to give rise to complex fluid-structure interactions (FSI). Moreover, in unhealthy cases, this variability is further emphasized due to the possible presence of calcifications, stenosis, aneurysms, or even prostheses such as grafts or stents \cite{formaggia2002,sherwin2003,ramella2022}.
To model the FSI occurring between blood and vessel walls, an appropriate constitutive model, which relates pressure to area variations, needs to be considered. To this end, it must be remembered that the smooth muscle cells that constitute the intermediate layer of vessels impart a viscoelastic behavior to the wall, which assumes a key role when high frequencies are dominant \cite{alastruey2012a,coccarelli2021}. In contrast, when stress is applied very slowly, viscous aspects do not occur, and the wall behaves mostly elastically. Although vessel viscosity is often neglected in blood flow models for simplicity, there is a growing number of contributions showing the advantages of modeling the mechanical behavior of the vessel wall using a viscoelastic rheological characterization, based on linear or quasilinear viscoelasticity or more complex nonlinear models \cite{valdez-jasso2009,alastruey2012a,ghigo2016,bertaglia2020a}.

In this complex setting, the paper aims to extend the treatment of the blood flow model first presented in \cite{bertaglia2020}, providing a multiscale framework of 1D blood flow modeling with a viscoelastic constitutive characterization of vessel walls. Indeed, by analyzing the asymptotic limits of the system, it will be shown that by adopting a linear viscoelastic constitutive relation of the Standard Linear Solid type, the proposed model can lead to a very flexible and accurate description of many possible mechanical behaviors of vessel walls, recovering different characterizations (e.g., the purely elastic model and other widely used linear viscoelastic models) with an appropriate choice of the scaling parameters involved. In addition, a new viscoelastic constitutive model will be derived by analyzing the perturbation of the local elastic equilibrium of the system. From a numerical point of view, a third order asymptotic-preserving (AP) Implicit-Explicit (IMEX) Runge-Kutta Finite Volume scheme is considered, which ensures consistency of the numerical method in all the asymptotic limits of the model (i.e., AP property) and whose time step size is not affected by the smallness of the scaling parameters.

The rest of the manuscript is organized as follows. In Section \ref{section_mathematicalmodel}, all modeling is presented and discussed, with emphasis toward the viscoelastic rheological characterization of vessel walls and the analysis of the asymptotic limits of the proposed augmented blood flow model. Section \ref{section_numericalmethod} is devoted to the presentation of the chosen numerical scheme and the proof of the fundamental AP property, together with the implementation of boundary conditions. In Section \ref{section_numericalresults}, several numerical tests are performed to validate the proposed methodology, including accuracy analysis, Riemann problems and a multiscale case study of a thoracic aorta with a stent implanted. Finally, some conclusions are drawn in Section \ref{section_conclusions}.

\section{One-dimensional blood flow modeling}
\label{section_mathematicalmodel}
The standard 1D mathematical model for blood flow, valid for medium to large-size vessels, is obtained averaging the incompressible Navier-Stokes equations over the cross-section, under the assumption of axial symmetry of the vessel and of the flow, obtaining the well established equations of conservation of mass and momentum \cite{formaggia2009}:
\begin{subequations}
\begin{align}
	&\frac{\partial A}{\partial t} + \frac{\partial(Au)}{\partial x} = 0 \label{eq.contST}\\ 
	&\frac{\partial (Au)}{\partial t} + \frac{\partial(Au^2)}{\partial x} + \dfrac{A}{\rho} \frac{\partial p}{\partial x} = 0 , \label{eq.momST}
\end{align}
\label{eq.cont&mom}	
\end{subequations}
with the choice, in the present work, to neglect friction losses.
Here \(A(x,t)\) is the cross-sectional area of the vessel, \(u(x,t)\) is the averaged fluid velocity, \(p(x,t)\) is the averaged fluid pressure, \(\rho\) is the density of the fluid and \(x\) and \(t\) are respectively space and time. 

To close this system of partial differential equations (PDEs), a tube law, representative of the interaction between vessel wall displacement (through the cross-sectional area \(A\)) and blood pressure \(p\), is required. To this end, the mechanical behavior of the vessel wall must be associated with a constitutive model, which relates stress and strain of the material as much realistically as possible.
\subsection{Elastic constitutive law}
\label{section_elastictubelaw}
In the simplest case, the pressure-area relationship is defined considering a perfectly elastic behavior of the vessel wall, hence the latter behaves like a simple linear spring characterized by its Young (elastic) modulus $E$. 
The constitutive equation of a linear elastic solid coincides with Hooke’s law, which is expressed as a linear relationship between stress $\sigma(t)$ and strain $\epsilon(t)$:
\begin{equation}
\label{hooke_simple}
\sigma = E \epsilon .
\end{equation}
We consider that the deformation of the material is geometrically related to the cross-sectional area through equation 
\begin{equation}
\label{epsilon-A}
\varepsilon = \alpha^m - \alpha^n ,
\end{equation}
where \(\alpha = A/A_0\) is the non-dimensional cross-sectional area scaled with respect to \(A_0(x)\), equilibrium cross-sectional area, and \(m\) and \(n\) are specific parameters related to the behavior of the vessel wall, whether artery or vein \cite{bertaglia2020}. Involving Barlow's formula, 
\begin{equation}
\label{Barlow}
\sigma = W(p-p_{0}) ,
\end{equation} 
where \(p_{0}(x)\) is the equilibrium pressure and $W(x)$ is a parameter depending on the wall thickness $h_0$ (here assumed to be always constant in space) and the equilibrium inner radius of the wall $R_0(x)$, which can have again different definitions if dealing with arteries or veins \cite{muller2013}, we can re-write Hooke's law in the following elastic constitutive tube law, widely adopted in literature \cite{formaggia2003,matthys2007,muller2013}:
\begin{equation}
p = p_{0} + \frac{E}{W} \left(\alpha^m - \alpha^n\right) .
\label{elastictubelaw}
\end{equation}
Here
\begin{equation*}
\psi = \frac{E}{W} \left(\alpha^m - \alpha^n\right)
\end{equation*}
is the elastic contribution of the transmural pressure.

If dealing with arteries, this constitutive law corresponds to the so-called Laplace law.
In contrast, when dealing with veins, their possible collapse in case of large negative transmural pressures needs to be considered \cite{toro2013,murillo2019}. The collapsed state for veins is identified by a cross-sectional area assuming a buckled, dumbbell shape configuration, in which opposite sides of the interior wall touch each other, still leaving some fluid flow in the two extremes. 
This particular aspect leads to the assumption of different parameters for the mechanical characterization of the wall behavior.
Following \cite{muller2013}, we have:
\begin{equation}
  W =
    \begin{cases}
      \frac{R_0}{h_0} & \text{if artery}\vspace{0.2cm}\\ 
      \frac{12 R_0^3}{h_0^3} & \text{if vein}
    \end{cases}       \qquad
 m =
    \begin{cases}
      1/2 & \text{if artery}\\
      10 & \text{if vein}
    \end{cases}   \qquad
 n =
    \begin{cases}
      0 & \text{if artery}\\
      -3/2 & \text{if vein}\,.
    \end{cases}  
\end{equation}
We remark here that, generally, one needs to choose $m > 0$ and $n \in [-2, 0]$ in order to preserve desirable mathematical properties of the PDE system \cite{toro2013,piccioli2021}.

\subsection{Linear viscoelastic constitutive laws}
\label{section_SLSM}
Even though mathematical models of blood circulation frequently neglect the viscous component of the vessel wall, it is well known that blood vessels (and living tissues in general) exhibit viscoelastic properties \cite{wang2016}. Viscoelastic effects are simulated in literature using different (more or less complex) rheological models, whether linear or not \cite{alastruey2011,holenstein1980,bessems2008,valdez-jasso2009,hasan2021,ghigo2016,kim2022}. With a simple but still effective choice, we can close system \eqref{eq.cont&mom} by considering a linear viscoelastic model as representative of the fluid-structure interaction mechanics of blood with the vessel wall. In general, a constitutive relation of linear viscoelasticity is built up considering the material as a sum of linear elastic springs, each one defined by a Young modulus $E$, and linear viscous dash-pots, characterized by a viscosity coefficient $\eta$, to take into account also the time dependent relaxation of the wall and its damping effect on pressure waves.
\begin{figure}
\centering
\begin{subfigure}[t]{0.25\textwidth}
\centering
\vspace{-2.05cm}
\includegraphics[width=\textwidth]{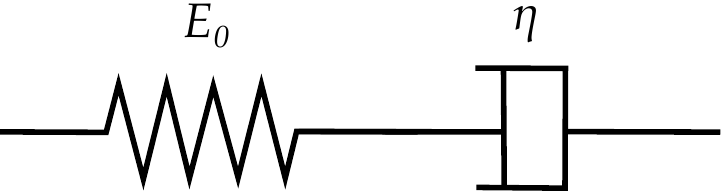}
\vspace{0.55cm}
\caption{Maxwell constitutive model}
\end{subfigure}
\hspace{0.7cm}
\begin{subfigure}[t]{0.24\textwidth}
\centering
\includegraphics[width=0.9\textwidth]{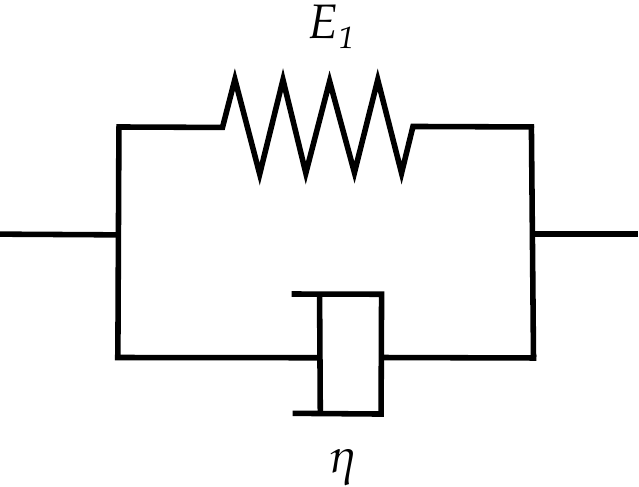}
\caption{Kelvin-Voigt constitutive model}
\end{subfigure}
\hspace{0.7cm}
\begin{subfigure}[t]{0.33\textwidth}
\centering
\includegraphics[width=\textwidth]{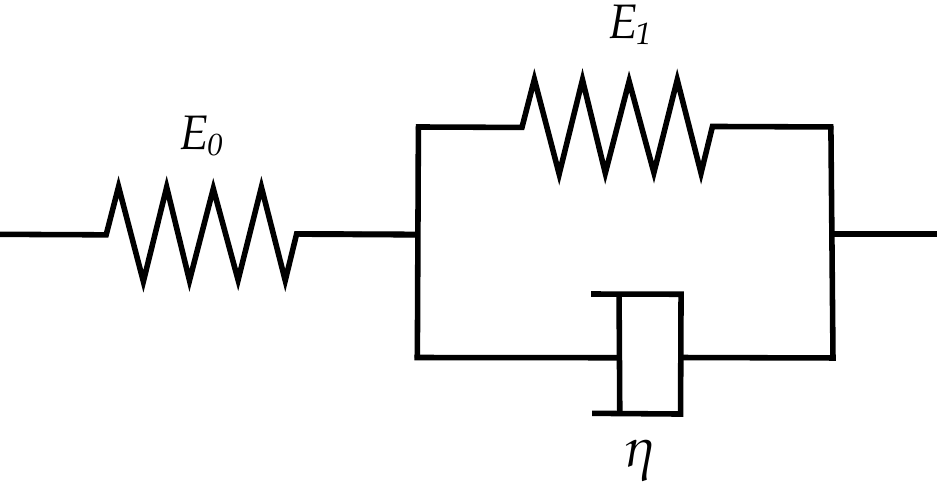}
\caption{Standard Linear Solid constitutive model}
\end{subfigure}
\caption{Scheme of the three simplest linear viscoelastic constitutive models. Coefficients $E$ represent Young moduli associated to springs, while $\eta$ identifies the viscosity coefficient characterizing the dash-pot.}
\label{fig.MXeKVeSLSM}
\end{figure}

\subsubsection{Maxwell constitutive law}
The Maxwell (MX) constitutive model consists on a spring and a dash-pot in series \cite{Lakes2009}, as presented in  Fig. \ref{fig.MXeKVeSLSM} (a). In addition to the contribution of the linear elastic solid, we need to recall that the behavior of a linear dash-pot follows the mechanics of a piston moving in an ideal incompressible viscous (Newtonian) fluid. For this dash-pot, the rheological law reads \cite{gurtin1962}:
\begin{equation}
\label{eq:dashpot}
\sigma = \eta \frac{\d \epsilon}{\d t} ,
\end{equation}
which means, the larger the stress, the faster the material deforms. Therefore, the constitutive law of the MX model results:
\begin{equation}
\label{eq.constitutiveMX}
\frac{\mathrm{d} \sigma}{\mathrm{d} t} = E_0 \frac{\mathrm{d} \epsilon}{\mathrm{d} t} - \frac{1}{\tau_r}\sigma ,
\end{equation}
where $\tau_r$ is the \emph{relaxation time} of the material, in this model defined as $\tau_r = \eta/E_0$.

The response of the model to a sudden load, maintained constant in time, reflects an instantaneous deformation of the spring, being the instantaneous Young modulus $E_0$ representative of the elastic response of the material (in our case, the vessel wall), when the viscous components are ``shortened'' and do not deform. Subsequently, the additional strain of the dash-pot, which takes time to react, manifests itself. Once the load is removed, the spring reacts again immediately, relaxing, but the dash-pot does not have any tendency to recover. Thus, the system remains with a ``creep'' strain due to the dash-pot, even though not very realistic resulting in a straight line in time, in contrast to curves that are observed experimentally \cite{Lakes2009}. Otherwise, if a step strain is applied to the unit, the relaxation response results:
\begin{equation}
E(t) = E_0 e^{-\frac{t}{\tau_r}} .
\label{relaxationfunction_MX}
\end{equation}
This equation, called \emph{relaxation function}, describes how the stiffness and the elastic behavior of the material change in time with respect to the value of the relaxation time $\tau_r$, starting from the instantaneous value of the Young modulus $E_0$. Let us notice, indeed, that $t\to0 \implies E(t) \to E_0$. On the other hand, $t\to\infty \implies E(t)\to0$, so with an elastic modulus that vanishes asymptotically in time. This last result implies that the so-called asymptotic Young modulus (representative of the elastic response when all effects of the viscosity of the viscoelastic material have manifested themselves) of the Maxwell rheological model is null, i.e., $E_{\infty}= 0$.

Taking into account a generic vessel (artery or vein), differentiating with respect to time eq.~\eqref{epsilon-A}, it follows that
\begin{equation}
\label{eq.depsilon}
\frac{\d\varepsilon}{\d t} = \frac1{A} \left(m \alpha^m -n \alpha^n\right) \frac{\d A}{\d t} ;
\end{equation}
while differentiating eq.~\eqref{Barlow} we obtain
\begin{equation}
\label{eq.dsigma}
\frac{\d \sigma}{\d t} = W \frac{\d p}{\d t} .
\end{equation}
With these expressions and using the continuity eq.~\eqref{eq.contST}, it is possible to rewrite the MX constitutive equation \eqref{eq.constitutiveMX} into the following PDE:
\begin{equation}
\label{eq.MX}
\frac{\partial p}{\partial t} + \frac{E_0}{WA} \left(m \alpha^m - n \alpha^n\right) \frac{\partial (Au)}{\partial x} = - \frac{1}{\tau_r}\left(p-p_{0}\right) .
\end{equation}

\subsubsection{Kelvin--Voigt constitutive law}
If we consider a single spring and a dash-pot connected in parallel, so that they both experience the same deformation or strain and the total stress is the sum of the stresses in each element, we have what is known as Kelvin--Voigt (KV) constitutive model \cite{Lakes2009}, represented in Fig. \ref{fig.MXeKVeSLSM} (b). A KV unit relates strain and stress as follows:
\begin{equation}
\sigma = E_1\varepsilon + \eta \frac{\d \varepsilon}{\d t} ,
\label{eq.constitutiveKV}
\end{equation}
It is worth to notice that, due to the mechanics of the model, the asymptotic Young modulus of the KV unit corresponds to the value $E_{\infty}=E_1$.

As previously applied to the MX constitutive law, using eqs.~\eqref{hooke_simple}--\eqref{Barlow}--\eqref{eq.dsigma} and the continuity eq.~\eqref{eq.contST}, it is possible to obtain the following PDE corresponding to the KV viscoelastic law:
\begin{equation}
p = p_0 + \frac{E_{\infty}}{W}\left(\alpha^m - \alpha^n\right) - \frac{\eta}{W A}\left(m\alpha^m - n\alpha^n\right)\frac{\partial(Au)}{\partial x}  .
\label{eq.KV}
\end{equation}

This viscoelastic law, which is widely adopted among literature's well recognized blood flow models, identifying with the parameter
\begin{equation*}
\Gamma = \frac{\eta h_0 \sqrt{\pi}}{2}
\label{gamma}
\end{equation*}
the viscous contribution of the material (in the case of arteries) \cite{alastruey2011,montecinos2014,mynard2015}, has the weakness of defining a relaxation response that is a constant plus a Dirac delta function. When the KV unit is placed at constant strain, indeed, the constitutive eq.~\eqref{eq.constitutiveKV} reduces to the simple Hooke's law, eq.~\eqref{hooke_simple}. In this way, the stress is taken up by the spring and is constant too. In fact, we observe that there is no stress relaxation over time and, therefore, it is not possible to define any relaxation function of the material \cite{Lakes2009}. On the other hand, when the KV unit is suddenly loaded with a constant stress over time, the spring cannot immediately deform because is held back by the dash-pot. Hence, this model is not able to describe an instantaneous elastic strain, being the stress initially totally absorbed by the dash-pot and transferred to the spring only successively in time. Also when unloading the unit, the dash-pot does not permit to the spring to instantaneously relax and no permanent strain is left.
\subsubsection{Standard Linear Solid constitutive law}
A richer behavior can be modeled by considering the Standard Linear Solid (SLS) constitutive model, represented in Fig. \ref{fig.MXeKVeSLSM} (c) in its version with a Kelvin--Voigt unit in series with an additional elastic spring \cite{Lakes2009}. Note that the same considerations that will follow also apply to the SLS law in its version with a Maxwell unit in parallel with an additional spring, since there is an exact correspondence between the parameters of the two versions, as already discussed in \cite{bertaglia2018}.
The constitutive equation of the SLS model reads
\begin{equation}
\label{constitutiveEq.}
	\frac{\d \sigma}{\d t} = E_0 \frac{\d \varepsilon}{\d t} - \frac{1}{\tau_r}(\sigma - E_\infty \varepsilon),
\end{equation}
where we have the instantaneous Young modulus \(E_0\), the asymptotic Young modulus \(E_\infty\) and the relaxation time \(\tau_r\), the last two defined respectively as
\begin{equation}
\label{SLSparameters}
E_\infty = \frac{E_0E_1}{E_0 + E_1}, \qquad \tau_r = \frac{\eta}{E_0 + E_1} = \frac{(E_0 - E_{\infty})\eta}{E_0^2},
\end{equation}
with $E_0$ being the Young modulus of the additional spring, in series with the KV unit, and $E_1$ Young modulus of the elastic spring of the KV element itself, as shown in Fig. \ref{fig.MXeKVeSLSM} (c). 

This model is the simplest linear viscoelastic model able to realistically exhibit all the three primary features of a viscoelastic material: creep, stress relaxation and hysteresis \cite{bertaglia2020,bertaglia2018}.
When the system is suddenly loaded, the instantaneous response is attributed solely to the first spring. The dash-pot then takes up the stress, transferring the load to the second spring as it slowly opens over time. If the load is maintained constant in time, the two springs collaborate as if there were only the two of them in series: $1/E_{\infty} = 1/E_{0} + 1/E_{1}$. While reaching this asymptotic state, the creep is attributed only to the spring in parallel with the dash-pot. Finally, when unloading the system, the first spring relaxes immediately while the second reacts slowly, being held back by the dash-pot.

In addition, the SLS law permits to define a relaxation function, describing how the stiffness of the material changes in time, starting from the instantaneous value and reaching the asymptotic one:
\begin{equation}
E(t) = E_0 e^{-\frac{t}{\tau_r}} + E_{\infty} \left(1- e^{-\frac{t}{\tau_r}}\right) .
\label{relaxationfunction}
\end{equation}

As previously presented for the MX and the KV constitutive laws, it is possible to write also the SLS constitutive equation in terms of pressure and area through a PDE \cite{bertaglia2018}.  
Introducing eqs.~\eqref{hooke_simple}--\eqref{Barlow}--\eqref{eq.depsilon}--\eqref{eq.dsigma} in the rheological law \eqref{constitutiveEq.} and using the continuity eq. \eqref{eq.contST}, the sought PDE is obtained:
\begin{equation}
\frac{\partial p}{\partial t} + \frac{E_0}{WA} \left(m \alpha^m - n \alpha^n\right) \frac{\partial (Au)}{\partial x} = -\frac{1}{\tau_r}\left[p-p_{0} - \frac{E_{\infty}}{W}\left( \alpha^m - \alpha^n\right) \right] .
\label{viscoPDEtubelaw}
\end{equation}
In the above equation, the coefficient of the transport term
\begin{equation}
E_0 G(A)\,, \qquad G(A)=\frac{1}{WA} \left(m \alpha^m - n \alpha^n\right)
\label{d_coeff}
\end{equation}
identifies the elastic contributions of the mechanics of the material, while the source term 
\begin{equation}
 -\frac{1}{\tau_r}\left(p-F(A) \right)\,, \qquad F(A) = p_{0} + \frac{E_{\infty}}{W}\left( \alpha^m - \alpha^n\right)
\label{source}
\end{equation}
takes into account the viscous property of the vessel wall.


\section{Asymptotic limits}
\label{asymptotics}
If we consider the SLS constitutive law \eqref{viscoPDEtubelaw} as closing equation for the governing system \eqref{eq.cont&mom}, we obtain an \textit{augmented} fluid-structure interaction (FSI) system of the cardiovascular bio-fluid dynamics, which reads \cite{bertaglia2020,bertaglia2020a,piccioli2021}:
\begin{subequations}
\begin{align}
	&\frac{\partial A}{\partial t} + \frac{\partial(Au)}{\partial x} = 0 \label{eq.cont}\\
	&\frac{\partial(Au)}{\partial t}+ \frac{\partial(Au^2)}{\partial x}  + \frac{A}{\rho} \, \frac{\partial p}{\partial x} = 0 		 \label{eq.mom}\\
	&\frac{\partial p}{\partial t} + E_0 G(A) \,\frac{\partial(Au)}{\partial x} = -\frac{1}{\tau_r}\left( p - F(A)\right) ,\label{eq.PDE}
\end{align}
\label{completesyst}
\end{subequations}
with $G(A)$ and $F(A)$ defined in eqs.~\eqref{d_coeff}-\eqref{source}, respectively. Thus, from now on, we will refer to the \textit{augmented blood flow model} meaning that it is the blood flow model in which the viscoelastic SLS constitutive law is used to close the problem.

This model is hyperbolic, being the Jacobian matrix diagonalizable, with a diagonal matrix $\Lambda$ containing all real eigenvalues and a complete set of linearly independent eigenvectors represented by the columns of the matrix $R$:\\
	\[\Lambda =
\begin{pmatrix} 
  	u-c &0 &0 \\ 0 &0 &0 \\ 0 &0 &u+c 
\end{pmatrix}, \quad
	R =
\begin{pmatrix} 
  	1 &1 &1 \\ u-c &0 &u+c \\ E_0 G(A) &\frac{\rho u^2}{A} &E_0 G(A)
\end{pmatrix}, \]
where $c$ is the wave speed,
\begin{equation}
c = \sqrt{\frac{A E_0}{\rho} \,  G(A) } .
\label{eq:cel}
\end{equation} 
The second eigenvector of the system is associated with a linearly degenerate (LD) characteristic field, while the first and the third define genuinely non-linear fields (leading to the formation of shocks or rarefaction waves) \cite{piccioli2021}. Concerning the Riemann Invariants of the system, those associated with the LD field are \cite{bertaglia2020,piccioli2021}
\begin{equation}
\Gamma_1^{LD} = Au , \qquad \Gamma_2^{LD} = p +\frac{1}{2}\rho u^2 ,
\label{eq.RI-LD}
\end{equation}
defining the quantities that remain constant across contact discontinuities. The Riemann Invariants associated to the genuinely non-linear fields are, instead,
\begin{subequations}
\begin{align}
&\Gamma_1 = u - \int\frac{c}{A}\,\d A, \\
&\Gamma_2 = u + \int\frac{c}{A}\,\d A, \\
&\Gamma_3 = p - \int E_0 G(A) \,\d A = p - \frac{E_0}{W} \left(\alpha^m - \alpha^n\right) .
\end{align}
\label{eq.RI}
\end{subequations}
Notice that when dealing with arteries, also integrals $\Gamma_1$ and $\Gamma_2$ can be analytically solved, resulting
\[\Gamma_{1,2} = u \mp 4c .\]

Finally, it is here remarked that to accommodate a correct numerical treatment of possible longitudinal discontinuities in space of geometrical and mechanical properties, such as equilibrium cross-sectional area, instantaneous Young modulus \(E_0\), asymptotic Young modulus $E_{\infty}$, viscosity coefficient $\eta$, and equilibrium pressure \(p_{0}\), it is necessary to introduce additional equations to system \eqref{completesyst} \cite{bertaglia2020,muller2013}. Considering these variables constant in time, the additional equations result: \(\partial_t A_0=0\), \(\partial_t E_0=0\), \(\partial_t E_{\infty}=0\), \(\partial_t \eta=0\) and \(\partial_t p_{0}=0\). 


Let us know analyze the asymptotic limits of system \eqref{completesyst} as the scaling parameter $\tau_r \to 0$, i.e., the so-called \emph{zero-relaxation} limits. We will prove that by choosing the Standard Linear Solid constitutive law as closing equation for the system, with an appropriate choice of the scaling parameters, all the different rheological characterizations previously discussed can be recovered. Hence, the proposed model can account for several mechanical behaviors of the vessel wall, from the elastic to different viscoelastic ones.

\subsection{Hyperbolic scaling}
If $\tau_r \to 0$ while $\eta \to 0$, from the relaxation function \eqref{relaxationfunction} we observe that 
$E(t) \to E_{\infty}=\hat{E}$, 
thus the stiffness of the material remains constant in time. This implies that, in this limit, the vessel wall tends to behave as a purely elastic material characterized by the Young modulus $\hat{E}=\frac{E_0E_1}{E_0+E_1}$, which is the resultant of the sum of the Young moduli of the two springs in series. From eq.~\eqref{eq.PDE} we then recover exactly the elastic constitutive law \eqref{elastictubelaw}:
\begin{equation}
p = F(A) = p_{0} + \frac{E_{\infty}}{W} \left(\alpha^m - \alpha^n\right) .
\label{elastictubelaw1}
\end{equation}
Using this equation into eq.~\eqref{eq.mom}, we observe that the proposed model recovers the classical blood flow elastic model \cite{formaggia2003,muller2013}:
\begin{subequations}
\begin{align}
	&\frac{\partial A}{\partial t} + \frac{\partial(Au)}{\partial x} = 0 \\
	&\frac{\partial(Au)}{\partial t}+ \frac{\partial(Au^2)}{\partial x}  + \frac{A}{\rho} \, \frac{\partial F(A) }{\partial x} = 0 .
\end{align}
\label{completesyst_el}
\end{subequations}

Notice that we recover the same equilibrium also for $\tau_r \to 0$ while $E_1 \to \infty$, with the only difference that, in this case, the material tends to behave like a spring with Young modulus $\hat{E}=E_0$.

\subsection{Diffusive scaling} 
If $\tau_r \to 0$ and $E_0 \to \infty$ while $\eta = \tau_r(E_0 + E_1) \to \tau_r E_0$ remains finite, we observe that, from eq.~\eqref{eq.PDE} we recover the diffusive behavior of the solution described by the Kelvin--Voigt constitutive law, hence eq.~\eqref{eq.KV}:
\[ p = F(A) - \eta G(A) \frac{\partial(Au)}{\partial x} = p_0 + \frac{E_{\infty}}{W}\left(\alpha^m - \alpha^n\right) - \frac{\eta}{W A}\left(m\alpha^m - n\alpha^n\right)\frac{\partial(Au)}{\partial x} ,\]
where $E_{\infty}=E_1.$
Inserting the above equilibrium into eq. \eqref{eq.mom}, we obtain the following parabolic, diffusive model:
\begin{subequations}
\begin{align}
	&\frac{\partial A}{\partial t} + \frac{\partial(Au)}{\partial x} = 0 \\
	&\frac{\partial(Au)}{\partial t}+ \frac{\partial(Au^2)}{\partial x}  + \frac{A}{\rho} \, \frac{\partial F(A)}{\partial x} = \frac{A}{\rho}\frac{\partial}{\partial x}\left(\eta\, G(A) \frac{\partial(Au)}{\partial x}\right) ,
\end{align}
\label{completesyst_KV}
\end{subequations}
where the presence of the additional parabolic term in the momentum equation is evident when comparing it to system \eqref{completesyst_el}. In fact, this asymptotic limit describes a strongly diffusive dynamic, in which viscous effects of the vessel wall occur over long rescaling times.

We highlight that system \eqref{completesyst_KV} reads exactly as literature blood flow models that adopts the Kelvin--Voigt rheological law \cite{alastruey2011,formaggia2003,montecinos2014,mynard2015}. 
It is also worth to emphasize that, in this limit, the celerity of the system $c \to \infty$ because of its dependence on the Young's modulus $E_0$ (see eq. \eqref{eq:cel} for the definition), which is indeed in agreement with the parabolic scaling obtained.

\subsection{Perturbation of the elastic local equilibrium}
Let us now analyze a first order perturbation in $\tau_r$ of the first local equilibrium, i.e., the elastic hyperbolic scaling, which has been demonstrated to be the leading order equation of the SLS model as $\tau_r~\to~0$ while $\eta \to 0$. The first order perturbation reads
\begin{equation}
p = p_0 + \frac{E_{\infty}}{W}\left(\alpha^m - \alpha^n\right) + \tau_r p_1 + \mathcal{O}\left(\tau_r^2\right) .
\label{SLS_expansion1}
\end{equation}
Substituting this expansion in eq.~\eqref{viscoPDEtubelaw}, we obtain
\[ p_1 = -\frac{\partial p}{\partial t} - \frac{E_0}{WA} \left(m \alpha^m - n \alpha^n\right) \frac{\partial (Au)}{\partial x} + \mathcal{O}\left(\tau_r^2\right).\]
Since, deriving with respect to time eq.~\eqref{elastictubelaw1},
\[ \frac{\partial p}{\partial t} = - \frac{E_{\infty}}{WA} \left(m \alpha^m - n \alpha^n\right) \frac{\partial (Au)}{\partial x},\]
we finally have
\[ p_1 = \frac{E_{\infty} - E_0}{WA}\left(m \alpha^m - n \alpha^n\right) \frac{\partial (Au)}{\partial x} + \mathcal{O}\left(\tau_r^2\right).\]
Substituting back in eq.~\eqref{SLS_expansion1} leads to
\[
p = p_0 + \frac{E_{\infty}}{W}\left(\alpha^m - \alpha^n\right) -  \frac{\tau_r\left(E_0 - E_{\infty}\right)}{WA}\left(m \alpha^m - n \alpha^n\right) \frac{\partial (Au)}{\partial x} + \mathcal{O}\left(\tau_r^2\right) .
\]
Considering the last equivalence in eq.~\eqref{SLSparameters} and omitting the second order terms in $\tau_r$, we can rewrite the above equation as follows, defining a new viscoelastic model that is a second-order accurate approximation, for $\tau_r \ll 1$ and $\eta \ll 1$ (namely, for mild viscous effects), of the Standard Linear Solid constitutive law:
\begin{equation}
p = p_0 + \frac{E_{\infty}}{W}\left(\alpha^m - \alpha^n\right) -  \frac{\left(E_0 - E_{\infty}\right)^2}{E_0^2}\frac{\eta}{WA}\left(m \alpha^m - n \alpha^n\right) \frac{\partial (Au)}{\partial x} .
\label{SLS_expansion1_final}
\end{equation}
Inserting the above law into eq. \eqref{eq.mom} leads to the following parabolic model:
\begin{subequations}
\begin{align}
	&\frac{\partial A}{\partial t} + \frac{\partial(Au)}{\partial x} = 0 \\
	&\frac{\partial(Au)}{\partial t}+ \frac{\partial(Au^2)}{\partial x}  + \frac{A}{\rho} \, \frac{\partial F(A)}{\partial x} = \frac{A}{\rho}\frac{\partial}{\partial x}\left(\frac{\left(E_0 - E_{\infty}\right)^2}{E_0^2}\eta\, G(A) \frac{\partial(Au)}{\partial x}\right) ,
\end{align}
\label{completesyst_newVisco}
\end{subequations}

Notice that this new viscoelastic constitutive law differs from the Kelvin--Voigt one because of the presence of a correction factor $\left(E_0 - E_{\infty}\right)^2/E_0^2$ in the viscosity coefficient of the material (which is simply $\eta$ in the KV unit). However, unlike the Kelvin--Voigt constitutive law, this new model describes a dynamic in which the viscous effects of the viscoelastic wall are not predominant, being indeed $\eta \ll 1$ and $0~<~\left(E_0 - E_{\infty}\right)^2/E_0^2~<~1$. In fact, the newly proposed viscoelastic model permits to capture the second order small viscosity effects, being a second-order accurate approximation (for small relaxation times $\tau_r$ and viscosity coefficient $\eta$) of the model closed with the SLS constitutive law.

Using eqs.~\eqref{hooke_simple}--\eqref{Barlow}--\eqref{eq.depsilon}, the rheological law can be written in terms of stress and strain, resulting
\[ \sigma = E_{\infty}\varepsilon + \frac{(E_0 - E_{\infty})^2}{E_0^2} \eta \frac{\d \varepsilon}{\d t}.\]

In addition, with the above equation, we observe that the elastic law itself is a good approximation of the SLS to the first order perturbation in $\tau_r$, provided that $E_0~\approx~E_{\infty} \implies p_1\approx0$. 
This mathematical result actually confirms what the mechanics of the SLS already states: if the Young modulus is constant in time, the material cannot express a relaxation process of the stress and cannot dissipate energy, so it behaves like a simple elastic spring.

\begin{figure}[!tbp]
\centering
\includegraphics[width=0.9\textwidth]{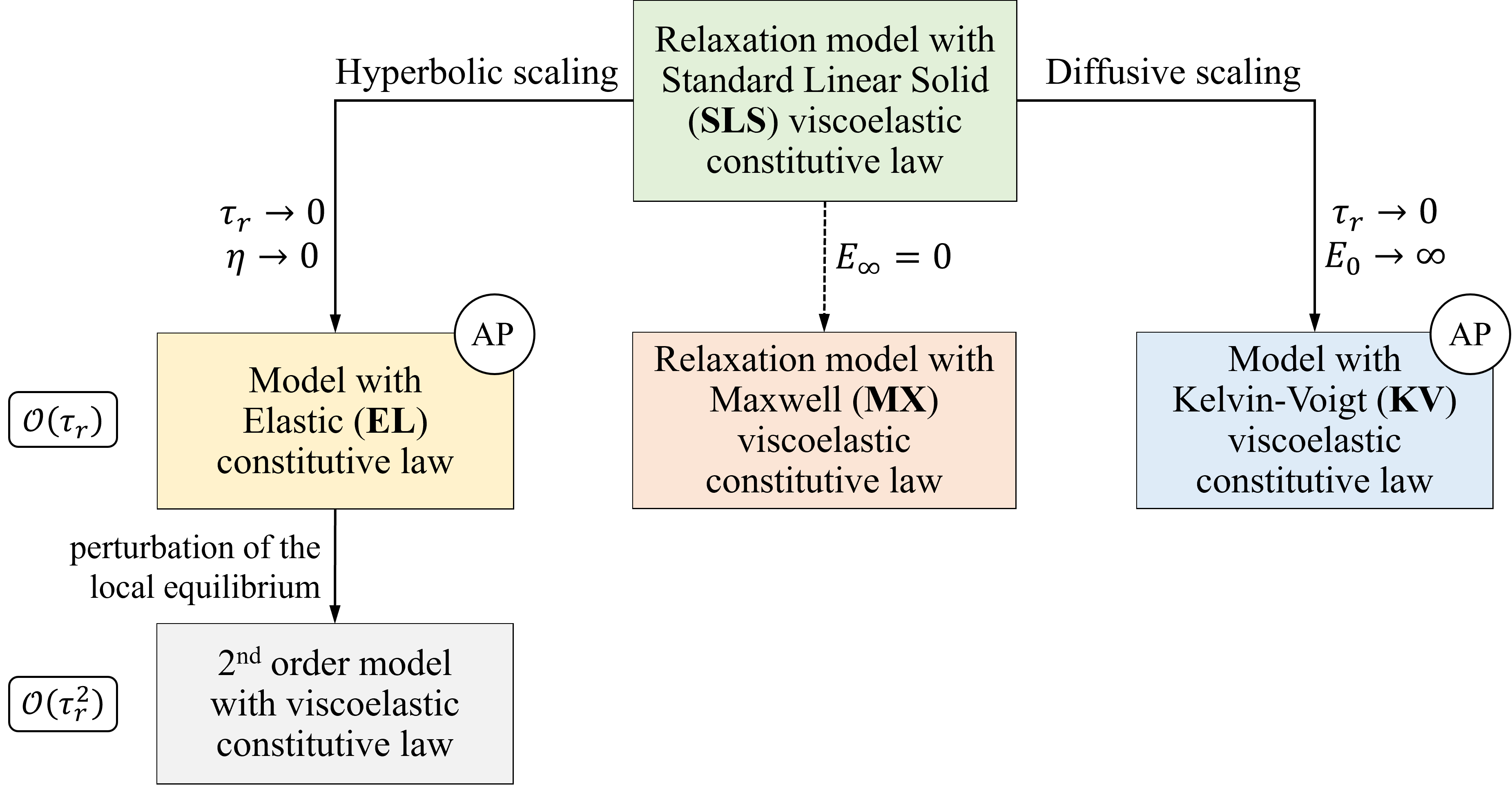}
\caption{Asymptotic limits of the multiscale constitutive framework. Using the SLS rheological law as closing equation for the blood flow model, it is possible to recover, under a suitable scaling, both the elastic EL and viscoelastic KV constitutive behaviors, leading to hyperbolic and diffusive behaviors, respectively. The need for an asymptotic-preserving (AP) numerical method to computationally simulate these limits is highlighted. A specific choice of the asymptotic Young modulus leads to the viscoelastic MX law. In addition, a new second-order viscoelastic model, for $\tau_r \ll 1$, can be derived through the perturbation of the local elastic equilibrium.}
\label{fig.APlimits}
\end{figure}

\begin{remark} 
In classical kinetic theory, the space-time scaling just discussed is related to the hydrodynamical limits of the Boltzmann equation \cite{albi2019,cercignani1994,lions1997}. In particular, the hyperbolic scaling corresponds to the compressible Euler scaling. In the case of the diffusive scaling, the dissipation effects become non-negligible and we get the incompressible Navier--Stokes scaling. The second-order accurate viscoelastic model here derived, instead, leads to a model linked to the compressible Navier--Stokes equations.
\end{remark}

\begin{remark}
Regarding the Maxwell viscoelastic characterization, let us point out that the Standard Linear Solid constitutive law exactly matches the Maxwell one when $E_1 = E_{\infty} = 0$. In fact, this implies that $F(A)=p_0$ and $\tau_r=\eta/E_0$, and from eq.~\eqref{eq.PDE} we recover eq.~\eqref{eq.MX}. Thus, the complete augmented model reads
\begin{subequations}
\begin{align}
	&\frac{\partial A}{\partial t} + \frac{\partial(Au)}{\partial x} = 0 \\
	&\frac{\partial(Au)}{\partial t}+ \frac{\partial(Au^2)}{\partial x}  + \frac{A}{\rho} \, \frac{\partial p}{\partial x} = 0 \\
	&\frac{\partial p}{\partial t} + E_0 G(A) \frac{\partial (Au)}{\partial x} = - \frac{1}{\tau_r}\left(p-p_{0}\right) .
\end{align}
\end{subequations}
In addition, the relaxation function \eqref{relaxationfunction} becomes equivalent to eq.~\eqref{relaxationfunction_MX}. 
\end{remark}

The summary of all the above limits and derivations is shown in the diagram in Fig. \ref{fig.APlimits} from a constitutive framework perspective.


\section{Asymptotic-preserving method}
\label{section_numericalmethod}
As pointed out in the previous section, the augmented blood flow model \eqref{completesyst} turns out to be a multiscale hyperbolic system, which, accordingly to the choice of the scaling parameters, can describe also diffusive-type phenomena, i.e., those associated with the Kelvin--Voigt characterization of the vessel wall material. Thus, selecting an appropriate numerical method is crucial to ensure the numerical discretization is consistent with all asymptotic behaviors of the model while maintaining the expected accuracy in the stiff limit. Moreover, we aim at working with a numerical scheme that does not have restrictions of the time step size related to the smallness of the scaling parameters. Indeed, generally, the latter could lead the time step size to be prohibitively small for stability reasons as $\tau_r \to 0$. These features are fulfilled if the chosen numerical method respects the \emph{asymptotic-preserving} (AP) property \cite{pareschi2005,albi2019,boscarino2017}. In the light of this, to solve the problem we consider a third order AP Implicit-Explicit (IMEX) Runge-Kutta Finite Volume method, following the partitioned approach proposed in \cite{boscarino2017a} for hyperbolic systems with multiscale relaxation.

\subsection{An IMEX Runge-Kutta Finite Volume scheme}
IMEX Runge-Kutta schemes can be easily represented by a double tableau (explicit on the left, implicit on the right) in the usual Butcher notation \cite{boscarino2017,pareschi2005}
\begin{center}
\begin{tabular}{c | c}
$\boldsymbol{\tilde c}$ & $\tilde{\mathcal{A}}$ \\  
\hline \\[-1.0em] 
 &  $\boldsymbol{\tilde b}^T$ \\ 
\end{tabular}
\hspace{2.0cm}
\begin{tabular}{c | c}
$\boldsymbol{c}$ & $\mathcal{A}$ \\  
\hline \\[-1.0em] 
 &  $\boldsymbol{b}^T$ \\ 
\end{tabular}\,.
\end{center}
Matrices $\tilde{\mathcal{A}} = (\tilde a_{kj})$, with $\tilde a_{kj} = 0 $ for $ j\geq k$, and $\mathcal{A} = (a_{kj})$ are $s \times s$ matrices, with $s$ number of Runge-Kutta stages. Being always preferable in terms of computational efficiency to deal with diagonally implicit Runge-Kutta (DIRK) schemes, we consider $a_{kj} = 0 $ for $ j > k$. The temporal steps coefficient vectors are $\boldsymbol{\tilde c}= (\tilde c_1, ...,\tilde c_s)^T$ and $\boldsymbol{c}= (c_1, ...,c_s)^T$, 
while vectors $\boldsymbol{\tilde b} = (\tilde b_1, ...,\tilde b_s)^T$ and $\boldsymbol{b} = (b_1, ...,b_s)^T$ are the quadrature weights that permit to combine the internal Runge-Kutta stages. 

In particular, we employ the third-order globally stiffly accurate (GSA) IMEX BPR(3,4,3) scheme proposed in \cite{boscarino2017}, which is characterized by 4 stages for the implicit part and 3 stages for the explicit part (see Appendix \ref{section_appendix}). In the following, we give recall the definition of the GSA property.
\begin{definition}
An IMEX Runge-Kutta method is said to be globally stiffly accurate (GSA) if the corresponding diagonally implicit Runge-Kutta (DIRK) method is stiffly accurate, namely
\begin{equation*}
a_{sj} = b_j, \qquad j = 1,\ldots,s ,
\end{equation*}
and the explicit method satisfies
\begin{equation*}
\tilde a_{sj} = \tilde b_j, \qquad j = 1,\ldots,s-1 .
\end{equation*}
\end{definition}
As a consequence of the above definition the numerical solution of a GSA IMEX method coincides exactly with the last internal stage of the scheme.

Notice that the third order IMEX Runge-Kutta method \eqref{eq:tableaux} is of type CK (see \cite{boscarino2017a,dimarco2013}), since the matrix $\mathcal{A}$ can be written as
\begin{equation}
\mathcal{A} =  \begin{pmatrix} 
0 &\boldsymbol{0}\\ \boldsymbol{a} &\hat{\mathcal{A}}
\end{pmatrix},
\label{A_im}
\end{equation}
with $\boldsymbol{a} = (a_{21}, ...,a_{s1})^T \in \mathbb{R}^{(s-1)}$ and the sub-matrix $\hat{\mathcal{A}}\in \mathbb{R}^{(s-1) \times (s-1)}$ is invertible, hence $a_{kk} \ne 0,\, k=2,\ldots,s$. We will also make use of the following representation of the matrix $\tilde{\mathcal{A}}$ in the explicit Runge-Kutta method:
\begin{equation}
\tilde{\mathcal{A}} =  \begin{pmatrix} 
0 &\boldsymbol{0}\\ \tilde{\boldsymbol{a}} &\hat{\tilde{\mathcal{A}}}
\end{pmatrix},
\label{A_ex}
\end{equation}
where $\tilde{\boldsymbol{a}} = (\tilde{a}_{21}, ...,\tilde{a}_{s1})^T \in \mathbb{R}^{(s-1)}$ and $\hat{\tilde{\mathcal{A}}}\in \mathbb{R}^{(s-1) \times (s-1)}$.

To obtain an AP scheme, the IMEX Runge-Kutta temporal discretization of system \eqref{completesyst}, written in semi-discrete form, consists in computing the internal stages
\begin{subequations}
\begin{align}
	&A^{(k)} = A^n -  \Delta t \sum_{j=1}^{k-1} \tilde a_{kj}\, \partial_x (Au)^{(j)}
	\\
	&(Au)^{(k)} = (Au)^n -  \Delta t \sum_{j=1}^{k-1} \tilde{a}_{kj} \partial_x \frac{\left[(Au)^{(j)}\right]^2}{A^{(j)}} - \Delta t \sum_{j=1}^{k-1} \tilde{a}_{kj} \frac{A^{(j)}}{\rho} \partial_x p^{(j)}
	\\
	&p^{(k)} = p^n - \Delta t \sum_{j=1}^{k} a_{kj} E_0 G(A^{(j)}) \partial_x(Au)^{(j)} - \Delta t \sum_{j=1}^{k} a_{kj} \frac{1}{\tau_r}\left( p^{(j)} - F(A^{(j)})\right),
\end{align}
\label{eq.iterIMEX}
\end{subequations}
for $k=1,\ldots,s$, and then the final numerical solution (even though, for definition, the numerical solution of a GSA IMEX Runge-Kutta scheme coincides exactly with the last internal stage of the scheme), which reads
\begin{subequations}
\begin{align}	
	&A^{n+1} = A^n -  \Delta t \sum_{k=1}^{s} \tilde b_{k}\, \partial_x (Au)^{(k)}
	\\
	&(Au)^{n+1} = (Au)^n -  \Delta t \sum_{k=1}^{s} \tilde{b}_{k} \partial_x \frac{\left[(Au)^{(k)}\right]^2}{A^{(k)}} - \Delta t \sum_{k=1}^{s} \tilde{b}_{k} \frac{A^{(k)}}{\rho} \partial_x p^{(k)}
	\\
	&p^{n+1} = p^n -  \Delta t \sum_{k=1}^{s} b_{k} E_0 G(A^{(k)})\partial_x(Au)^{(k)}  - \Delta t \sum_{k=1}^{s} b_{k} \frac{1}{\tau_r}\left( p^{(k)} - F(A^{(k)})\right).
\end{align}
\label{eq.finalIMEX}
\end{subequations}
The time step size \(\Delta t = t^{n+1}-t^{n}\) follows the less restrictive between the standard hyperbolic Courant-Friedrichs-Levy condition, $\Delta t \leq \CFL \Delta x/\max\left\{|u\pm c|\right\}$, and the parabolic stability restriction, $\Delta t \leq \nu \Delta x^2$, where $\Delta x  = x_{i+1/2} - x_{i-1/2}$, $i=1,\ldots,N_x$, with $N_x$ number of cells, is the size of the uniform space grid and $\CFL$ and $\nu$ are suitable stability constants \cite{boscarino2017}.

For the space derivatives, we consider a Finite Volume discretization.
To ensure the correct treatment of both conservative and non-conservative terms of system \eqref{completesyst} even in case of variables that are discontinuous in space, to evaluate numerical fluxes and non-conservative jump terms we employ the Dumbser-Osher-Toro (DOT) Riemann solver \cite{dumbser2011,dumbser2011a,bertaglia2020,bertaglia2018}. Boundary-extrapolated values at both interfaces of cell $i$ are computed through a third order Weighted Essentially Non-Oscillatory (WENO) reconstruction \cite{shu1998,qiu2002}.

\subsection{AP property}
To write the IMEX scheme in compact form, let us denote 
\[\boldsymbol{A}~=~\left(A^{(1)}, \ldots, A^{(s)}\right)^T,\quad 
\boldsymbol{Au} = \left((Au)^{(1)}, \ldots, (Au)^{(s)}\right)^T,\quad
\boldsymbol{p} = \left(p^{(1)}, \ldots, p^{(s)}\right)^T,\]
\[\boldsymbol{Au^2} = \left([(Au)^{(1)}]^2/A^{(1)}, \ldots, [(Au)^{(s)}]^2/A^{(s)}\right)^T,\quad
\boldsymbol{e}~=~\left(1, \ldots, 1\right)^T \in \mathbb{R}^s,\]
\[\boldsymbol{F}(\boldsymbol{A})~=~\left(F(A^{(1)}), \ldots, F(A^{(s)})\right)^T,\quad
\boldsymbol{G}(\boldsymbol{A})=\left(G(A^{(1)}), \ldots, G(A^{(s)})\right)^T,\]
obtaining
\begin{subequations}
\begin{align}
	\boldsymbol{A} &= A^n \boldsymbol{e} -  \Delta t \,\tilde{\mathcal{A}} \, \partial_x( \boldsymbol{Au}) \label{eq.cont_compact}
	\\
	\boldsymbol{Au} &= (Au)^n \boldsymbol{e} -  \Delta t \,\tilde{\mathcal{A}} \,\partial_x\left(\boldsymbol{Au^2}\right) -\Delta t \, \tilde{\mathcal{A}}\, \frac{\boldsymbol{A}}{\rho} \partial_x \boldsymbol{p} \label{eq.mom_compact}
	\\
	\boldsymbol{p} &= p^n \boldsymbol{e} - \Delta t \,\mathcal{A} \, E_0 \boldsymbol{G}(\boldsymbol{A}) \partial_x(\boldsymbol{Au}) - \frac{\Delta t \,\mathcal{A}}{\tau_r}\left(\boldsymbol{p} - \boldsymbol{F}(\boldsymbol{A})\right) \label{eq.PDE_compact}
\end{align}
\label{eq.iterIMEX_compact}
\end{subequations}
and
\begin{subequations}
\begin{align}	
	A^{n+1} &= A^n -  \Delta t \,\boldsymbol{\tilde b}^{T}\, \partial_x (\boldsymbol{Au})
	\\
	(Au)^{n+1} &= (Au)^n -  \Delta t \,\boldsymbol{\tilde{b}}^{T} \, \partial_x \left(\boldsymbol{Au^2}\right) -  \Delta t \, \boldsymbol{\tilde b}^{T}\, \frac{\boldsymbol{A}}{\rho} \partial_x \boldsymbol{p}
	\\
	p^{n+1} &= p^n -  \Delta t \,\boldsymbol{b}^{T} E_0\boldsymbol{G}(\boldsymbol{A}) \partial_x(\boldsymbol{Au}) - \frac{\Delta t \,\boldsymbol{b}^{T}}{\tau_r}\left(\boldsymbol{p} - \boldsymbol{F}(\boldsymbol{A})\right) .
\end{align}
\label{eq.finalIMEX_compact}
\end{subequations}
Now, recalling eqs. \eqref{A_im} and \eqref{A_ex}, the IMEX scheme \eqref{eq.iterIMEX_compact}-\eqref{eq.finalIMEX_compact} can be written as follows,
\begin{subequations}
\begin{align}
A^{(1)} &= A^n  \label{eq.cont_compact1}
\\
	\boldsymbol{\hat{A}} &= A^n \boldsymbol{\hat{e}} - \Delta t \,\boldsymbol{\tilde{a}} \, \partial_x (Au)^{(1)} -  \Delta t \,\hat{\tilde{\mathcal{A}}} \, \partial_x( \boldsymbol{\hat{Au}})\label{eq.cont_compact2}
	\\
	Au^{(1)} &= (Au)^n \label{eq.mom_compact1}
	\\
	\boldsymbol{\hat{Au}} &= (Au)^n \boldsymbol{\hat{e}} -  \Delta t \,\boldsymbol{\tilde{a}} \,\partial_x(Au^2)^{(1)} -  \Delta t \,\hat{\tilde{\mathcal{A}}} \,\partial_x(\boldsymbol{\hat{Au}^2})  \label{eq.mom_compact2}\\
	&-\Delta t \, \boldsymbol{\tilde{a}}\, \frac{A^{(1)}}{\rho} \partial_x p^{(1)} -\Delta t \, \hat{\tilde{\mathcal{A}}}\, \frac{\boldsymbol{\hat{A}}}{\rho} \partial_x \boldsymbol{\hat{p}}\nonumber
	\\
	p^{(1)} &= p^n \label{eq.PDE_compact1}
	\\
	\boldsymbol{\hat{p}} &= p^n \boldsymbol{\hat{e}} - \Delta t \,\boldsymbol{a} \, E_0 G(A^{(1)}) \partial_x(Au)^{(1)} - \Delta t \,\hat{\mathcal{A}} \, E_0 \boldsymbol{G}(\boldsymbol{\hat{A}}) \partial_x(\boldsymbol{\hat{Au}}) \label{eq.PDE_compact2}\\
	&- \frac{\Delta t \,\boldsymbol{a}}{\tau_r}\left(p^{(1)} - F(A^{(1)})\right) - \frac{\Delta t \,\hat{\mathcal{A}}}{\tau_r}\left(\boldsymbol{\hat{p}} - \boldsymbol{F}(\boldsymbol{\hat{A}})\right) .\nonumber
\end{align}
\label{eq.iterIMEX_compact1}
\end{subequations}
\begin{subequations}
\begin{align}	
	A^{n+1} &= A^n -  \Delta t \,\tilde{b}_1\, \partial_x (Au^{(1)}) -  \Delta t \,\boldsymbol{\hat{\tilde b}}^{T}\, \partial_x (\boldsymbol{\hat{Au}})
	\\
	(Au)^{n+1} &= (Au)^n -  \Delta t \,\tilde{b}_1 \, \partial_x (Au^2)^{(1)} -  \Delta t \,\boldsymbol{\hat{\tilde{b}}}^{T} \, \partial_x (\boldsymbol{\hat{Au}^2}) \\
	&-  \Delta t \, \tilde{b}_1\, \frac{A^{(1)}}{\rho} \partial_x p^{(1)} -  \Delta t \, \boldsymbol{\hat{\tilde{b}}}^{T}\, \frac{\boldsymbol{\hat{A}}}{\rho} \partial_x \boldsymbol{\hat{p}}\nonumber
	\\
	p^{n+1} &= p^n -  \Delta t \, b_1 E_0 G(A^{(1)}) \partial_x (Au)^{(1)} -  \Delta t \,\boldsymbol{\hat{b}}^{T} E_0\boldsymbol{G}(\boldsymbol{\hat{A}}) \partial_x(\boldsymbol{\hat{Au}}) \\
	&- \frac{\Delta t \,b_1}{\tau_r}\left(p^{(1)} - F(A^{(1)})\right) - \frac{\Delta t \,\boldsymbol{\hat{b}}^{T}}{\tau_r}\left(\boldsymbol{\hat{p}} - \boldsymbol{F}(\boldsymbol{\hat{A}})\right) ,\nonumber
\end{align}
\label{eq.finalIMEX_compact1}
\end{subequations}
where $\boldsymbol{\hat{e}}=(1,\ldots,1)^T \in \mathbb{R}^{(s-1)}$.
To solve the problem, from eq.~\eqref{eq.cont_compact2} we can directly compute explicitly $\boldsymbol{\hat{A}}$, thus, also $\boldsymbol{F}(\boldsymbol{\hat{A}})$ and $\boldsymbol{G}(\boldsymbol{\hat{A}})$, while from eq.~\eqref{eq.mom_compact2} we obtain explicitly $\boldsymbol{\hat{Au}}$. Then, inverting eq. \eqref{eq.PDE_compact2}, we obtain an explicit expression to compute the pressure:
\begin{subequations}
\begin{align}
\hat{\boldsymbol{p}} = \left( \frac{\tau_r}{\Delta t} \boldsymbol{I} + \hat{\mathcal{A}}\right)^{-1} &\Big(\frac{\tau_r}{\Delta t}p^n \boldsymbol{\hat{e}} - \tau_r\boldsymbol{a} \,E_0 G(A^{(1)}) \partial_x(Au)^{(1)} - \tau_r\hat{\mathcal{A}} \,E_0\boldsymbol{G}(\boldsymbol{\hat{A}}) \partial_x(\boldsymbol{\hat{Au}}) \label{eq.IMEX_p}\\
&- \boldsymbol{a}\left(p^{(1)} - F(A^{(1)})\right) + \hat{\mathcal{A}}\boldsymbol{F}(\boldsymbol{\hat{A}})\Big)\,, \nonumber
	\end{align}
\end{subequations}
and definitely solve system \eqref{eq.finalIMEX_compact1}.

Before analyzing in details the asymptotic behavior of the method, we need to introduce the notion of \emph{well prepared} initial data, or, equivalently, initial data \emph{consistent} with the limit problem \cite{boscarino2017a}.
\begin{definition}
\label{def.wellIC}
The initial data for system \eqref{completesyst} is said to be \emph{consistent} or \emph{well prepared} if
\begin{equation}
p(x,0) = F(A(x,0)) - \tau_r \,E_0 \, G(A(x,0))\, \partial_x (A(x,0)u(x,0)) + \mathcal{O}(\tau_r).
\label{wellIC}
\end{equation}
\end{definition}

Let us now show the AP property of the IMEX method in the two asymptotic limits.

\begin{proof}[Proof of the AP property for the hyperbolic scaling]
As $\tau_r \to 0$ while $\eta \to 0$, when considering consistent initial data as in \eqref{wellIC}, hence at the initial iteration $p^n = F(A^n)$, from eq.~\eqref{eq.IMEX_p}, recalling also eq. \eqref{eq.cont_compact1} and \eqref{eq.PDE_compact1}, we obtain
\begin{subequations}	
\begin{align*}
\hat{\boldsymbol{p}} &= -\hat{\mathcal{A}}^{-1}\boldsymbol{a}\left(p^{(1)} - F(A^{(1)})\right)+ \hat{\mathcal{A}}^{-1}\hat{\mathcal{A}}\boldsymbol{F}(\boldsymbol{\hat{A}})\,\\
&= -\hat{\mathcal{A}}^{-1}\boldsymbol{a}\left(p^n - F(A^n)\right)+ \boldsymbol{F}(\boldsymbol{\hat{A}})\,\\
&= \boldsymbol{F}(\boldsymbol{\hat{A}})\,.
\end{align*}
\end{subequations}
Moreover, since the scheme is GSA, we also have that $p^{n+1} = F(A^{n+1})$, so at the next time step the initial value remains consistent.
If we now insert this result in eq.~\eqref{eq.mom_compact2}, we exactly recover a consistent explicit Runge-Kutta discretization of the elastic system \eqref{completesyst_el}, having internal stages
\begin{subequations}
\begin{align}
	\boldsymbol{A} &= A^n \boldsymbol{e} -  \Delta t \,\tilde{\mathcal{A}}  \, \partial_x( \boldsymbol{Au})
	\\
	\boldsymbol{Au} &= (Au)^n \boldsymbol{e} -  \Delta t \,\tilde{\mathcal{A}} \,\partial_x\left(\boldsymbol{Au^2}\right) -\Delta t \, \tilde{\mathcal{A}}\, \frac{\boldsymbol{A}}{\rho}\, \partial_x \boldsymbol{F}(\boldsymbol{A})\,,
\end{align}
\label{eq.iterIMEX_compact_el}
\end{subequations}
and final update
\begin{subequations}
\begin{align}
	A^{n+1} &= A^n  -  \Delta t \,\boldsymbol{\tilde b }^T  \, \partial_x( \boldsymbol{Au})
	\\
	(Au)^{n+1} &= (Au)^n  -  \Delta t \,\boldsymbol{\tilde b }^T \,\partial_x\left(\boldsymbol{Au^2}\right) -\Delta t \, \boldsymbol{\tilde b }^T\, \frac{\boldsymbol{A}}{\rho}\, \partial_x \boldsymbol{F}(\boldsymbol{A})\,.
\end{align}
\label{eq.finalIMEX_compact_el}
\end{subequations}
\end{proof}

\begin{proof}[Proof of the AP property for the diffusive scaling]
If $\tau_r \to 0$ while $E_0 \to \infty$ and $\eta = \tau_r E_0 $ remains finite, when considering consistent initial data as in \eqref{wellIC}, from eq.~\eqref{eq.IMEX_p} and recalling again also eq. \eqref{eq.cont_compact1} and \eqref{eq.PDE_compact1}, we can compute
\begin{subequations}
\begin{align*}
\boldsymbol{\hat{p}} &= \hat{\mathcal{A}}^{-1}\Big[ -\boldsymbol{a}\left(p^{(1)} - F(A^{(1)}) +\eta\, G(A^{(1)}) \partial_x(Au)^{(1)}\right)
- \hat{\mathcal{A}}\, \eta\,\boldsymbol{G}(\boldsymbol{\hat{A}}) \partial_x(\boldsymbol{\hat{Au}})\Big] + \boldsymbol{F}(\boldsymbol{\hat{A}})\, \\
&= \hat{\mathcal{A}}^{-1}\Big[-\boldsymbol{a}\left(p^n - F(A^n) + \eta\, G(A^n) \partial_x(Au)^n\right)
- \hat{\mathcal{A}}\, \eta\,\boldsymbol{G}(\boldsymbol{\hat{A}}) \partial_x(\boldsymbol{\hat{Au}})\Big] + \boldsymbol{F}(\boldsymbol{\hat{A}})\,\\
&= - \eta\,\boldsymbol{G}(\boldsymbol{\hat{A}}) \partial_x(\boldsymbol{\hat{Au}}) + \boldsymbol{F}(\boldsymbol{\hat{A}})\,.
\end{align*}
\end{subequations}
Also in this case, being the scheme GSA, $p^{n+1}=F(A^{n+1}) - \eta\, G(A^{n+1}) \partial_x (Au)^{n+1}$, which permits to maintain consistent initial data also for the further time steps.
Substituting the above result in eq.~\eqref{eq.mom_compact2}, we correctly recover a consistent explicit discretization for the parabolic system \eqref{completesyst_KV}, having internal stages
\begin{subequations}
\begin{align}
	\boldsymbol{A} &= A^n \boldsymbol{e} -  \Delta t \,\tilde{\mathcal{A}} \, \partial_x( \boldsymbol{Au})
	\\
	\boldsymbol{Au} &= (Au)^n \boldsymbol{e} -  \Delta t \,\tilde{\mathcal{A}}  \,\partial_x\left(\boldsymbol{Au^2}\right) -\Delta t \, \tilde{\mathcal{A}} \,\frac{\boldsymbol{A}}{\rho}\, \partial_x \boldsymbol{F}(\boldsymbol{A}) + \Delta t \,\tilde{\mathcal{A}}\,\frac{\boldsymbol{A}}{\rho}\,\partial_x \left(\eta\,\boldsymbol{G}(\boldsymbol{A}) \partial_x(\boldsymbol{Au})\right)\,,
\end{align}
\label{eq.iterIMEX_compact_KV}
\end{subequations}
and final solution
\begin{subequations}
\begin{align}
	A^{n+1} &= A^n -  \Delta t \,\boldsymbol{\tilde b }^T \, \partial_x(\boldsymbol{Au})
	\\
	(Au)^{n+1} &= (Au)^n -  \Delta t\,\boldsymbol{\tilde b }^T  \partial_x\left(\boldsymbol{Au^2}\right) -\Delta t\, \boldsymbol{\tilde b}^T \frac{\boldsymbol{A}}{\rho} \partial_x \boldsymbol{F}(\boldsymbol{A}) + \Delta t\, \boldsymbol{\tilde b }^T\,\frac{\boldsymbol{A}}{\rho}\,\partial_x \left(\eta\,\boldsymbol{G}(\boldsymbol{A}) \partial_x(\boldsymbol{Au})\right).
\end{align}
\label{eq.finalIMEX_compact_KV}
\end{subequations}
\end{proof}

The above analyses show that the proposed IMEX Runge-Kutta scheme provides a high order temporal discretization that is consistent with the behavior of the continuous model even in the asymptotic regimes, meaning that the scheme is satisfying not only the AP property but also the property of being \emph{asymptotically accurate} \cite{pareschi2005}.

\begin{remark}^^>
\label{remarkIMEX}
\begin{itemize}
\item 
We emphasize that one could choose a different IMEX partitioning of the problem, following the methodology proposed in \cite{boscarino2017}. This method treats partly explicitly and partly implicitly the terms in the second and third equation of system \eqref{completesyst}, reading in compact form as
\begin{subequations}
\begin{align}
	\boldsymbol{A} &= A^n \boldsymbol{e} -  \Delta t \,\tilde{\mathcal{A}} \, \partial_x( \boldsymbol{Au})
	\\
	\boldsymbol{Au} &= (Au)^n \boldsymbol{e} -  \Delta t \,\tilde{\mathcal{A}} \,\partial_x\left(\boldsymbol{Au^2}\right) -\Delta t \, \mathcal{A}\, \frac{\boldsymbol{A}}{\rho} \partial_x \boldsymbol{p}
	\\
	\boldsymbol{p} &= p^n \boldsymbol{e} - \Delta t \,\tilde{\mathcal{A}} \, E_0 \boldsymbol{G}(\boldsymbol{A}) \partial_x(\boldsymbol{Au}) - \frac{\Delta t \,\mathcal{A}}{\tau_r}\left(\boldsymbol{p} - \boldsymbol{F}(\boldsymbol{A})\right).
\end{align}
\end{subequations}
We have tested also the above approach by obtaining analogous results to those here presented. Note, however, that the above method requires the additional condition $\boldsymbol{\tilde c}=\boldsymbol{c}$ in order to preserve the stationary solutions of the problem.
We refer to \cite{boscarino2017a,boscarino2017} for more details.

\item It is here stressed that the construction of numerical methods that preserve the order of accuracy also in regimes described by model \eqref{completesyst_newVisco}, arising from the perturbation of the local elastic equilibrium, requires additional conditions, which are not addressed in the context of this paper. The reader may refer to \cite{boscarino2017a} for further discussion of these issues.
\end{itemize}
\end{remark}

\subsection{Well-balanced property}
For consistent initial data, we can show that the partitioned IMEX scheme \eqref{eq.iterIMEX}-\eqref{eq.finalIMEX} results \emph{well-balanced in time}, meaning that, by ignoring the space discretization error, it preserves stationary solutions \cite{bertaglia2020,muller2013,pimentel-garcia2023}. 
More precisely, we have
\begin{subequations}
\begin{align}
\partial_x(Au)^n &=0\\
\rho\,\partial_x(Au^2)^n + A^n \, \partial_x p^n &= 0 \quad \implies \quad \boldsymbol{Q}^{n+1} = \boldsymbol{Q}^{n}\\
p^n - F(A^n) + \tau_r \,E_0\, G(A^n) \,\partial_x(Au)^n&= 0\,,
\end{align}
\label{WBTh}
\end{subequations}
where $\boldsymbol{Q}=(A,Au,p)^T$. Note that the last equation is equivalent to assume that the initial condition is well prepared, as in Definition \ref{def.wellIC}.

\begin{proof}[Proof of the well-balanced property]
Let us consider an initial condition as in \eqref{WBTh} for the model \eqref{completesyst}. From the scheme \eqref{eq.finalIMEX_compact}, we have $\boldsymbol{Q}^{n+1} = \boldsymbol{Q}^{n}$ if
\begin{subequations}
\begin{align*}
	\Delta t \,\boldsymbol{\tilde b}^{T}\, \partial_x (\boldsymbol{Au}) &= 0
	\\
	\Delta t \,\boldsymbol{\tilde{b}}^{T} \left( \partial_x \left(\boldsymbol{Au^2}\right) + \frac{\boldsymbol{A}}{\rho} \partial_x \boldsymbol{p} \right) &= 0
	\\
	\Delta t \,\boldsymbol{b}^{T} \left( E_0\boldsymbol{G}(\boldsymbol{A}) \partial_x(\boldsymbol{Au}) + \frac{1}{\tau_r}\left(\boldsymbol{p} - \boldsymbol{F}(\boldsymbol{A})\right) \right) &= 0.
\end{align*}
\end{subequations}
This is guaranteed if $\boldsymbol{A}=A^{n}\boldsymbol{e}$, $\boldsymbol{Au}=(Au)^{n}\boldsymbol{e}$, $\boldsymbol{p}=p^{n}\boldsymbol{e}$. We can verify that the just mentioned solution is indeed the unique one for system \eqref{eq.iterIMEX_compact1}. In fact, knowing already that the first step of the method reads $A^{(1)} = A^n$, $Au^{(1)} = (Au)^n$, $p^{(1)} = p^n$, we only need to satisfy additionally
\begin{subequations}
\begin{align*}
	\boldsymbol{\tilde{a}} \, \partial_x (Au)^{(1)} +  \hat{\tilde{\mathcal{A}}} \, \partial_x( \boldsymbol{\hat{Au}}) &= 0\\
	\boldsymbol{\tilde{a}} \left( \partial_x(Au^2)^{(1)} + \frac{A^{(1)}}{\rho} \partial_x p^{(1)} \right) +  \hat{\tilde{\mathcal{A}}} \left(\partial_x(\boldsymbol{\hat{Au}^2}) + \frac{\boldsymbol{\hat{A}}}{\rho} \partial_x \boldsymbol{\hat{p}} \right)&= 0\\
	\boldsymbol{a} \left( \frac{p^{(1)} - F(A^{(1)})}{\tau_r} + E_0\, G(A^{(1)})\, \partial_x(Au)^{(1)} \right) + \hat{\mathcal{A}} \left(\frac{ \boldsymbol{\hat{p}} - \boldsymbol{F}(\boldsymbol{\hat{A}})}{\tau_r} + E_0 \,\boldsymbol{G}(\boldsymbol{\hat{A}})\, \partial_x(\boldsymbol{\hat{Au}}) \right) &=0.
\end{align*}
\end{subequations}
To this end, we observe that, under the assumptions \eqref{WBTh}, it is always verified that
\begin{subequations}
\begin{align*}
	\partial_x(Au)^{(1)} &= 0
	\\
	\rho\, \partial_x(Au^2)^{(1)} + A^{(1)}\, \partial_x p^{(1)} &= 0
	\\
	p^{(1)} - F(A^{(1)}) + \tau_r\,E_0\, G(A^{(1)}) \,\partial_x(Au)^{(1)} &= 0
	\\
	\partial_x( \boldsymbol{\hat{Au}}) &= 0
	\\
	 \rho\,\partial_x(\boldsymbol{\hat{Au}^2}) + \boldsymbol{\hat{A}}\, \partial_x \boldsymbol{\hat{p}} &= 0
	\\
	\boldsymbol{\hat{p}} - \boldsymbol{F}(\boldsymbol{\hat{A}}) + \tau_r\, E_0\, \boldsymbol{G}(\boldsymbol{\hat{A}})\, \partial_x(\boldsymbol{\hat{Au}}) &= 0\,,
\end{align*}
\end{subequations}
which concludes the proof.
\end{proof}


\subsection{Boundary conditions}
In order to simulate realistic scenarios, inflow and outflow boundary conditions are imposed by prescribing an input flow rate waveform in time and coupling the output of the 1D model to a 0D RCR model (a lumped-parameter model, also known as 3-element Windkessel model) \cite{quarteroni2003,quarteroni2017}. For simplicity, the treatment is presented for a first order of accuracy. For extensions to higher orders, the reader can refer to \cite{piccioli2021,muller2015}.

To evaluate the fluxes at the inlet of the domain, knowing at each Runge-Kutta time steps the inlet flow rate $q_{in}$, we recur to the Riemann Invariants $\Gamma_1$ and $\Gamma_3$ in eq. \eqref{eq.RI} and solve the system
\begin{subequations}
\begin{align}
&q_{in}=A_{in}u_{in}\\
&u_{in} - \int\frac{c(A_{in},E_{0,in})}{A_{in}}\,\d A = u_{1} - \int\frac{c(A_{1},E_{0,1})}{A_{1}}\,\d A\\
&p_{in} - \frac{E_{\infty,in}}{W} \left[\left(\frac{A_{in}}{A_{0,in}}\right)^m - \left(\frac{A_{in}}{A_{0,in}}\right)^n\right] = p_{1} - \frac{E_{\infty,1}}{W}  \left[\left(\frac{A_{1}}{A_{0,1}}\right)^m - \left(\frac{A_{1}}{A_{0,1}}\right)^n\right]\,,
\end{align}
\end{subequations}
with unknowns $A_{in}$, $u_{in}$, $p_{in}$ and known variables in the first cell of the domain, $A_{1}$, $u_{1}$, $p_{1}$, $A_{0,1}$, $E_{0,1}$, $E_{\infty,1}$. We remark that the definition of the celerity $c$ is given in eq. \eqref{eq:cel} and involves the variables $A$ and $E_0$.
Under the assumptions $A_{0,in} = A_{0,1}$, $E_{0,in} = E_{0,1}$ and $E_{\infty,in} = E_{\infty,1}$, the system can be numerically solved recurring to a Newton-Raphson iterative procedure to find $A_{in}$ and then straightforwardly compute $u_{in}$ and $p_{in}$ with $\Gamma_1$ and $\Gamma_3$, respectively.

At the outlet of the domain, an RCR circuit (analogous to the corresponding electrical circuit), consisting of a resistor with resistance $R_1$ connected in series with a parallel combination of a second resistor with resistance $R_2$ and a capacitor with compliance $C$, is used to simulate the effects of resistance and compliance of the terminal arteries on the propagation of pressure waves (see Fig. \ref{fig.TCstent_scheme} for a sketch). The RCR ordinary differential equation (ODE) model reads \cite{piccioli2021}:
\begin{subequations}
\begin{align}
&\frac{\d p_C}{\d t} = \frac1{C R_1}\left( p_* - p_C \right) - \frac1{C R_2}\left( p_C - p_{out} \right)\\
&A_*u_* = \frac{p_* - p_C}{R_1}\,,
\end{align}
\label{eq.RCR}
\end{subequations}
where $p_C$ is the pressure at the capacitor, $A_*$, $u_*$ and $p_*$ are the unknown area, velocity and pressure, respectively, at the interface between 1D and 0D model, and $p_{out}$ is the pressure at the outlet of the RCR unit, fixed to be $p_{out}=0$ to mimic the blood pressure when the flux reaches the venous system. We discretize in time the ODE explicitly due to the absence of stiff terms, and couple the problem with the Riemann Invariants $\Gamma_2$ and $\Gamma_3$ obtaining the following system to be solved for each Runge-Kutta time step $\Delta t = t^{k+1}-t^k$ of the numerical method:
\begin{subequations}
\begin{align}
&p_{C}^{k+1}=p_C^k + \frac{\Delta t}{CR_1}\left( p_*^{k+1} - p_C^k\right) - \frac{\Delta t}{CR_2}\left( p_C^k - p_{out}\right)\\
&A_*^{k+1}u_*^{k+1} = \frac{p_*^{k+1} - p_C^{k+1}}{R_1} \\
&u_*^{k+1} + \int\frac{c(A_*^{k+1},E_{0,*})}{A_*^{k+1}}\,\d A = u_{N}^{k+1} + \int\frac{c(A_{N}^{k+1},E_{0,N}^{k+1})}{A_{N}^{k+1}}\,\d A\\
&p_*^{k+1} - \frac{E_{\infty,*}}{W} \left[\left(\frac{A_*^{k+1}}{A_{0,*}^{k+1}}\right)^m - \left(\frac{A_*^{k+1}}{A_{0,*}^{k+1}}\right)^n\right] = p_{N}^{k+1} - \frac{E_{\infty,N}}{W} \left[\left(\frac{A_N^{k+1}}{A_{0,N}^{k+1}}\right)^m - \left(\frac{A_N^{k+1}}{A_{0,N}^{k+1}}\right)^n\right]\,.
\end{align}
\end{subequations}
Here, $A_{N}$, $u_{N}$, $p_{N}$, $A_{0,N}$, $E_{0,N}$ and $E_{\infty,N}$ are the known values of variables in the last cell of the 1D domain.
Similarly to the inlet boundary, under the assumptions $A_{0,*} = A_{0,N}$, $E_{0,*} = E_{0,N}$ and $E_{\infty,*} = E_{\infty,N}$, the above non-linear system can be numerically solved recurring to a Newton-Raphson method to compute $A_*^{k+1}$, and then directly derive $u_*^{k+1}$ and $p_*^{k+1}$ through $\Gamma_2$ and $\Gamma_3$, respectively, as well as $p_{C}^{k+1}$, which will be used for the following time step \cite{piccioli2021,bertaglia2020a}.

\section{Numerical tests and applications}
\label{section_numericalresults}
In this section, we present several numerical tests that permit to validate the proposed methodology. First, an accuracy analysis of the method is performed considering three different configurations of the scaling parameters of the augmented blood flow model, accounting for the different asymptotic behaviors discussed in Section \ref{asymptotics}. Then, five Riemann problems (RP) are executed, to test the methodology in presence of variables discontinuous in space, again taking into account different constitutive settings. The first of these problems aims to verify also numerically the well-balancing of the method. Finally, a more applied test case is performed taking into account a multiscale configuration of the rheological parameters. In this test, the haemodynamics of a thoracic aorta is simulated under normal, healthy conditions, and compared with that in the presence of a stent/prosthesis in the center of the vessel. The latter constitutes an element of increased, localized, wall stiffening that it is considered to lead the model to the parabolic scaling in the stretch where the stent is situated.
In all the simulations, for stability we fix $\CFL = 0.9$ and $\nu = 0.5$. Finally, if not otherwise stated, in the WENO reconstruction to compute the nonlinear weights we fix the positive parameter $\varepsilon = 10^{-6}$ \cite{shu1998}.

\begin{table}[!htb]
\caption{$L_2$ error and empirical order of accuracy in the state variables $A$, $Au$, $p$ for different values of the scaling parameters, while leaving $\bar{E}_{\infty} = 8\cdot10^{5}$ Pa fixed. Hyperbolic viscoelastic case (first column), parabolic viscoelastic case (second column), and purely elastic case (third column). $N_x$ is the number of cells in the computational domain.} \label{tab:accuracy} 
\centering
\begin{tabular}{| c | l | c c | c c | c c |}
\toprule
\multirow{5}{*}{Variable} &\multirow{5}{*}{$N_x$} &\multicolumn{2}{c|}{SLS} &\multicolumn{2}{c|}{KV} &\multicolumn{2}{c|}{EL}\\
& &\multicolumn{2}{c|}{$\tau = 10^{-1}$ s} &\multicolumn{2}{c|}{$\tau = 10^{-4}$ s} &\multicolumn{2}{c|}{$\tau = 0$ s}\\
& &\multicolumn{2}{c|}{$\eta = 5\cdot10^{5}$ Pa$\cdot$s} &\multicolumn{2}{c|}{$\eta = 5\cdot10^{4}$ Pa$\cdot$s} &\multicolumn{2}{c|}{$\eta = 0$ Pa$\cdot$s}\\
& &\multicolumn{2}{c|}{$\bar{E}_0 = 10^{6}$ Pa} &\multicolumn{2}{c|}{$\bar{E}_0 = 5\cdot10^{8}$ Pa} &\multicolumn{2}{c|}{$\bar{E}_0 = 8\cdot10^{5}$ Pa}\\
\cline{3-8}
& & $L_2$ error &order & $L_2$ error &order  & $L_2$ error &order \\
\midrule
\multirow{4}{*}{$A$}
     &15 & 4.21e-03 &        & 2.37e-02 &         & 4.37e-03 &                  \\ 
     &45 & 1.98e-04 &    2.78 & 2.29e-03 &    2.13 & 2.32e-04 &    2.67  \\ 
    &135 & 8.64e-06 &    2.85 & 7.83e-05 &    3.07 &  9.69e-06 &    2.89 \\ 
    &405 & 2.26e-07 &    3.32 & 1.04e-06 &    3.94 & 2.50e-07 &    3.33 \\ 
\hline 
\multirow{4}{*}{$Au$}
     &15 & 3.64e-02 &         & 1.84e-01 &         & 4.16e-02 &                  \\ 
     &45 & 1.77e-03 &    2.75 & 7.88e-03 &    2.87 & 1.86e-03 &    2.83 \\ 
    &135 & 5.46e-05 &    3.16 & 1.81e-04 &    3.43 & 5.98e-05 &    3.13 \\ 
    &405 & 1.47e-06 &    3.29 &  3.66e-06 &    3.55 & 1.66e-06 &    3.26  \\ 
\hline
\multirow{4}{*}{$p$}
	&15 & 1.05e-03 &        &  6.89e-03 &         & 9.86e-04 &          \\ 
     &45 & 5.37e-05 &    2.71 & 6.28e-04 &    2.18  & 5.35e-05 &    2.65 \\ 
    &135 & 2.34e-06 &    2.85 & 2.08e-05 &    3.10& 2.23e-06 &    2.89 \\ 
    &405 & 4.22e-08 &    3.66 &  2.69e-07 &    3.96 & 4.28e-08 &    3.60 \\ 
\bottomrule
\end{tabular} 
\end{table}
\begin{figure}[!htb]
\centering
\subfloat{\includegraphics[width=0.45\textwidth]{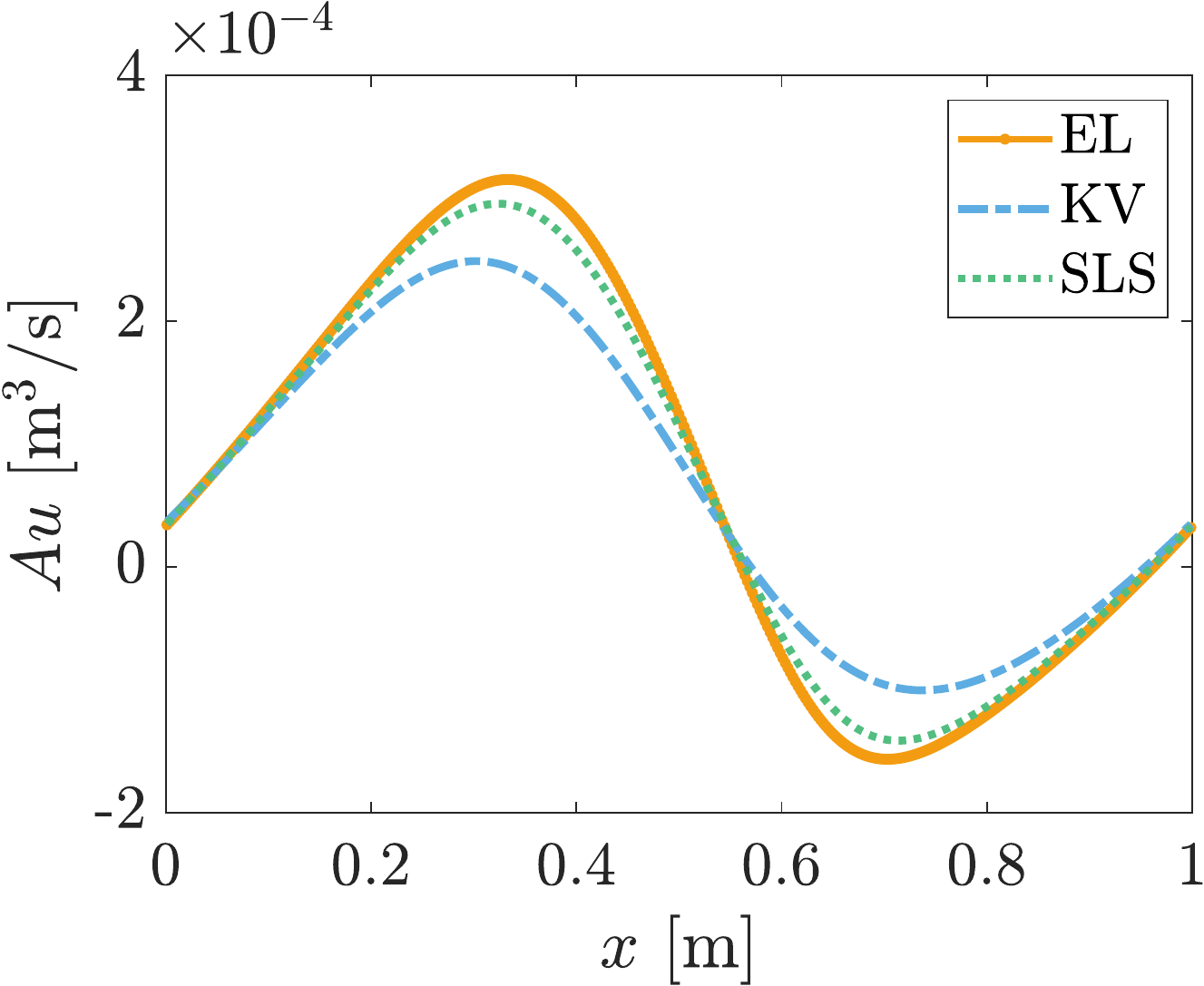}}
\hfil
\subfloat{\includegraphics[width=0.45\textwidth]{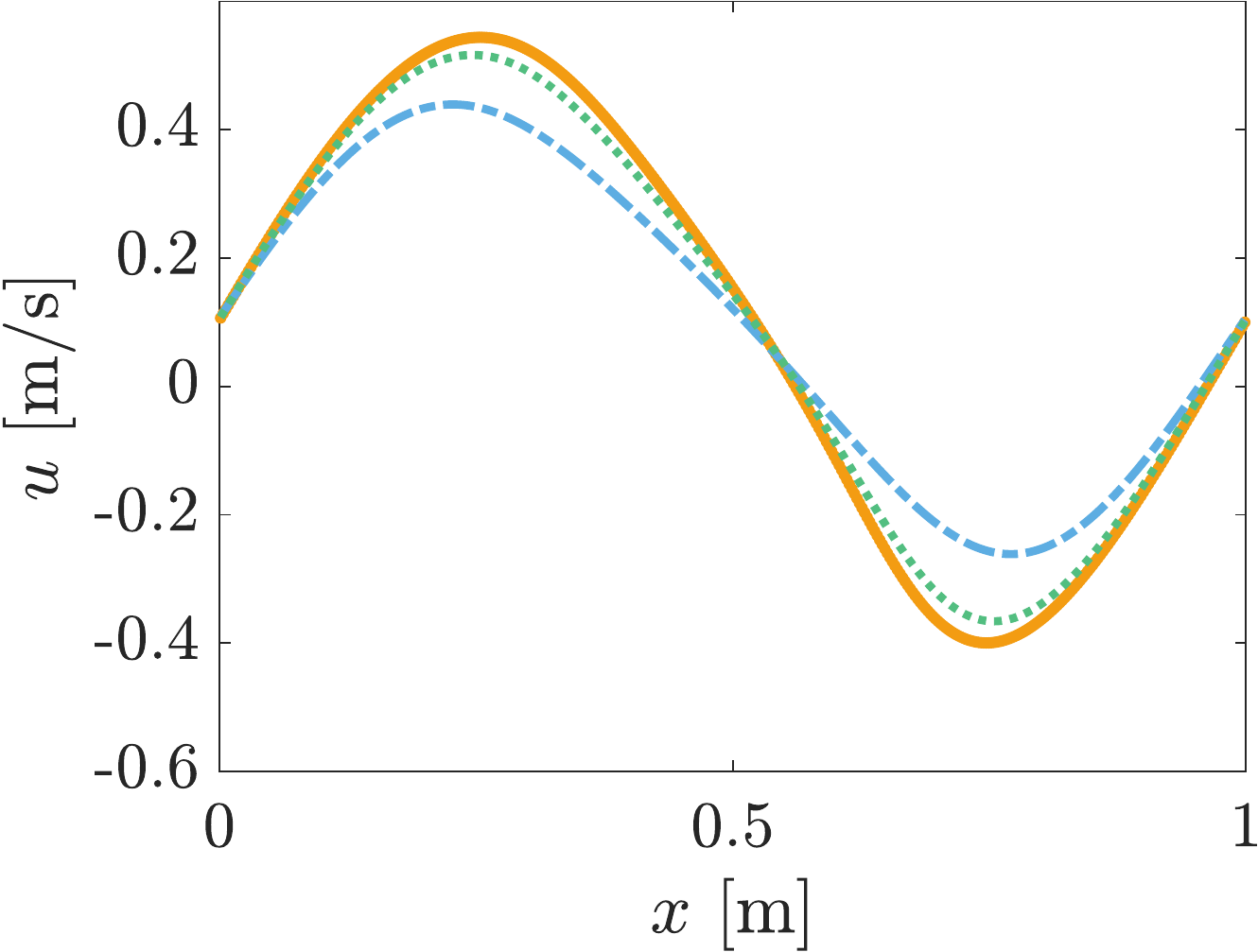}}
\hfil
\subfloat{\includegraphics[width=0.45\textwidth]{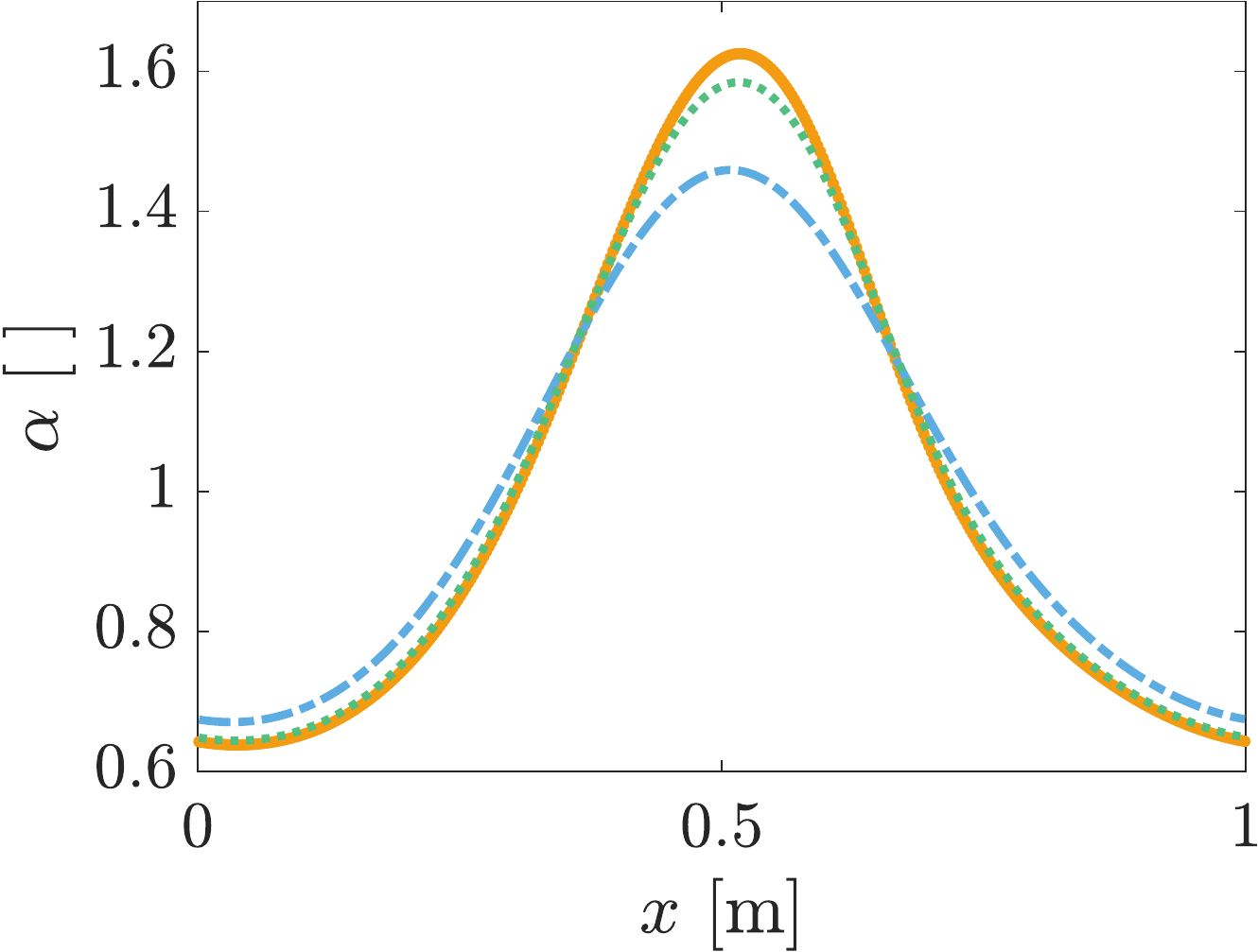}}
\hfil
\subfloat{\includegraphics[width=0.45\textwidth]{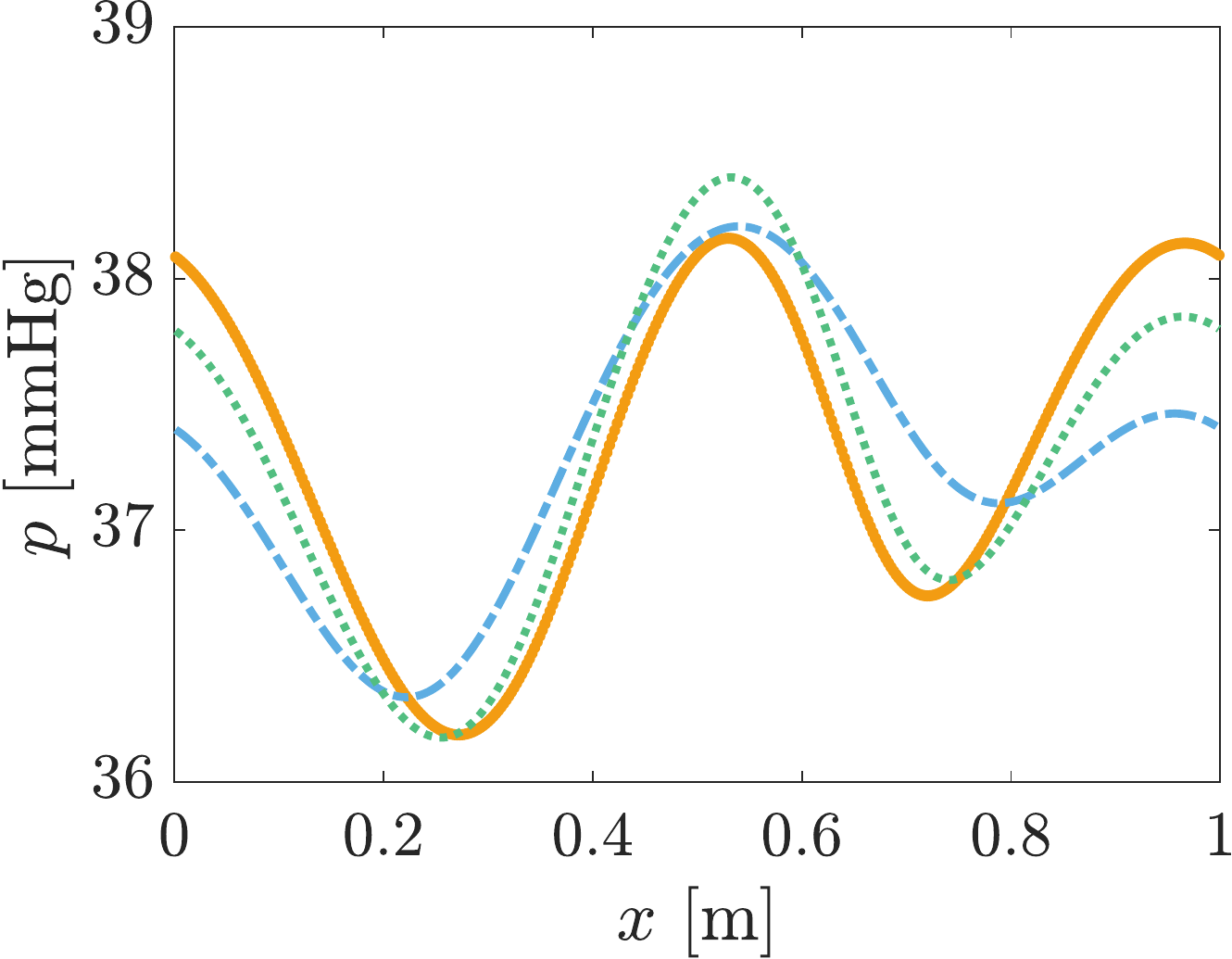}}
\caption{Comparison of the results at $t_{end}=0.25$ obtained in the accuracy analysis with the three different settings: hyperbolic viscoelastic case (SLS), parabolic viscoelastic case (KV), and elastic case (EL).}
\label{fig.AA}
\end{figure}
\subsection{Accuracy analysis}
To empirically verify the accuracy of the proposed method, it is important to remember that, due to the presence of variables in the system \eqref{completesyst} that have widely varying orders of magnitude, it is necessary to switch to the dimensionless form of the equations to avoid numerical errors, especially when using high-order methods. \cite{pimentel-garcia2023}. For details on the derivation of the dimensionless form the reader can refer to the Appendix \ref{section_appendix}.

We design a numerical test with periodic boundary conditions and smooth initial distributions of the variables:
\[ A_0(x,0) = \bar A + \bar{a} \, \sin\left(\frac{2\pi x}{L}\right)\,, \quad
p_0(x,0) = \bar P_0 + \bar p \, \sin\left(\frac{2\pi x}{L}\right)\,, \]
\[ E_0(x,0) = \bar E_0 + \bar e \, \sin\left(\frac{2\pi x}{L}\right)\,, \quad
E_{\infty}(x,0) = \bar E_{\infty} + \bar e \, \sin\left(\frac{2\pi x}{L}\right)\,. \]
Then, we fix $Au(x,0) = 5 \cdot 10^{-5}$ m$^3$/s, $A(x,0)=A_0(x,0)$ and evaluate $p(x,0)$ through the elastic tube law, eq. \eqref{elastictubelaw}, knowing the initial condition of the rest of the variables.
We consider a hypothetical artery of length $L=1$ m and wall thickness $h_0 = 1.5$ mm, with $\bar A = 5$ cm$^2$, $\bar a = 1$ cm$^2$, $\bar P_0 = 5$ kPa, $\bar p = 500$ Pa, $\bar e = 0.2$ MPa and $\bar E_{\infty} =0.8$ MPa, and blood density $\rho = 1050$ kg/m$^3$. The value of $\bar E_0$, together with that of $\tau_r$ and, consequently, $\eta$, are set accordingly to three different combinations, in order to perform an accuracy analysis of the method in all the configurations given by the asymptotic limits of the proposed blood flow model in terms of constitutive behavior of the vessel wall: a hyperbolic viscoelastic case (SLS), a parabolic viscoelastic case (KV), and a hyperbolic elastic case (EL).
The values of the scaling parameters are listed in Table \ref{tab:accuracy}, where the $L_2$ relative error norms and empirical order of accuracy obtained at $t_{end}=0.25$ are also presented.
Referring to \cite{boscarino2017}, for each state variable $q$ of the system, we compute the relative $L_2$ norm for results obtained with mesh size $\Delta x$ as follows:
\[L_2^{\Delta x} = \sqrt{\frac{\sum_i^{N_x}\left(q_i^{\Delta x} - q_i^{\Delta x/3}\right)^2 }{\sum_i^{N_x} \left(q_i^{\Delta x/3}\right)^2 }},\]
where $q_i^{\Delta x}$ is the value of the state variable resulting from the simulation with mesh size $\Delta_x$ in the $i-$th cell. Thus, the order of convergence is computed as
\[p^{\Delta x/3} = \log_3 \left( \frac{L_2^{\Delta x}}{L_2^{\Delta x/3}}\right).\]
In Table \ref{tab:accuracy}, we can notice that the expected order of accuracy is plainly confirmed, even if especially the KV and EL configurations define very stiff problems. 
In Fig. \ref{fig.AA}, a comparison of the final results obtained in the three constitutive settings is presented for the main variables of the system, where the different response of the viscoelastic configurations can be appreciated.

\begin{table}[!tbp]
\caption{Domain data and initial conditions for the Riemann problems. Subscripts ``\(L\)'' and ``\(R\)'' stand for left and right state, respectively, of the piece-wise constant initial values, while $x_0$ identifies the position of the initial discontinuity. For each parameter in this list, the same value is used for the three cases (a), (b) and (c). For all the tests, the vessel wall thickness is $h_0 = 0.3$ mm and the blood density $\rho=1050$ kg/m$^3$.}
\centering
\begin{footnotesize}
\begin{tabular}{|l | ccccc |}
\toprule
Variable &RP1 &RP2 &RP3 &RP4 &RP5\\
\midrule
	$L$~[m] &0.2 &0.2 &0.2 &0.2 &0.5  \\
	\(x_0\)~[m]  &0.10 &0.05 &0.05 &0.10 &0.25  \\
	$t_{end}$~[s]  &0.100 &0.007  &0.015 &0.010 &0.050  \\
	\(A_{0,L}\)~[mm$^2$] &627.06 &156.77 &110.00 &313.53 &28.274  \\
	\(A_{0,R}\)~[mm$^2$] &313.53 &313.53  &130.00 &313.53 &29.688  \\
	\(A_L\)~[mm$^2$]  &641.38 &250.82 &99.00 &470.30 &31.00  \\
	\(A_R\)~[mm$^2$]  &312.82 &329.21  &208.00 &219.47 &31.00  \\
	\(u_L\)~[m/s]  &0.00 &1.00  &0.00 &0.00 &$- 0.20$ \\
	\(u_R\)~[m/s]  &0.00 &0.00  &0.00 &0.00 &0.10  \\
	\(p_L\)~[mmHg]  &80.00 &146.67  &9.97 &178.99 &0.9099  \\
	\(p_R\)~[mmHg]  &80.00 &108.78  &46.05 &8.05 &5.0303  \\
	\(p_{0,L}\)~[mmHg]  &75.00 &30.00  &10.00 &80.00 &0.50  \\
	\(p_{0,R}\)~[mmHg]  &85.00 &0.00   &5.00 &80.00 &0.50  \\
	\(E_{\infty,L}\)~[MPa]  &2.7655 &1.3828  &0.4604 &1.9555 &0.4000  \\
	\(E_{\infty,R}\)~[MPa]  &19.555 &19.555  &5.9153 &1.9555 &12.911  \\
\bottomrule
\end{tabular}
\label{tab.RPdata1}
\end{footnotesize}
\end{table}
\begin{table}[!tbp]
\caption{Initial conditions of mechanical parameters eventually having different values in the cases (a), (b) and (c) of the Riemann problems. Subscripts ``\(L\)'' and ``\(R\)'' stand for left and right state, respectively, of the piece-wise constant initial values.}
\centering
\begin{footnotesize}
\begin{tabular}{|l c | c c c c c |}
\toprule
Test &Case &$E_{0,L}$~[MPa] &$E_{0,R}$~[MPa] &$\eta_L$~[kPa$\cdot$s] &$\eta_R$~[kPa$\cdot$s] &$\tau_r$~[s]\\
\midrule
RP1 & &3.4569 &24.444 &8.6423 &61.111 &0.0005\\
\midrule
\multirow{3}{*}{RP2} &(a) &1.3828 &19.555 &0.00 &0.00 &0.00\\
&(b) &1.7285 &24.444 &4.3212 &61.111 &0.0005\\
&(c) &1.7285 &24.444 &86.423 &1222.2 &0.01\\
\midrule
\multirow{3}{*}{RP3} &(a) &0.4604 &5.9153 &0.00 &0.00 &0.00\\
&(b) &0.5755 &7.3941 &1.4388 &18.485 &0.0005\\
&(c) &0.5755 &7.3941 &5.7552 &73.941 &0.002\\
\midrule
\multirow{3}{*}{RP4} &(a) &1.9555 &1.9555 &0.00 &0.00 &0.00\\
&(b) &2.4444 &2.4444 &6.1111 &6.1111 &0.0005\\
&(c) &2.4444 &2.4444 &24.444 &24.444 &0.002\\
\midrule
\multirow{3}{*}{RP5} &(a) &0.400 &12.911 &0.00 &0.00 &0.00\\
&(b) &0.500 &16.139 &2.500 &80.693 &0.001\\
&(c) &0.500 &16.139 &250.00 &8069.3 &0.10\\
\bottomrule
\end{tabular}
\label{tab.RPdata2}
\end{footnotesize}
\end{table}

\subsection{Riemann problems}
Five Riemann problems, for which an exact solution is available when considering a simple elastic behavior of the vessel wall \cite{toro2013}, have been selected with reference to \cite{bertaglia2020,muller2013,pimentel-garcia2023} to test the methodology in presence of variables' discontinuities. The first RP aims to numerically verify the well-balancing of the scheme in a \emph{blood at rest} condition for a generic arterial setting. The rest of the problems, in addition to being simulated in the case of elastic wall (case (a)), have also been simulated considering two different viscoelastic characterizations (case (b) and (c)) to highlight the impact of the viscous damping of the vessel. The complete set of data and initial conditions is listed, for each RP and each constitutive framework, in Tables \ref{tab.RPdata1} and \ref{tab.RPdata2}. Each test has been run with $N_x = 100$ cells in the computational domain and setting $\varepsilon = 10^{-14}$ in the WENO reconstruction.

\paragraph{RP1} In this test, we consider a generic arterial setting and verify the well-balancing of the scheme, namely we confirm also numerically that the method preserves stationary solutions. To do so, we simulate the particular case of \emph{blood at rest} (i.e., the stationary case at zero flow rate), as taken from \cite{muller2013,pimentel-garcia2023}. The relative $L_2$ norms computed for the three main state variables $A$, $Au$ and $p$ at time $t_{end}=0.1$ (after 1040 time iterations) result $1.8392\cdot 10^{-17}$, $3.4120\cdot 10^{-17}$ and $8.5748\cdot 10^{-17}$, respectively, confirming the well-balance of the method.

\begin{figure}[!tb]
\centering
\subfloat{\includegraphics[width=0.45\textwidth]{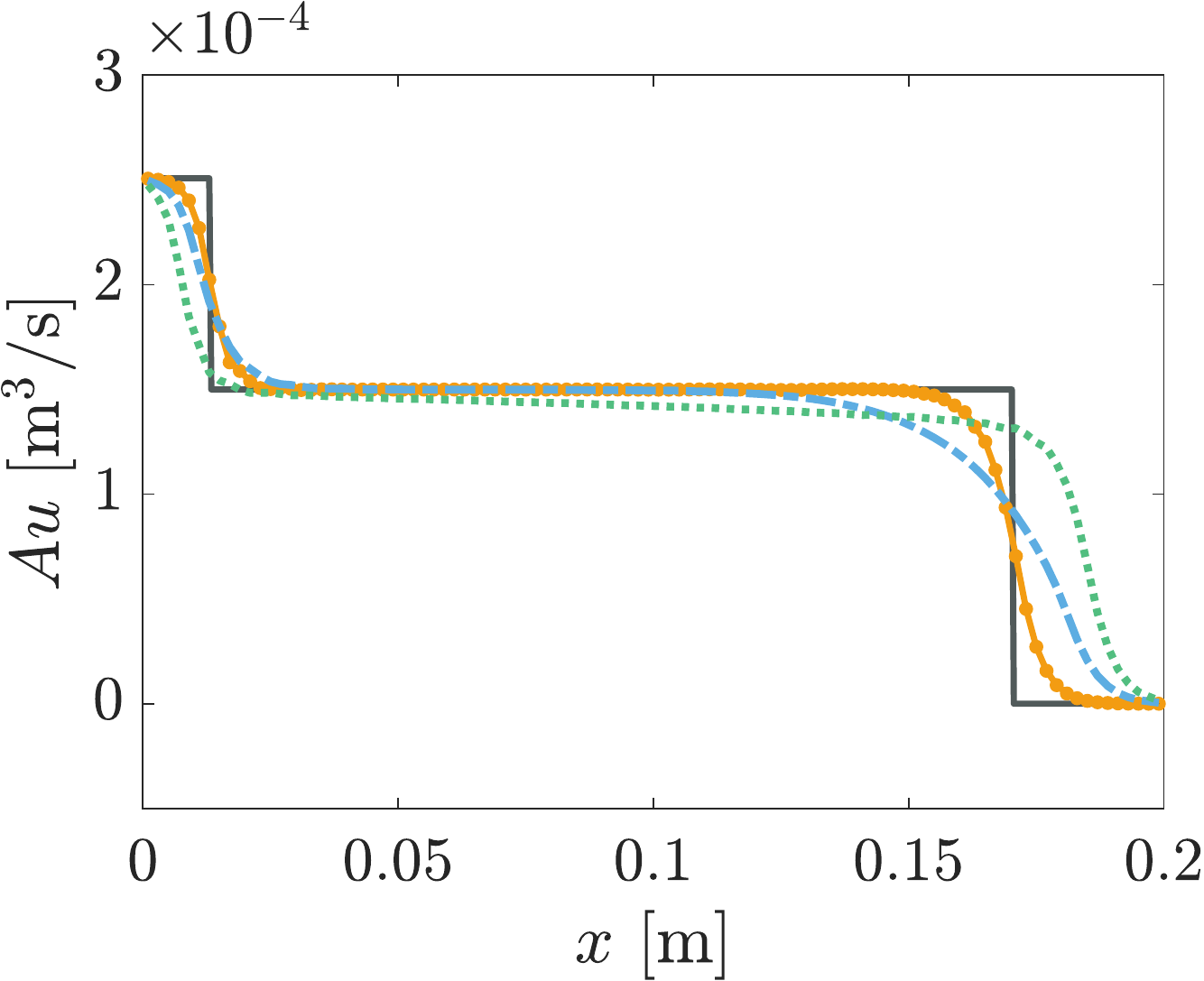}}
\hfil
\subfloat{\includegraphics[width=0.45\textwidth]{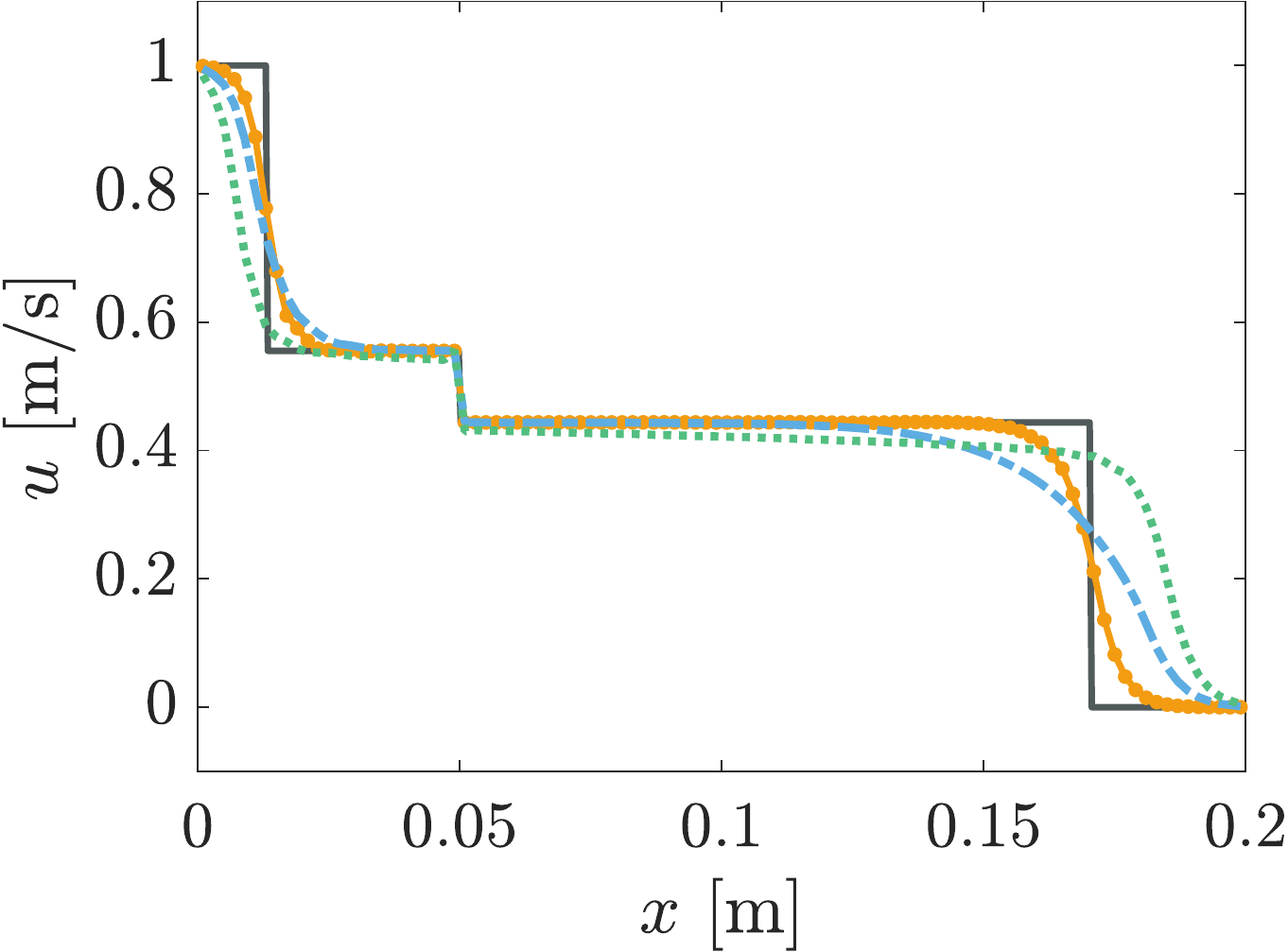}}\\
\subfloat{\includegraphics[width=0.45\textwidth]{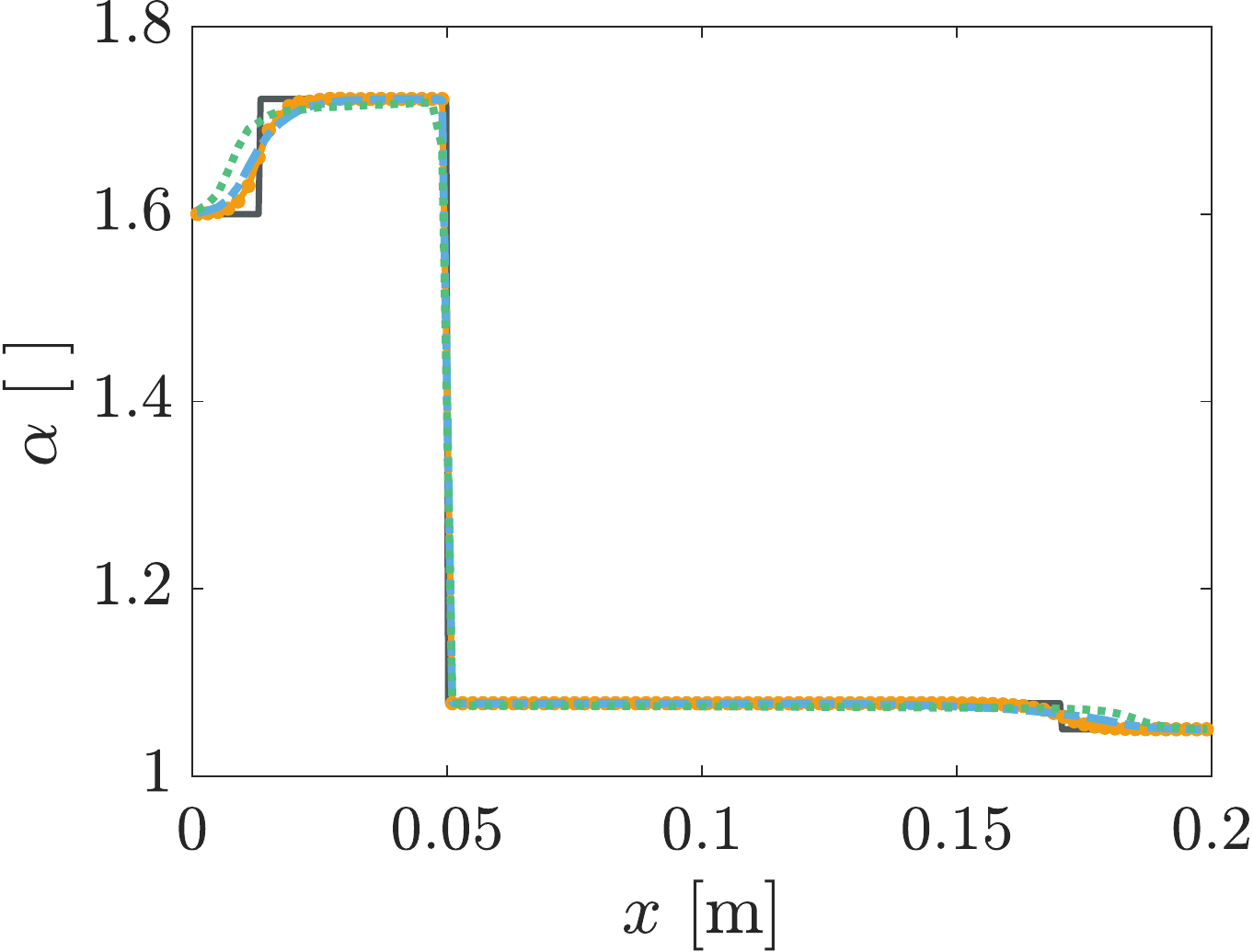}}
\hfil
\subfloat{\includegraphics[width=0.45\textwidth]{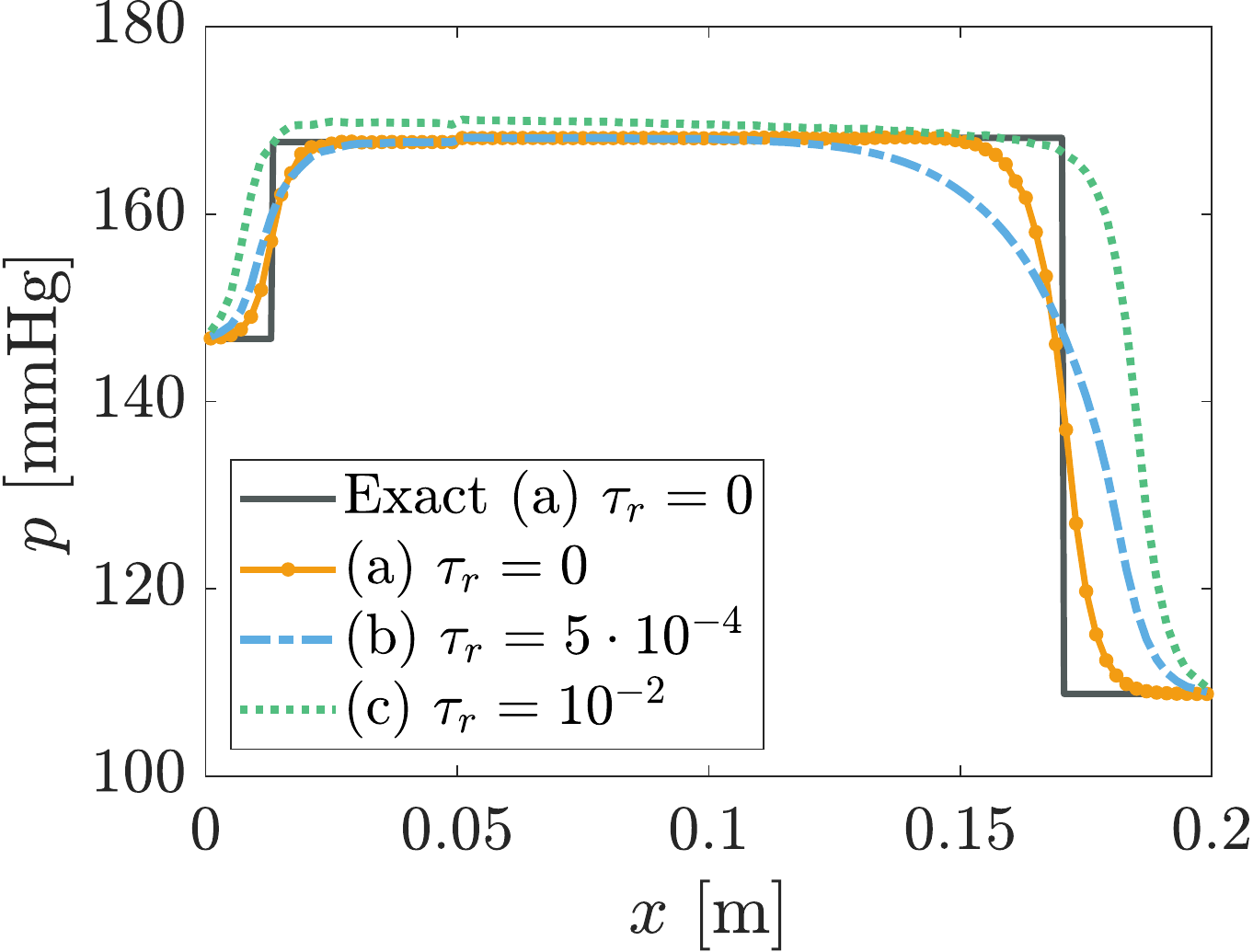}}
\caption{Comparison of the results obtained in RP2 with 3 different choices of the relaxation time $\tau_r$, which correspond to 3 different configurations of the problem (one elastic and two viscoelastic). Numerical results are plotted also against the exact solution of the elastic case ($\tau_r = 0$).}
\label{fig.RP2}
\end{figure}
\paragraph{RP2} The second Riemann problem (RP2) schematically represents the arrival of a systolic pulse pressure and, consequently, a spike in blood flow in a portion of the thoracic aorta. In this problem, the left side of the aorta, thus the part that in the initial state is reached by the systolic peak, is compressed, while to the right of the initial discontinuity the aorta is 10 times stiffer than the part to the left. This idealized configuration leads to partial reflection of the incoming wave, which can be seen in Fig. \ref{fig.RP2} by the presence of the shock wave on the left. The solution, indeed, consists of a left shock and a right shock traveling in opposite directions and separated by a stationary contact discontinuity. The elastic numerical solution results in very good agreement with the exact one. At the same time, the viscous damping effects are well visible when comparing solutions (b) and (c) with the elastic one (a), especially in the right tract of the vessel, characterized by a higher wall viscosity. In particular, an evident forward shift of the shock position is here observed.

\begin{figure}[!tb]
\centering
\subfloat{\includegraphics[width=0.45\textwidth]{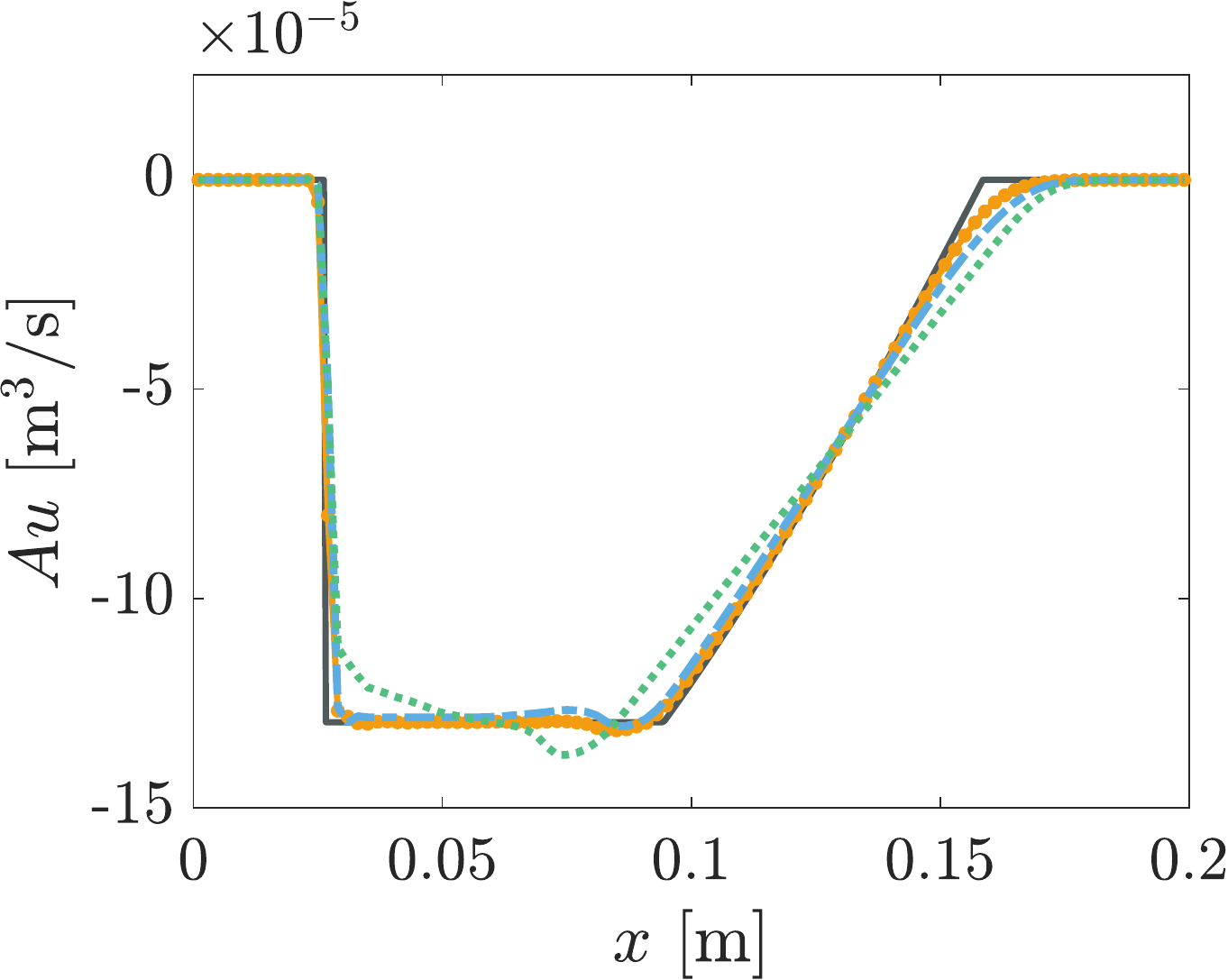}}
\hfil
\subfloat{\includegraphics[width=0.45\textwidth]{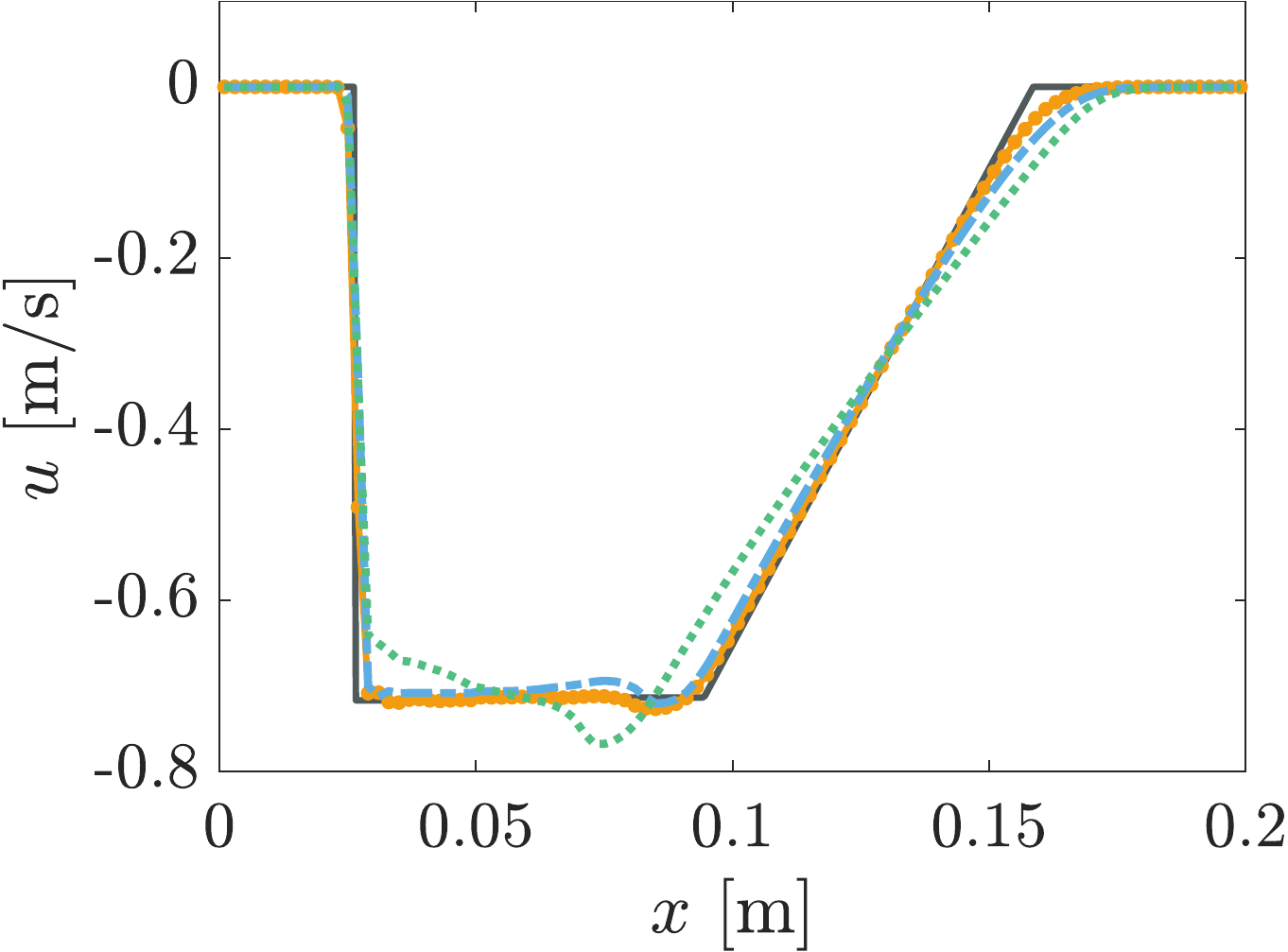}}\\
\subfloat{\includegraphics[width=0.45\textwidth]{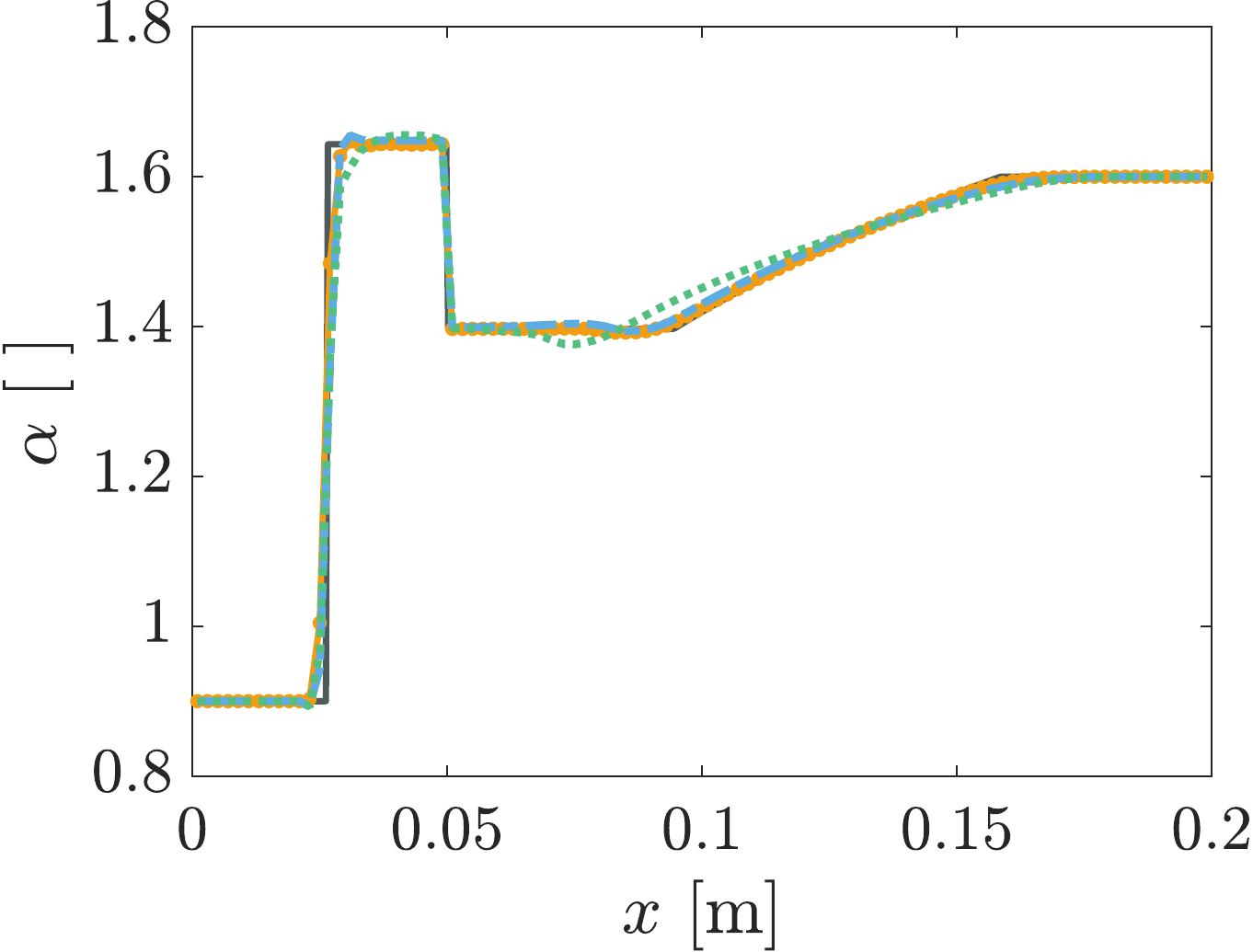}}
\hfil
\subfloat{\includegraphics[width=0.45\textwidth]{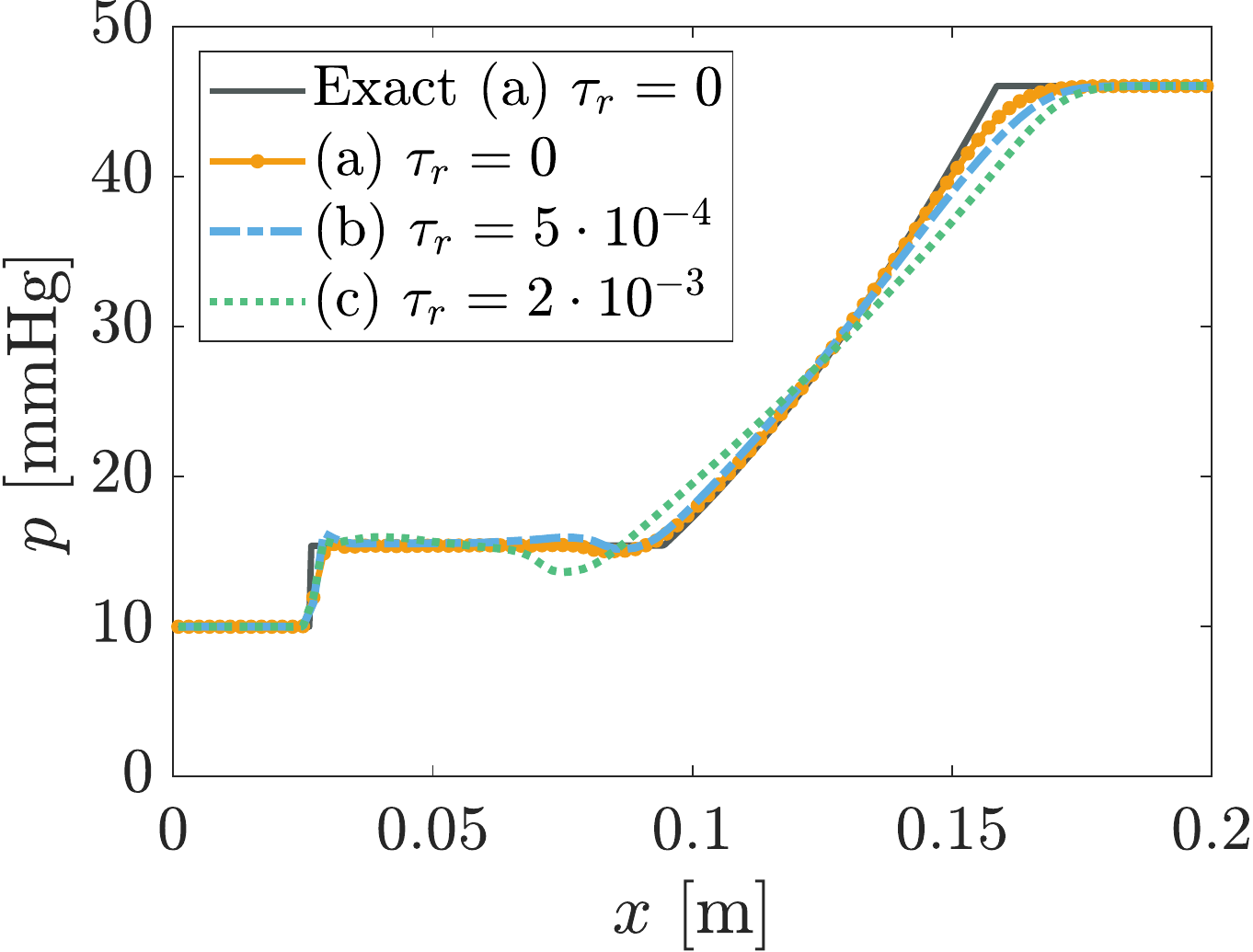}}
\caption{Comparison of the results obtained in RP3 with 3 different choices of the relaxation time $\tau_r$, which correspond to 3 different configurations of the problem (one elastic and two viscoelastic). Numerical results are plotted also against the exact solution of the elastic case ($\tau_r = 0$).}
\label{fig.RP3}
\end{figure}
\paragraph{RP3} In the third Riemann problem (RP3), the effects of a Valsalva maneuver on a portion of the internal jugular vein are schematically depicted. The Valsalva maneuver consists of forced exhalation with a closed glottis, an operation that produces a large increase in the subject's central venous pressure. In addition, in this test, an incompetent valve is considered downstream at the heart, causing venous reflux to the head. As presented in Fig. \ref{fig.RP3}, the solution consists of a left shock and a right rarefaction traveling in opposite directions and separated by a stationary contact discontinuity. Also in this very challenging test case, a good agreement of the elastic numerical result is observed with respect to the exact solution. The solution of configuration (b) differs only slightly from that of the elastic configuration (a), while solution (c) shows the effects of wall viscosity better, especially with regard to the velocity trend.

\begin{figure}[!tb]
\centering
\subfloat{\includegraphics[width=0.45\textwidth]{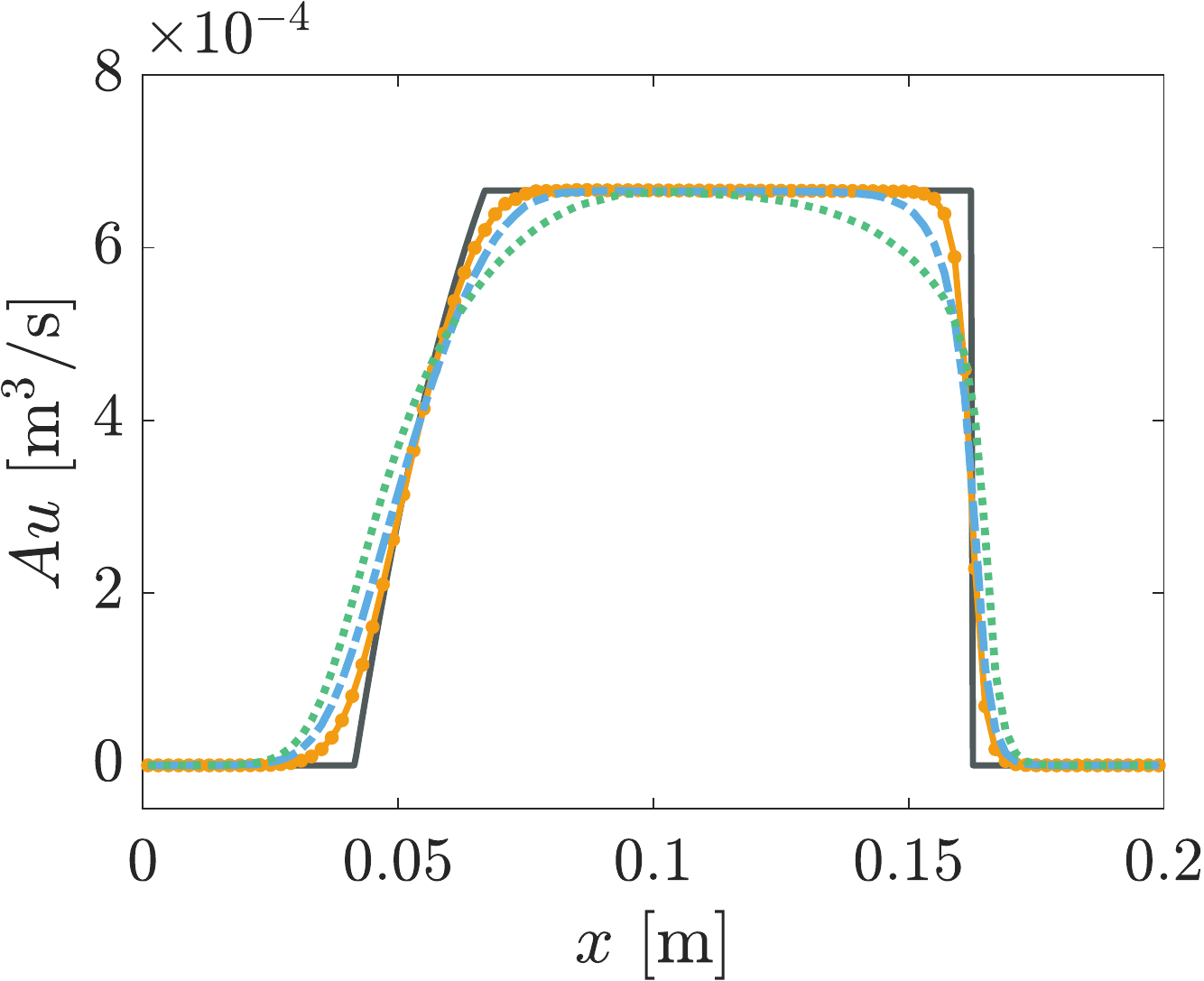}}
\hfil
\subfloat{\includegraphics[width=0.45\textwidth]{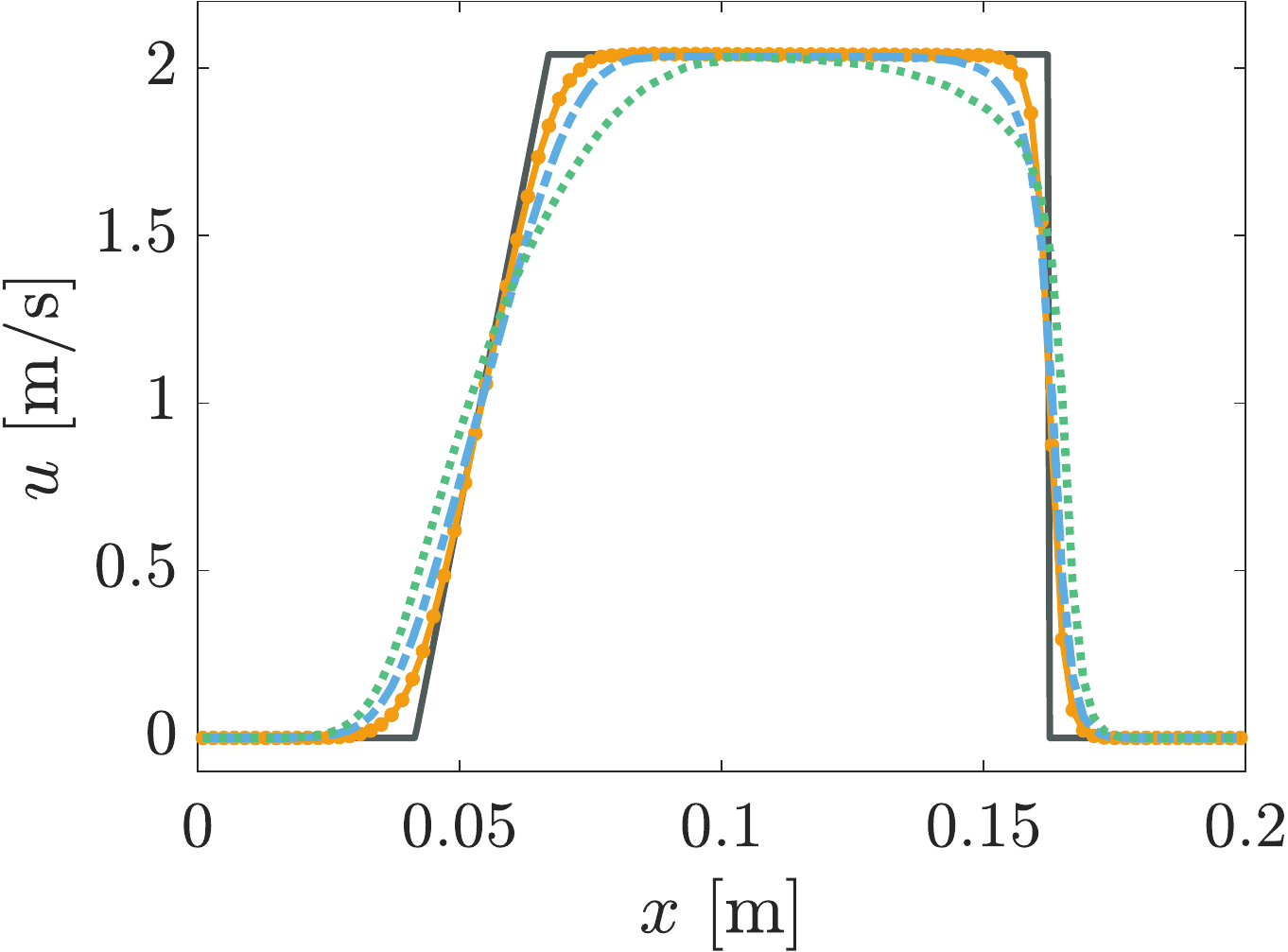}}\\
\subfloat{\includegraphics[width=0.45\textwidth]{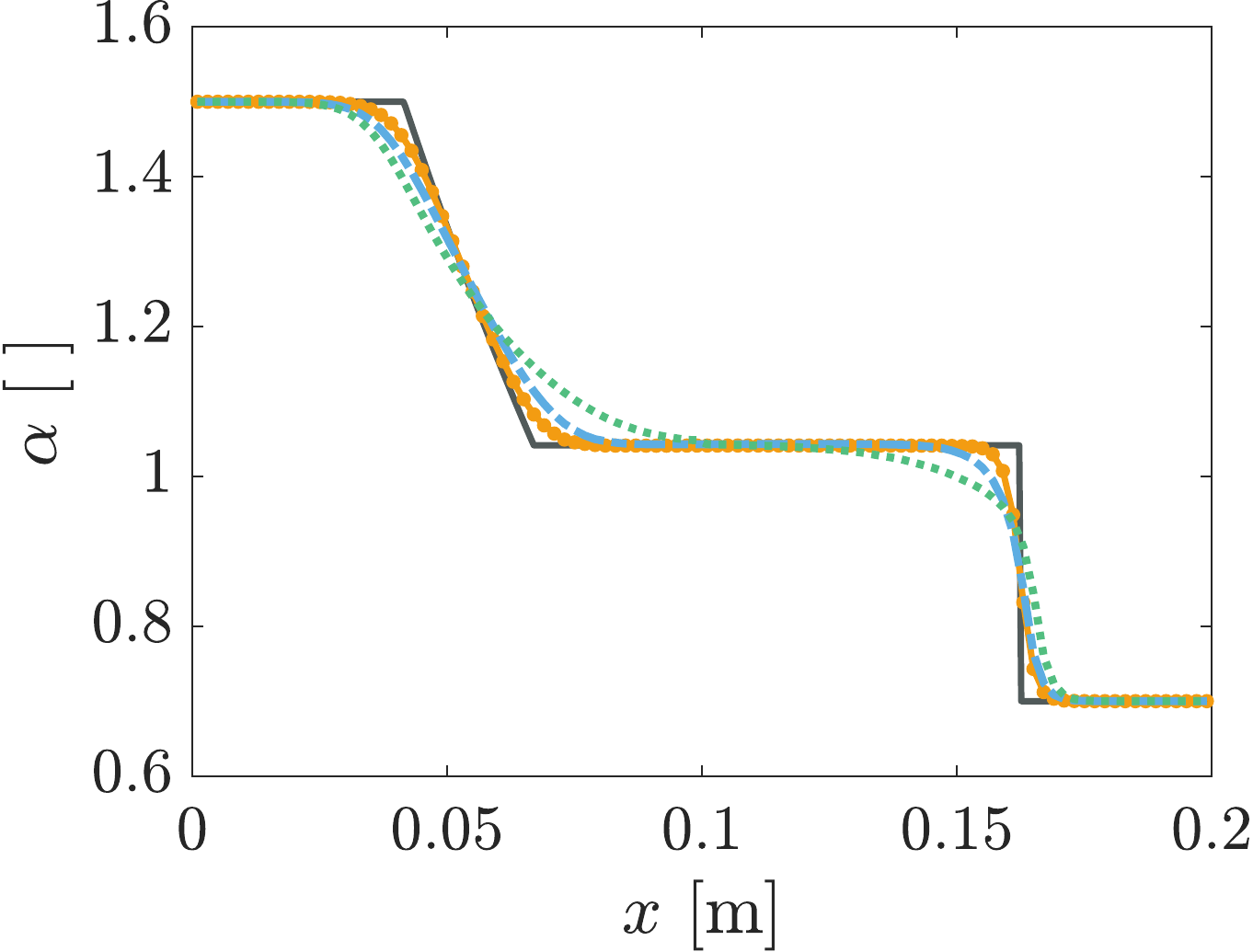}}
\hfil
\subfloat{\includegraphics[width=0.45\textwidth]{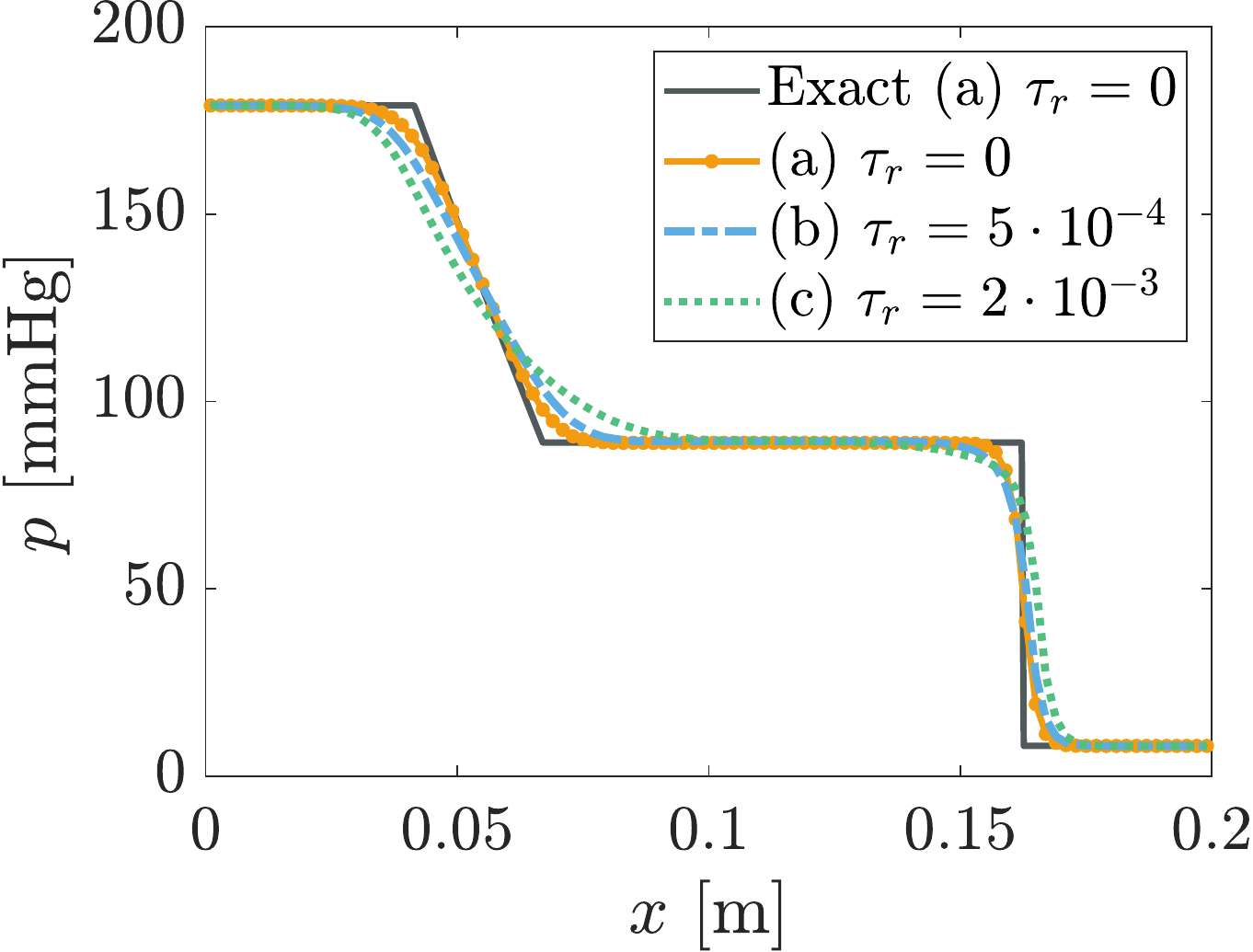}}
\caption{Comparison of the results obtained in RP4 with 3 different choices of the relaxation time $\tau_r$, which correspond to 3 different configurations of the problem (one elastic and two viscoelastic). Numerical results are plotted also against the exact solution of the elastic case ($\tau_r = 0$).}
\label{fig.RP4}
\end{figure}
\paragraph{RP4} In the fourth Riemann Problem (RP4) considered, we simulate the dynamics of a tract of thoracic aorta that initially presents all the model's variables continuous in space, except for the cross-sectional area and, consequently, the pressure. The solution of the problem, presented in Fig. \ref{fig.RP4}, shows the propagation of a left rarefaction and a right shock wave, very well captured by the model in the elastic case when comparing numerical results with the exact solution. Finally, the different effects of the chosen viscosity configurations can also be appreciated.

\begin{figure}[!htb]
\centering
\subfloat{\includegraphics[width=0.45\textwidth]{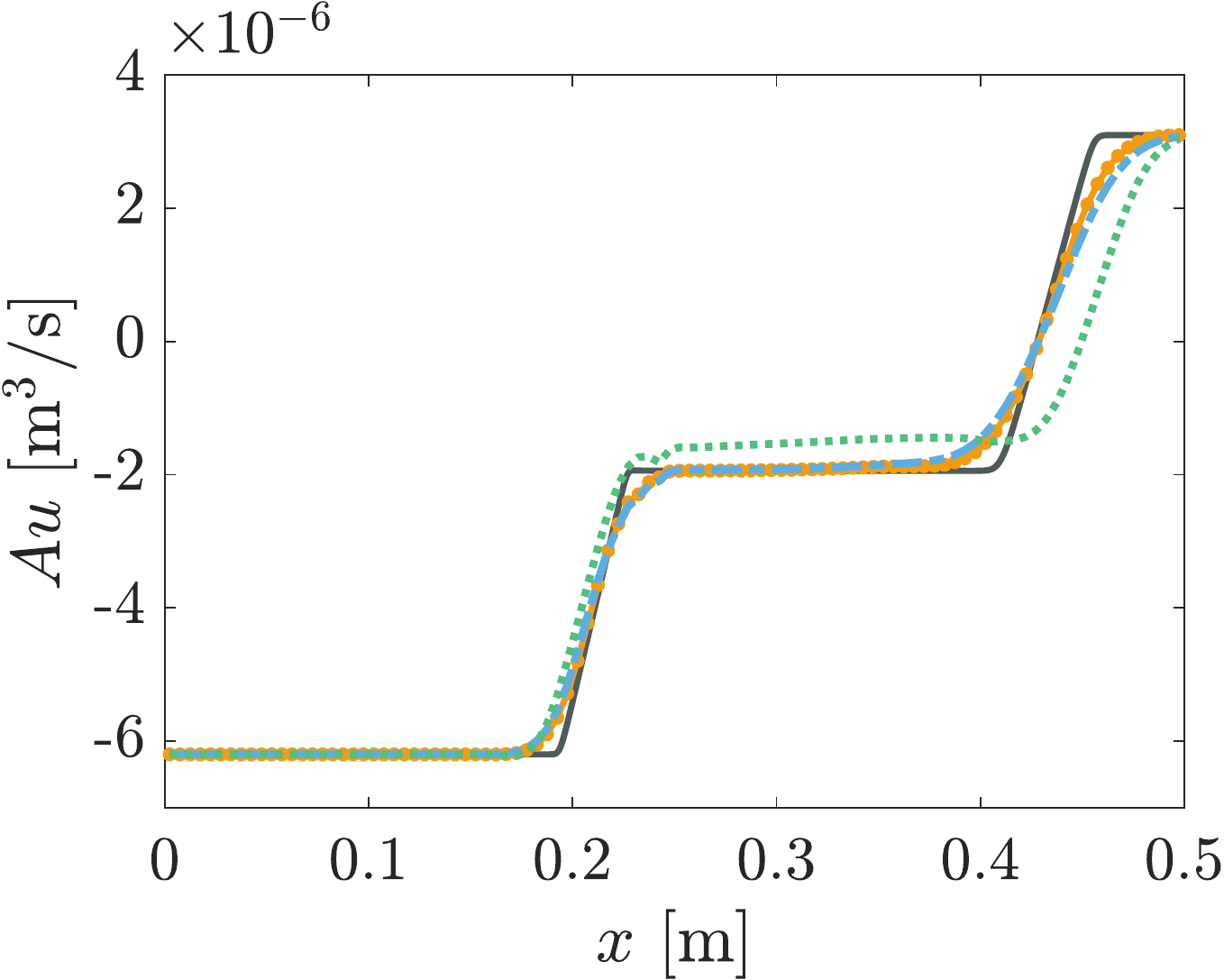}}
\hfil
\subfloat{\includegraphics[width=0.45\textwidth]{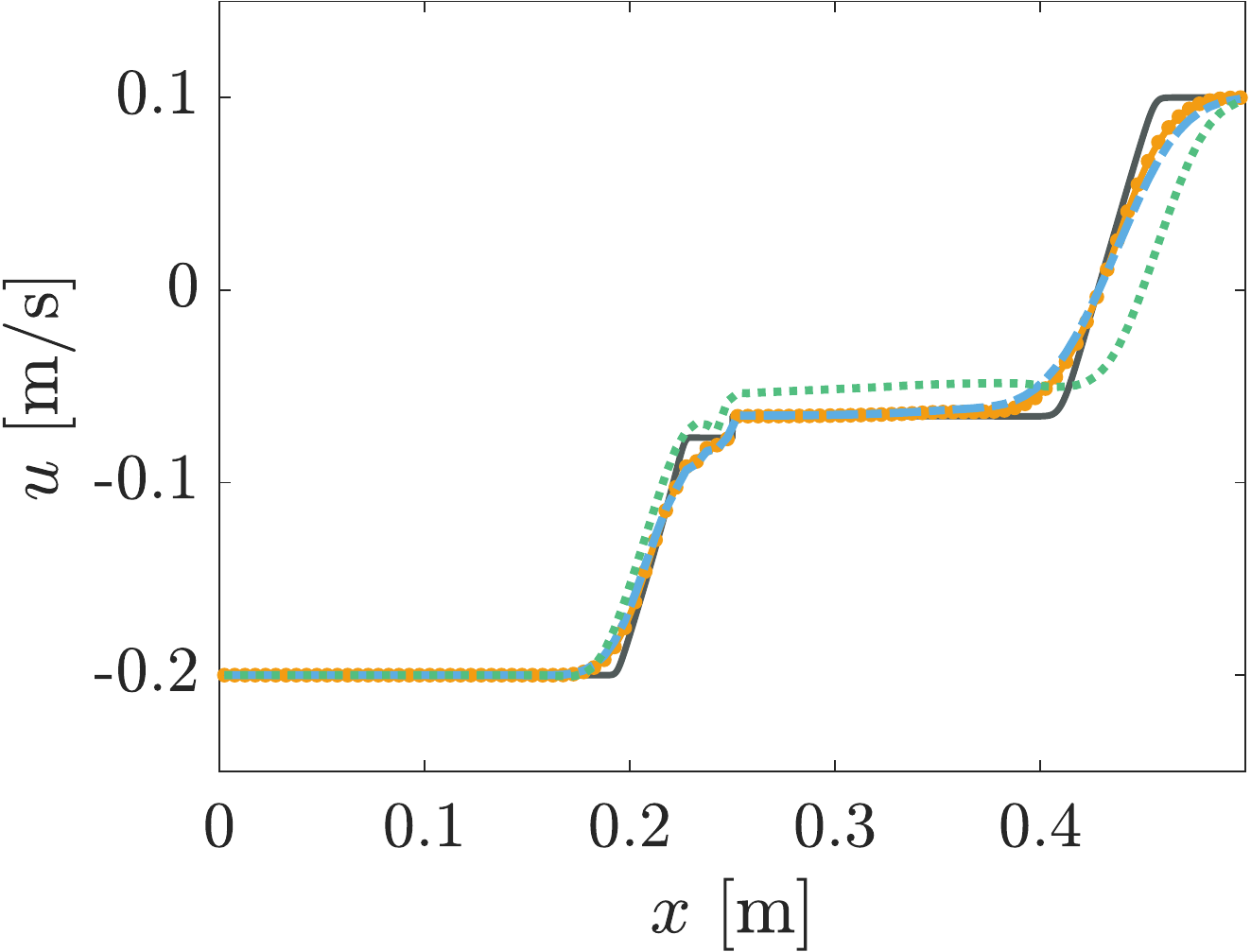}}\\
\subfloat{\includegraphics[width=0.45\textwidth]{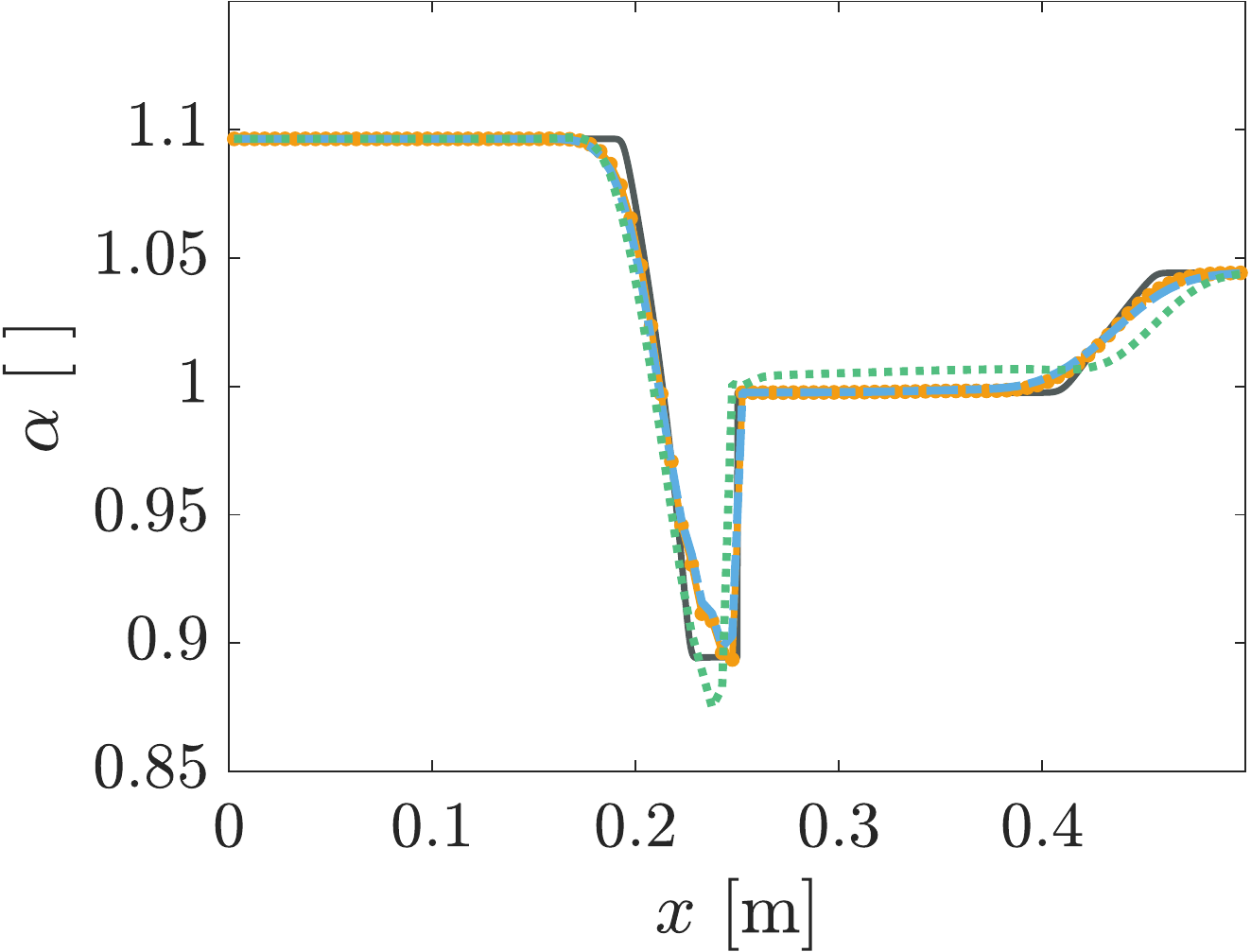}}
\hfil
\subfloat{\includegraphics[width=0.45\textwidth]{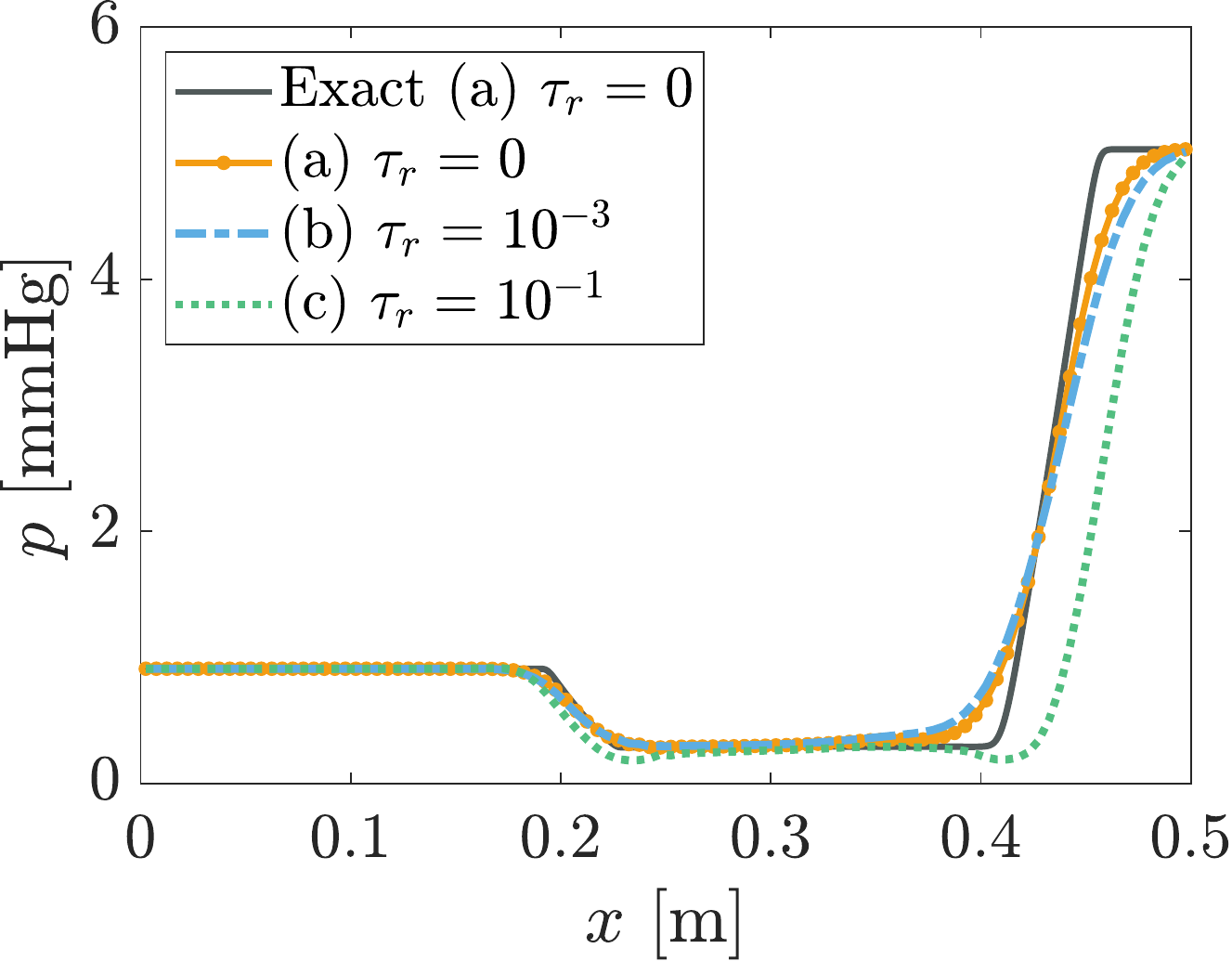}}
\caption{Comparison of the results obtained in RP5 with 3 different choices of the relaxation time $\tau_r$, which correspond to 3 different configurations of the problem (one elastic and two viscoelastic). Numerical results are plotted also against the exact solution of the elastic case ($\tau_r = 0$).}
\label{fig.RP5}
\end{figure}
\paragraph{RP5} In the last Riemann problem (RP5), we consider a generic vein whose wall in the second half is 30 times stiffer than the first half and thus subject to greater internal pressure. Moreover, an initial condition of reflux in the first half of the vessel only is considered. The solution of this problem consists of left and right rarefactions traveling in opposite directions and separated by a stationary contact discontinuity, as depicted in Fig. \ref{fig.RP5}. Once more, the elastic numerical result is in line with the reference solution. Configuration (b) shows only a small damping effect, especially in the right rarefaction wave, while configuration (c) results in a substantial forward shift of the position of the right rarefaction and a raising of the central plateau of flow rate, velocity and area ratio solutions.

\subsection{Multiscale case study: Thoracic aorta with a stent}
As a result of some vascular diseases, a section of an artery might be replaced by a prosthesis or reinforced by the application of a stent (a metal net), causing a sudden change in the mechanical properties of the vessel \cite{ramella2022}.

\begin{figure}[!tb]
\centering
\includegraphics[width=\textwidth]{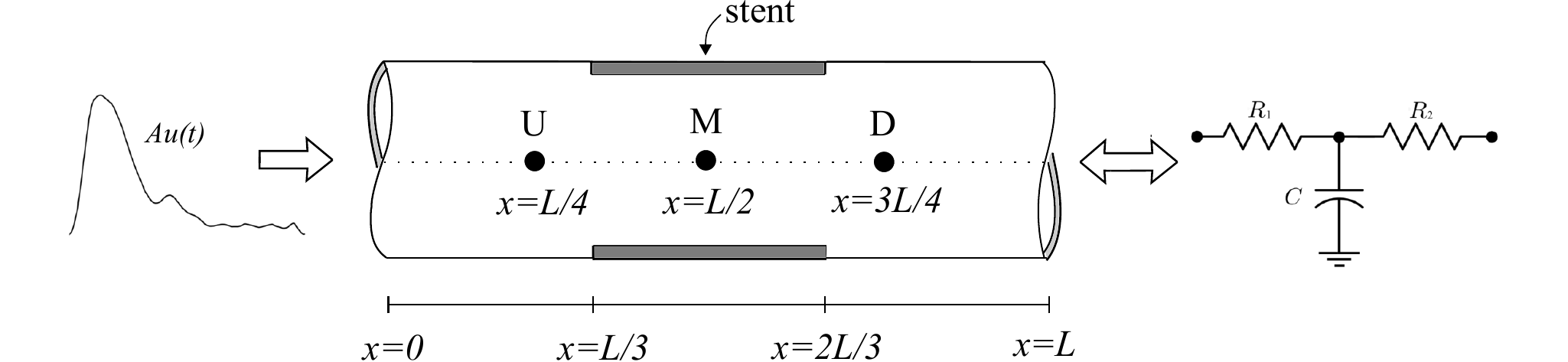}
\caption{Layout of the multiscale test, relating to a generic thoracic aorta with a stent in the middle, with the identification of the 3 control points: upstream of the stent location (U), in the middle of the stent (M) and downstream of the stent (D). Inlet and outlet boundary conditions are prescribed through a given inlet flow rate waveform $Au(t)$ and an RCR model, respectively.}
\label{fig.TCstent_scheme}
\end{figure}

Inspired by \cite{sherwin2003,formaggia2002}, we have designed a test case relating to a section of thoracic aorta in which a stent has been inserted in the middle of it, causing a stiffening of the wall in the center of the vessel. This represents an effective multiscale case study, in which different values of the scaling parameters of the model are attributed to the vessel wall, which lead the system to tend towards the asymptotic diffusive limit (see Section \ref{asymptotics}) in the tract affected by the stent. Indeed, we consider that in the portion of the vessel where the stent is present the Young moduli of the wall are increased by a factor of 100 and the relaxation time of the material is reduced by the same factor. In Table \ref{tab.TCstent}, all the parameters and initial conditions of the test are given, while in Fig. \ref{fig.TCstent_scheme} a schematic representation of the test layout is shown. The parameters of the stentless part of the vessel are set referring to \cite{xiao2014,bertaglia2020a}. Referring to these same works, a realistic input flow rate waveform is considered and RCR model parameters are fixed at the output of the 1D domain, allowing a plausible simulation of the effects of peripheral resistance and compliance on the pulse wave propagation.
 \begin{table}[!htb]
\caption{Initial conditions and mechanical parameters of the multiscale test. For both the sections, the vessel wall thickness is $h_0 = 1.2$ mm. The total length of the artery is $L=24$ cm and the blood density is fixed to be $\rho=1060$ kg/m$^3$. Parameters of the RCR model are $R_1 = 14.047$ MPa$\cdot$s$\cdot$m$^{-3}$, $R_2 = 111.67$ MPa$\cdot$s$\cdot$m$^{-3}$, $C = 14.238$ m$^3\cdot$GPa$^{-1}$.}
\centering
\begin{tabular}{| l | c c |}
\toprule
	 Variable &Stentless tracts &Stented tract\\
\midrule
	$A_0$~[mm$^2$] &452.39 &452.39\\
	$A$~[mm$^2$] &306.04 &450.78\\
	$u$~[m/s] &0.0 &0.0\\
	$p$~[mmHg] &0.0 &0.0\\
	$p_0$~[mmHg] &71.0 &71.0\\
	$E_{\infty}$~[MPa] &0.5333 &53.333\\
	$E_{0}$~[MPa] &0.7619 &76.190\\
	$\eta$~[kPa$\cdot$s] &50.794 &50.794\\
	$\tau_r$~[s] &0.02 &0.0002\\	
\bottomrule
\end{tabular}
\label{tab.TCstent}
\end{table}
\begin{figure}[!tb]
\centering
\subfloat{\includegraphics[width=0.45\textwidth]{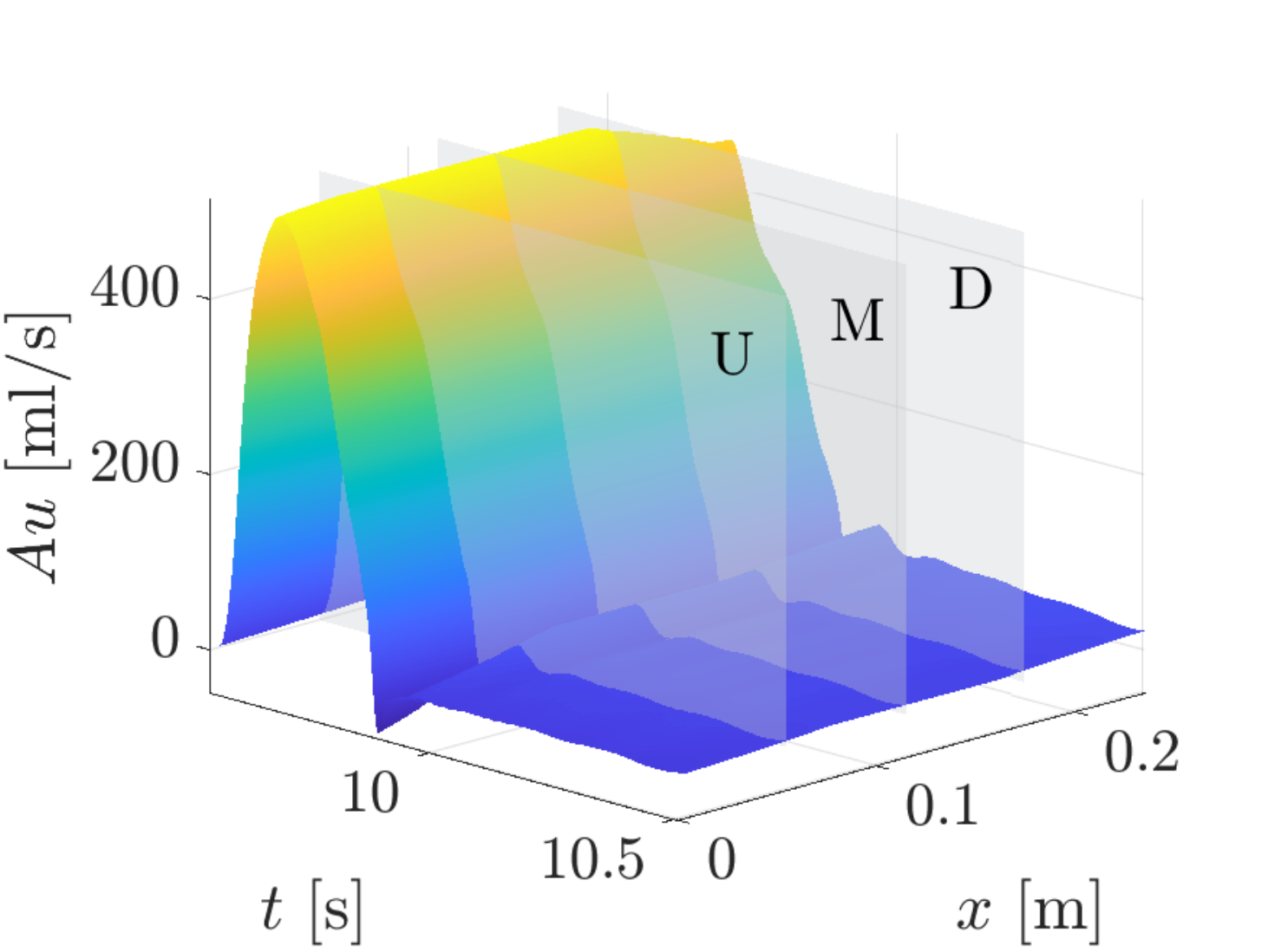}}
\hfil
\subfloat{\includegraphics[width=0.45\textwidth]{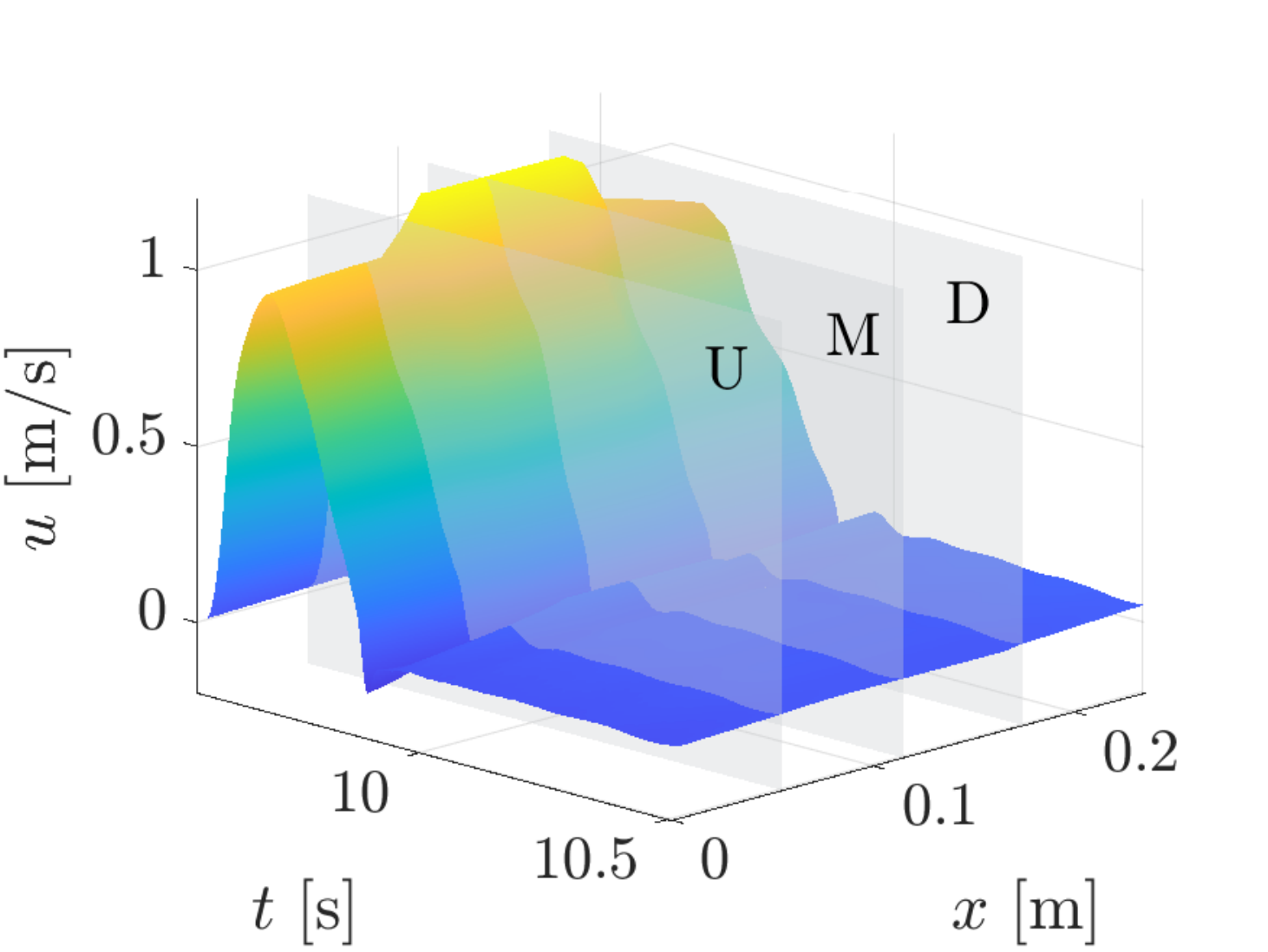}}\\
\subfloat{\includegraphics[width=0.45\textwidth]{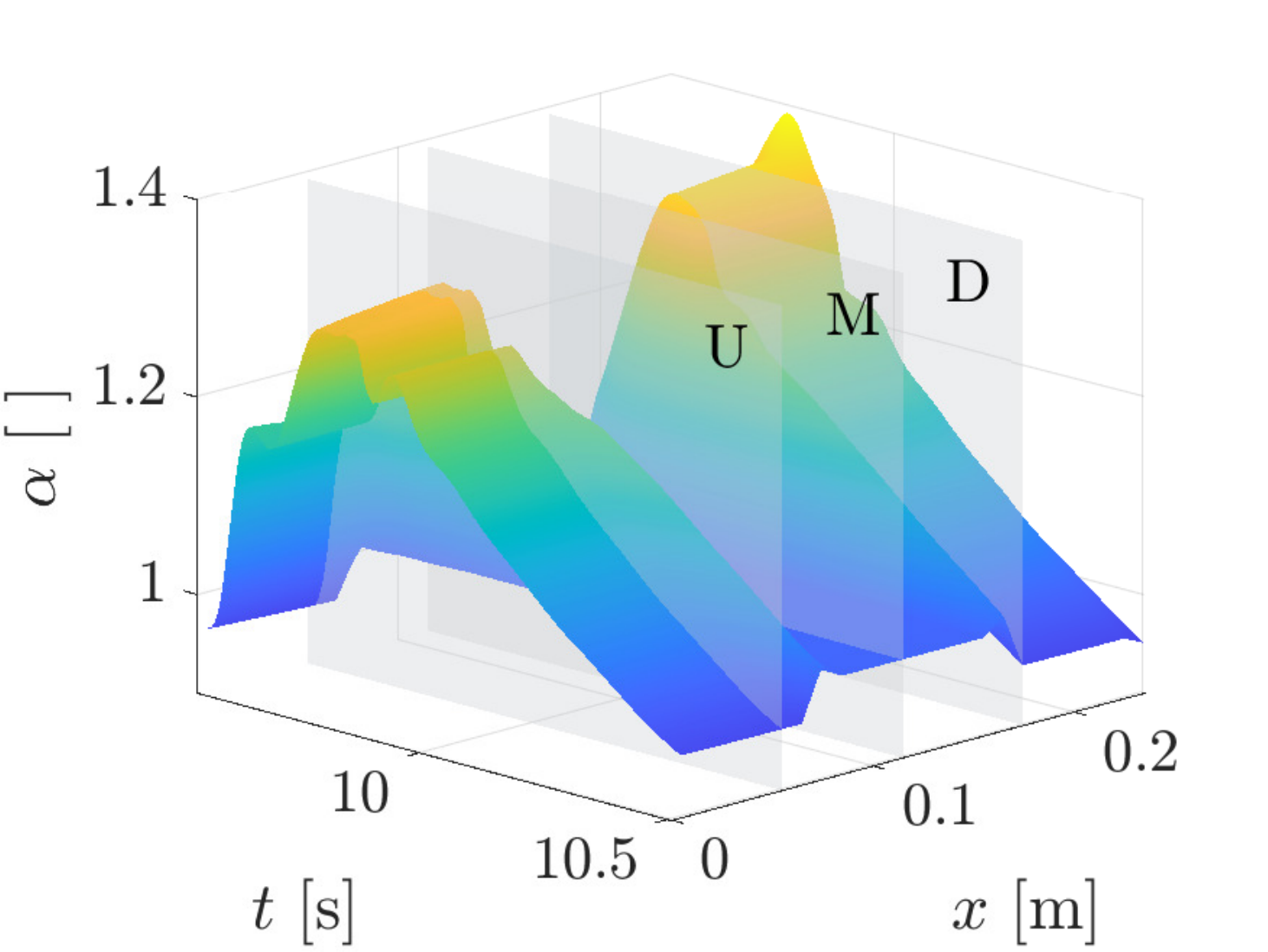}}
\hfil
\subfloat{\includegraphics[width=0.45\textwidth]{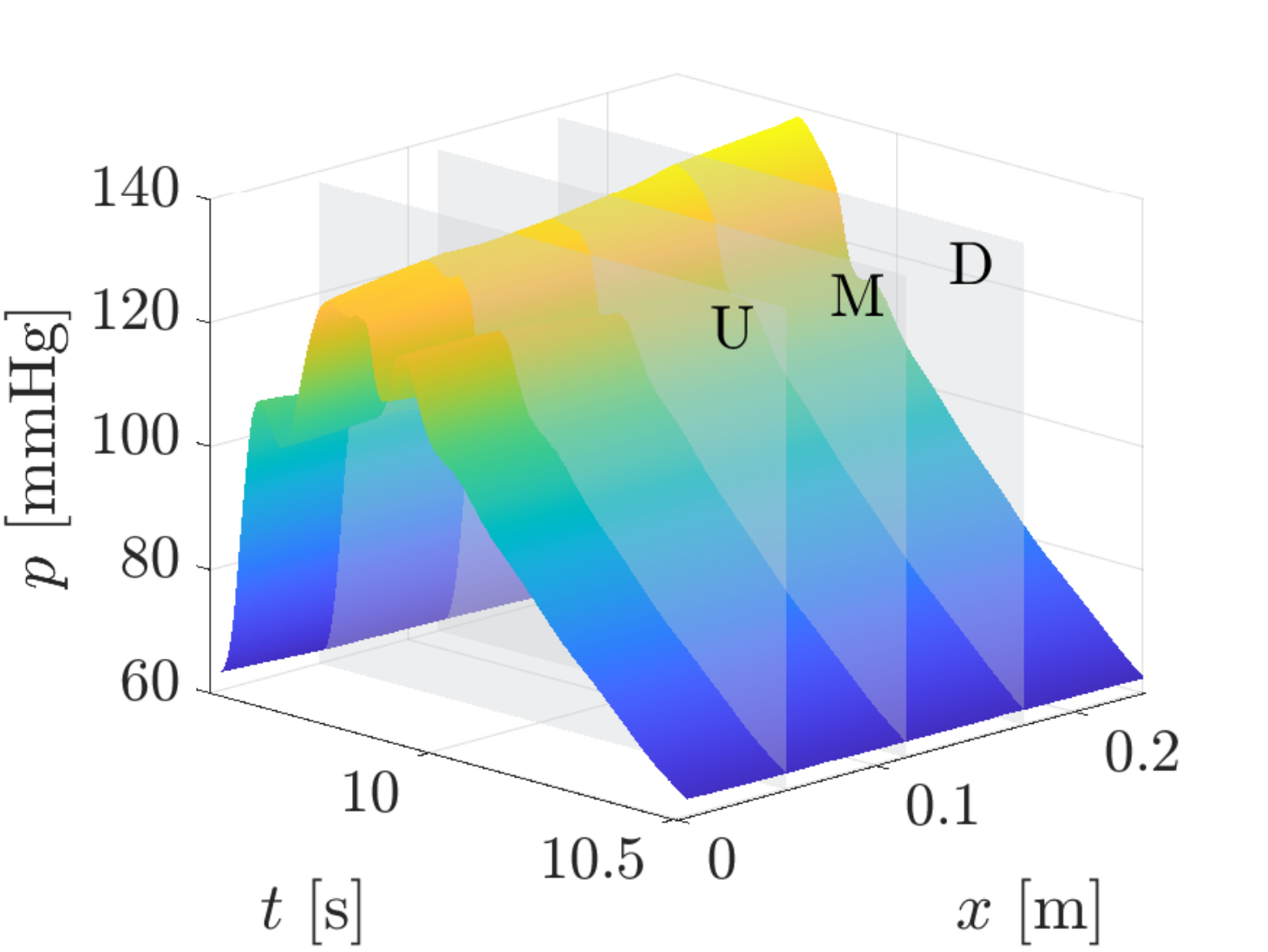}}
\caption{Multiscale test. Space-time evolution in one cardiac cycle of flow rate, velocity, cross-sectional area ratio and pressure in a thoracic aorta having a stent in the middle. The gray planes identify the control sections: upstream of the stent location (U), in the middle of the vessel (M) and downstream of the stent (D).}
\label{fig.TCstent_xt}
\end{figure}
\begin{figure}[!p]
\centering
\captionsetup[subfigure]{labelformat=empty}
\begin{subfigure}{0.31\textwidth}
\centering
\subcaption{$\qquad$Upstream (U)}
\includegraphics[height=3.7cm]{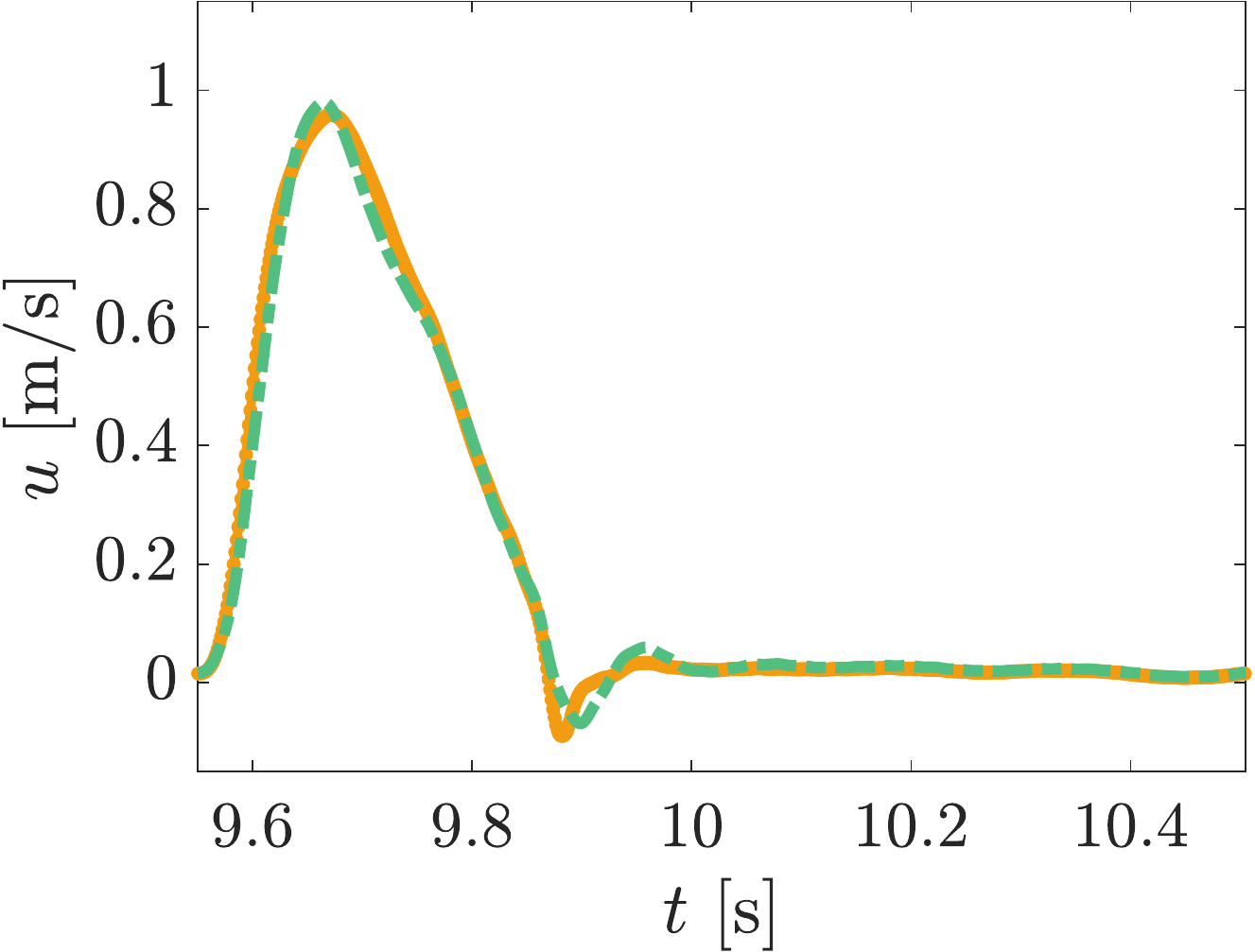}
\end{subfigure}
\hspace{0.1cm}
\begin{subfigure}{0.31\textwidth}
\centering
\subcaption{$\quad$Midpoint (M)}
\includegraphics[height=3.7cm]{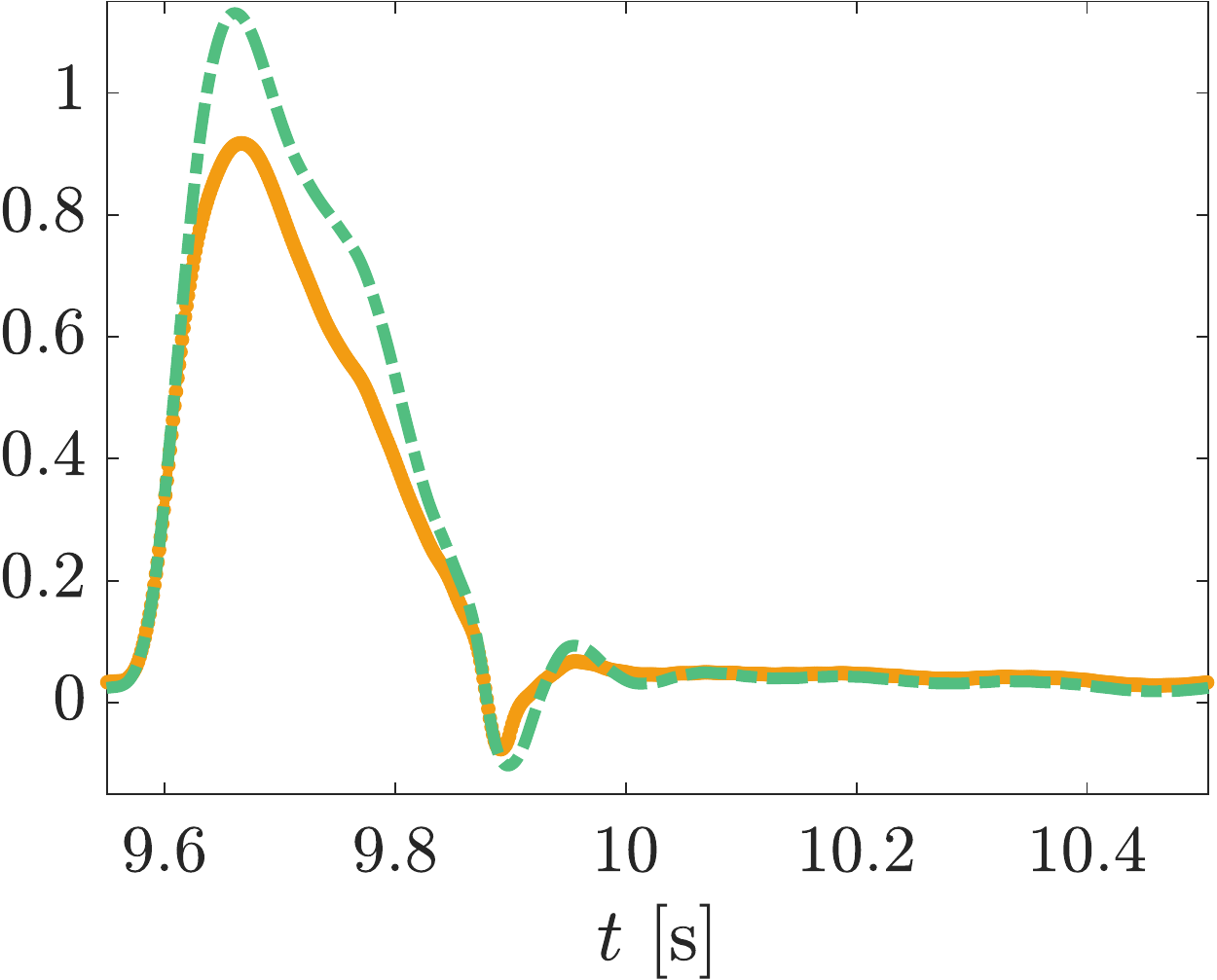}
\end{subfigure}
\hspace{0.1cm}
\begin{subfigure}{0.31\textwidth}
\centering
\subcaption{$\,\,$Downstream (D)}
\includegraphics[height=3.7cm]{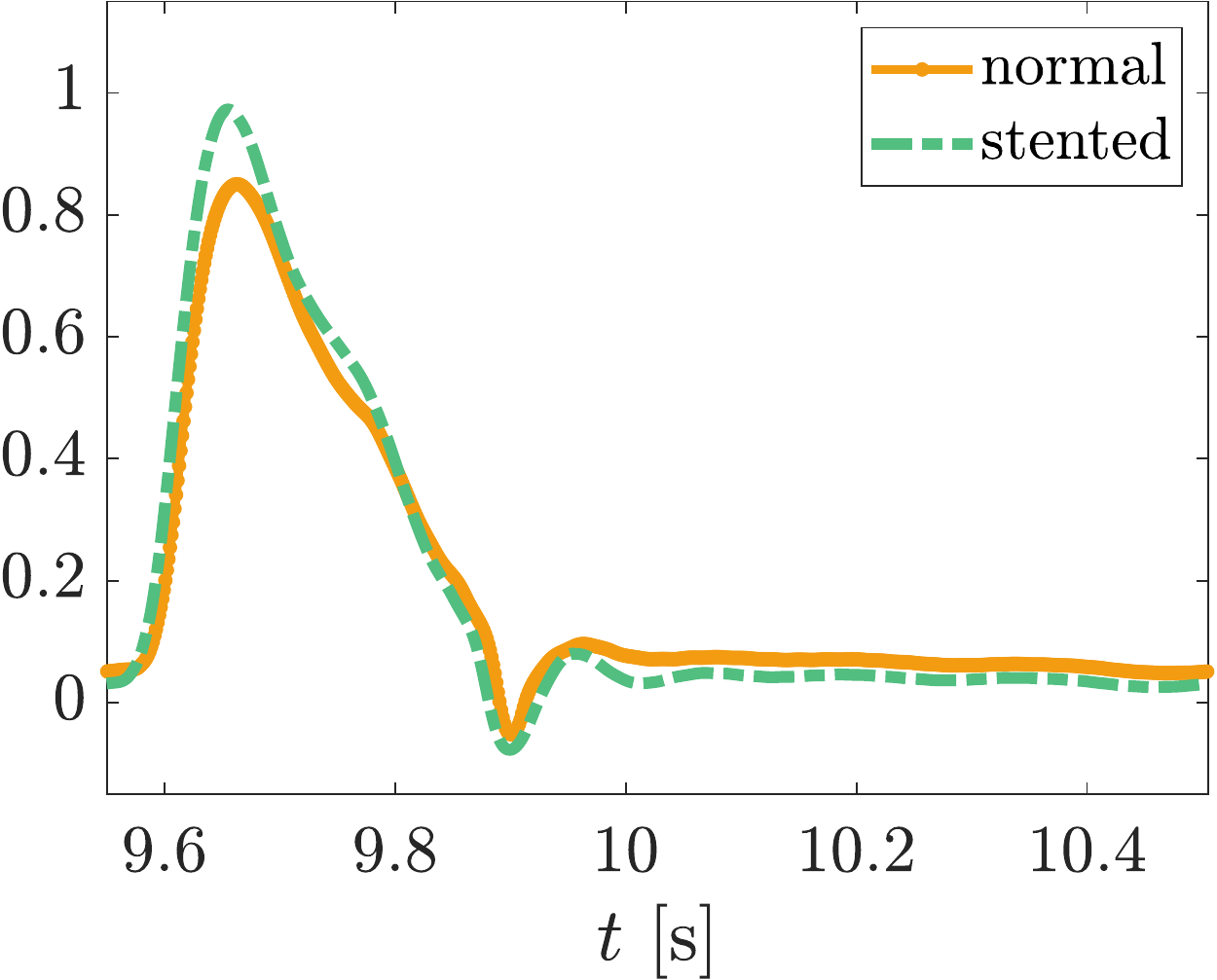}
\end{subfigure}
\begin{subfigure}{0.31\textwidth}
\centering
\includegraphics[height=3.7cm]{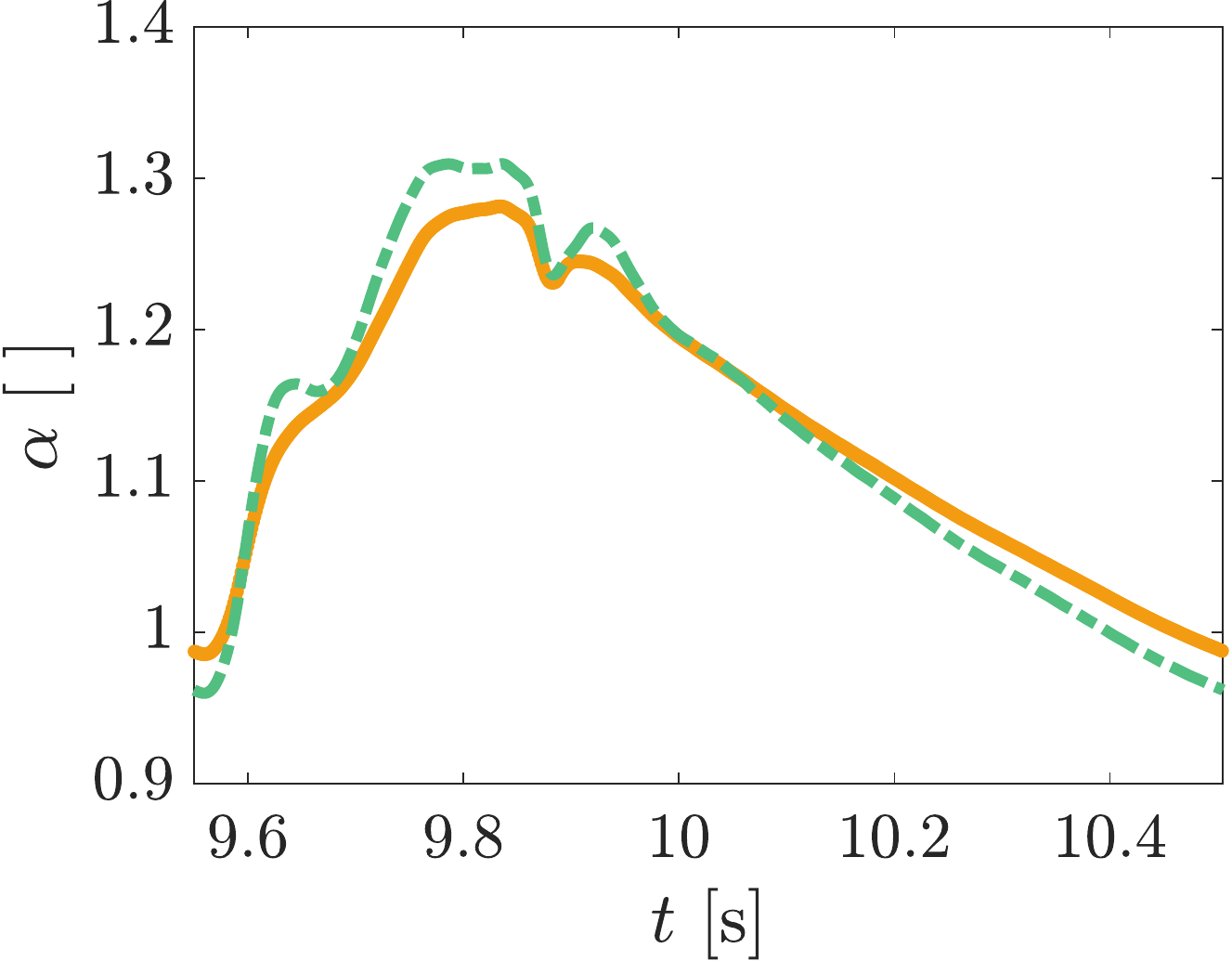}
\end{subfigure}
\hspace{0.1cm}
\begin{subfigure}{0.31\textwidth}
\centering
\includegraphics[height=3.7cm]{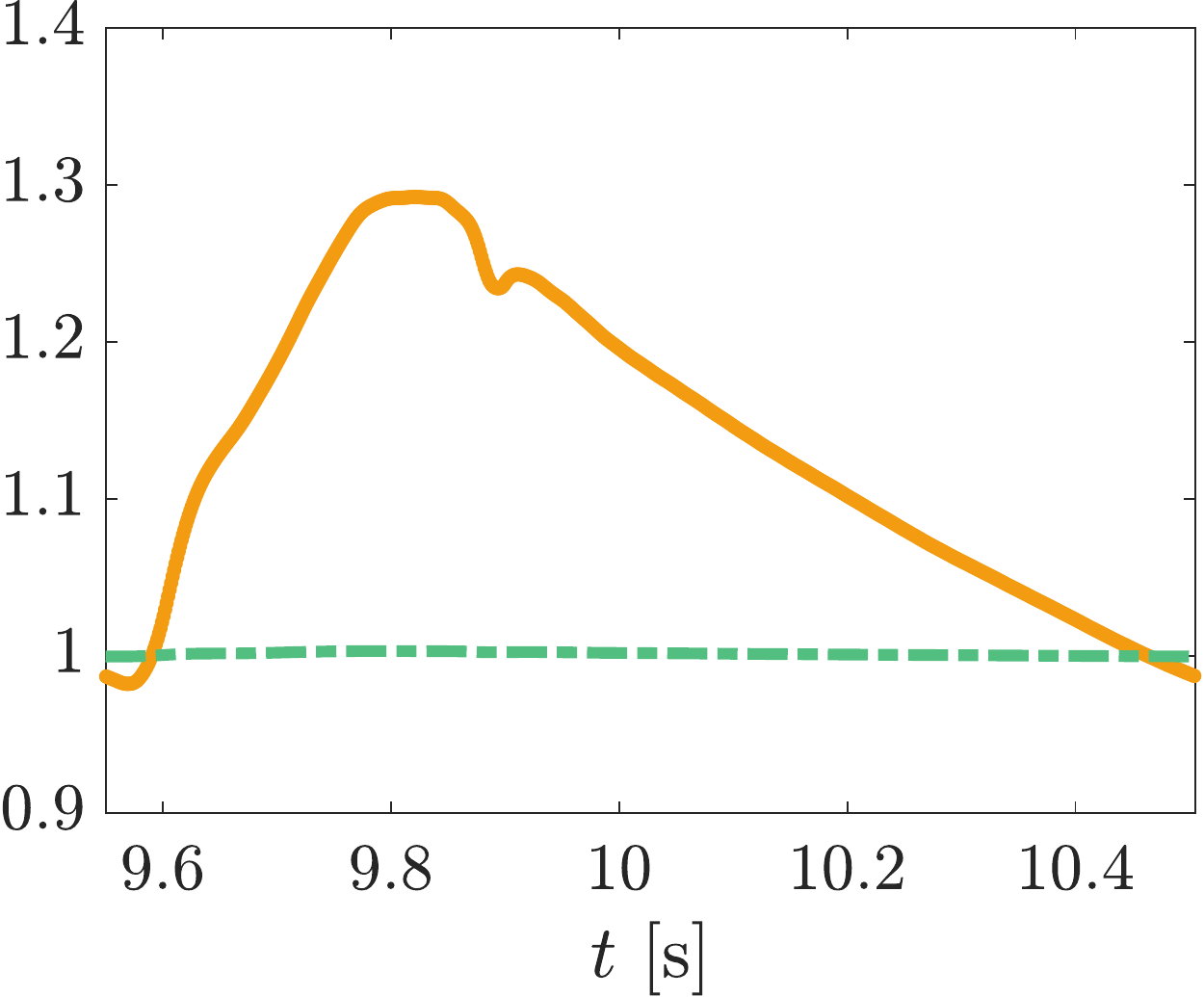}
\end{subfigure}
\hspace{0.1cm}
\begin{subfigure}{0.31\textwidth}
\centering
\includegraphics[height=3.7cm]{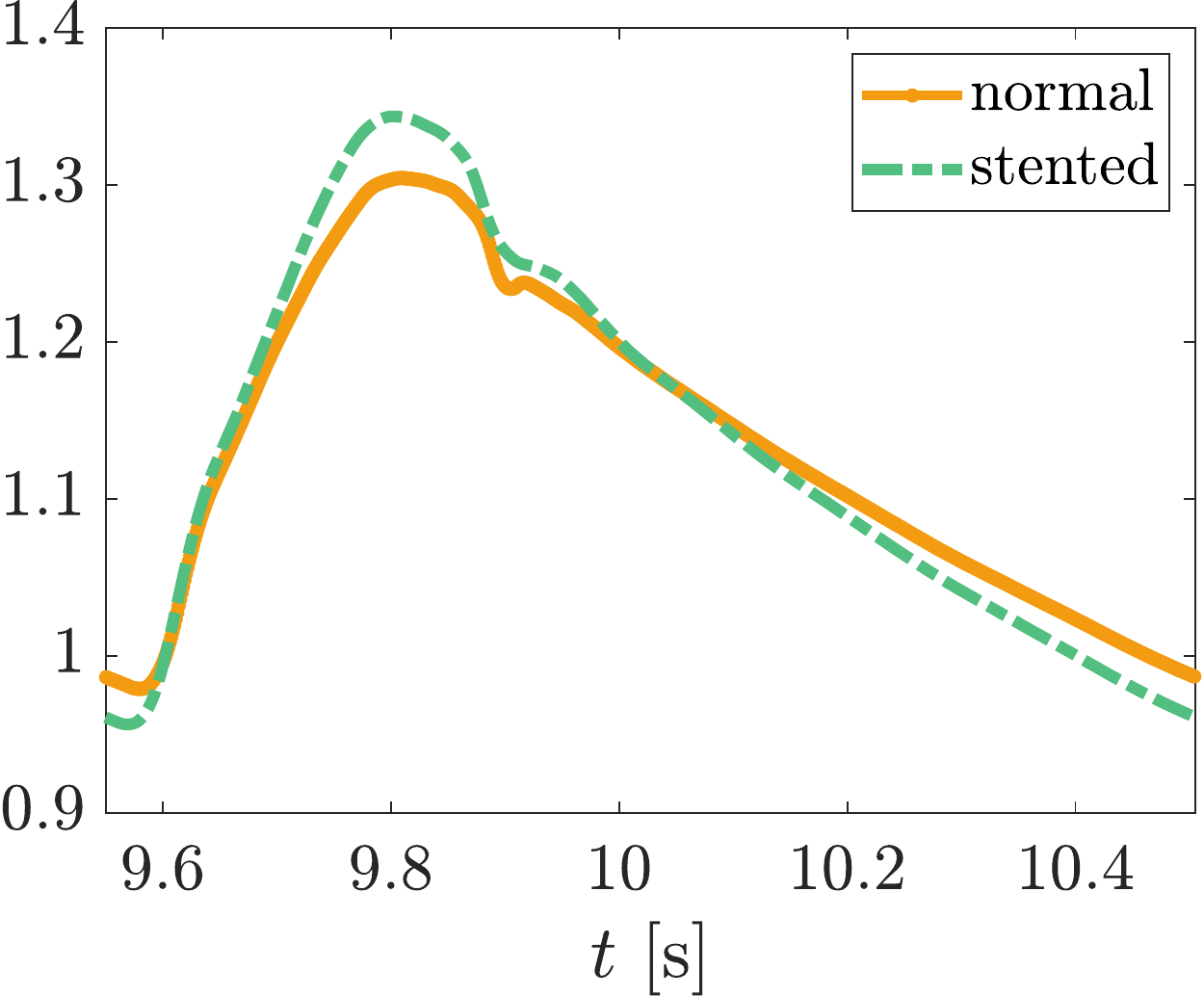}
\end{subfigure}
\begin{subfigure}{0.31\textwidth}
\centering
\includegraphics[height=3.7cm]{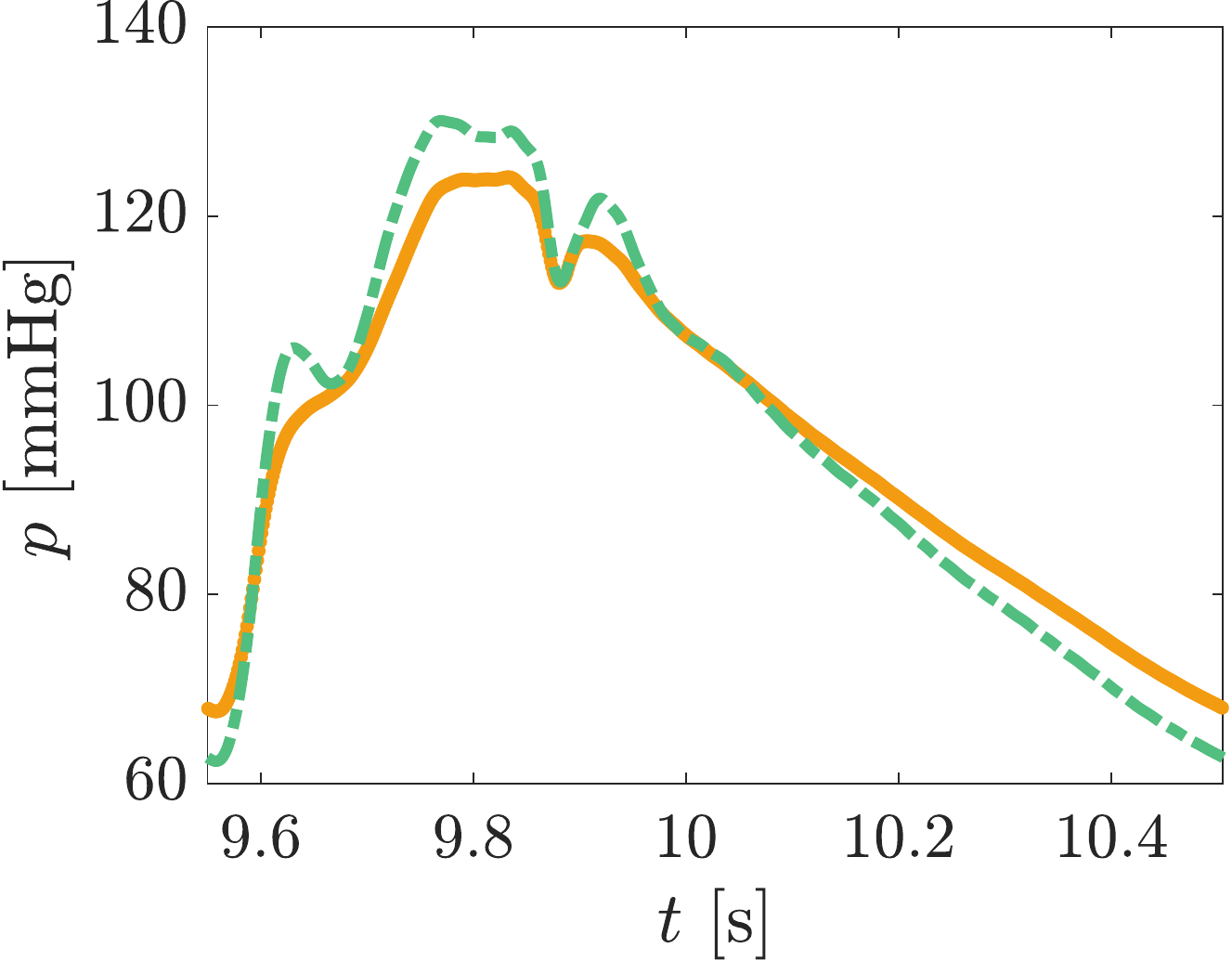}
\end{subfigure}
\hspace{0.1cm}
\begin{subfigure}{0.31\textwidth}
\centering
\includegraphics[height=3.7cm]{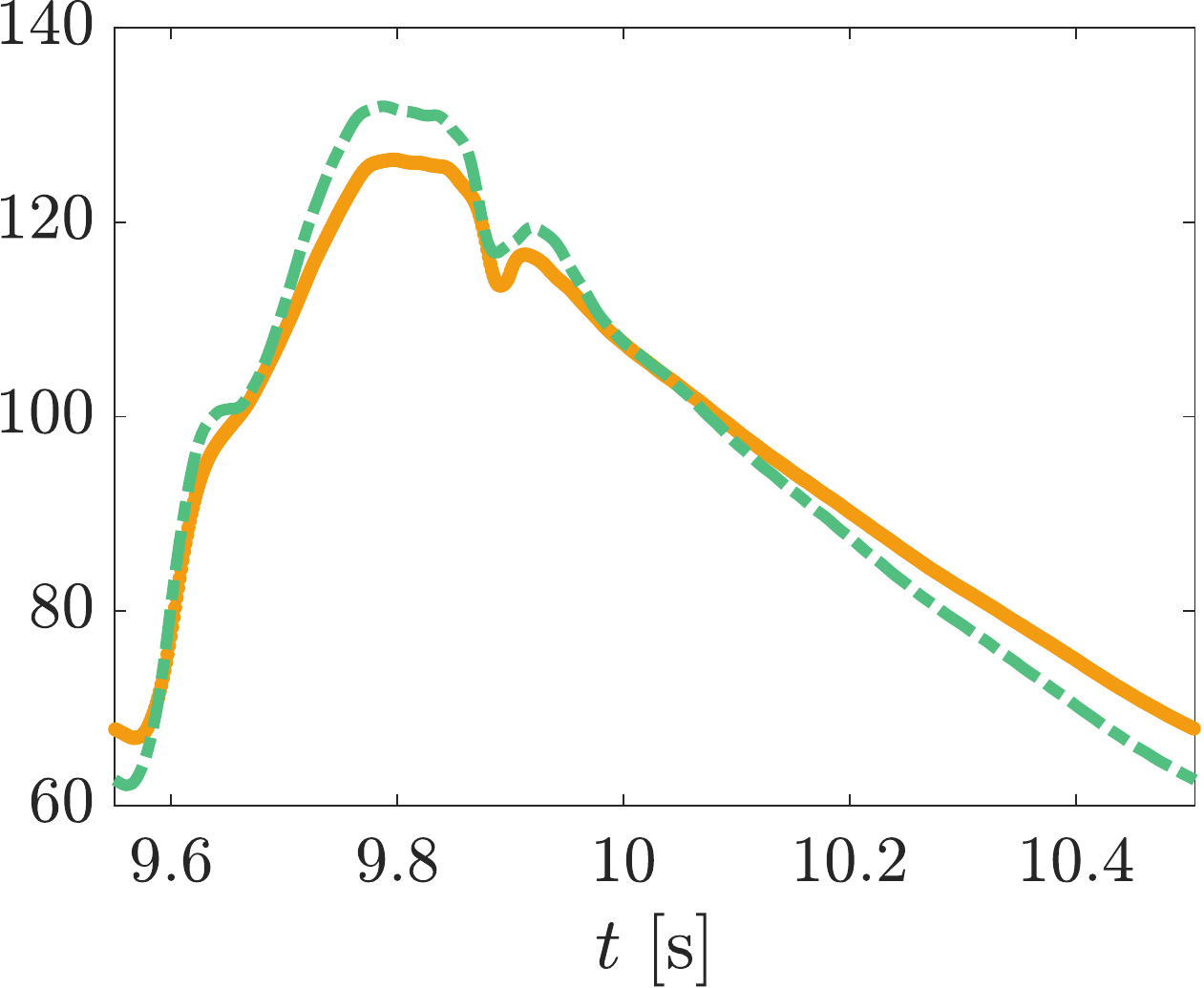}
\end{subfigure}
\hspace{0.1cm}
\begin{subfigure}{0.31\textwidth}
\centering
\includegraphics[height=3.7cm]{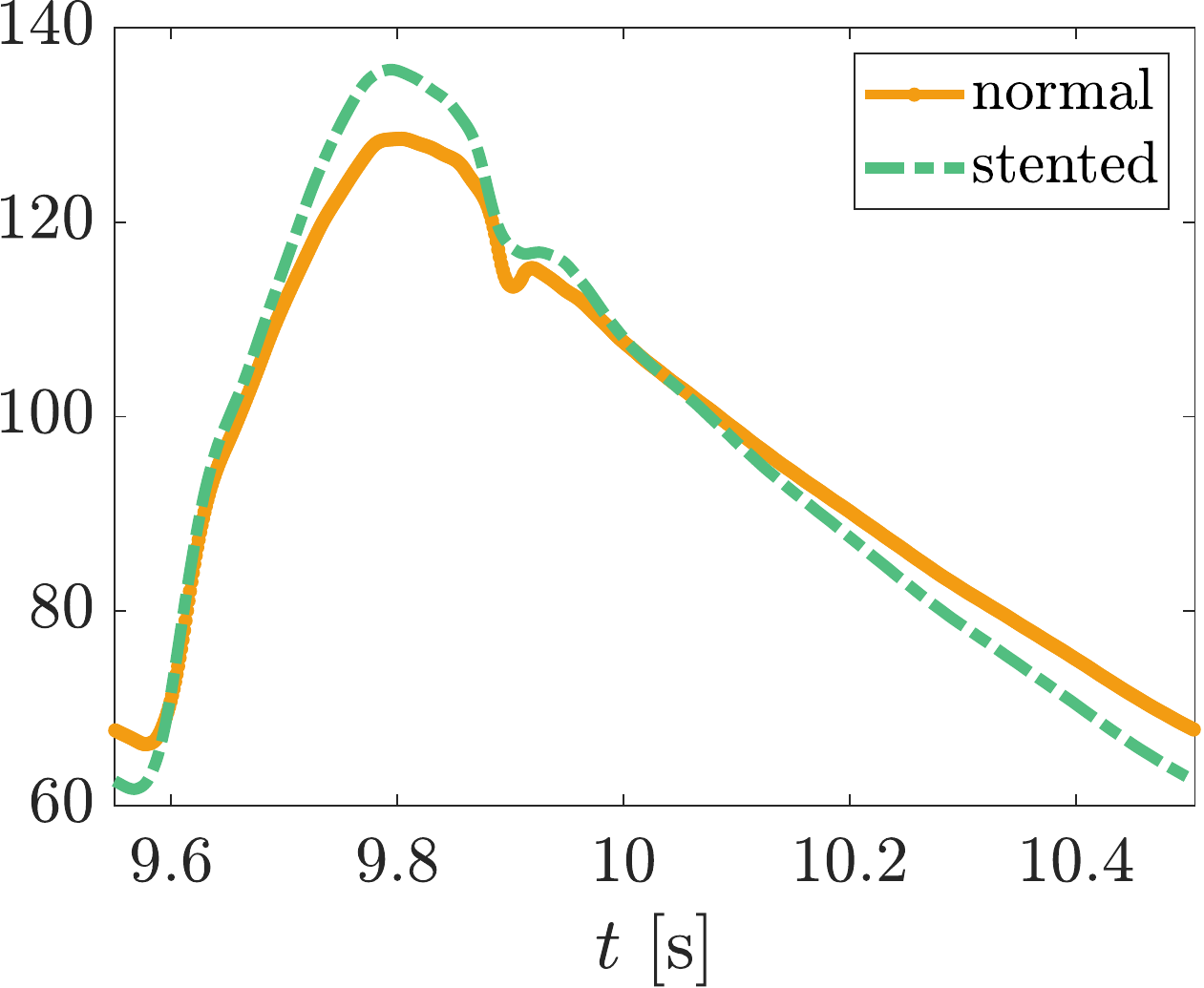}
\end{subfigure}
\begin{subfigure}{0.31\textwidth}
\centering
\includegraphics[height=3.7cm]{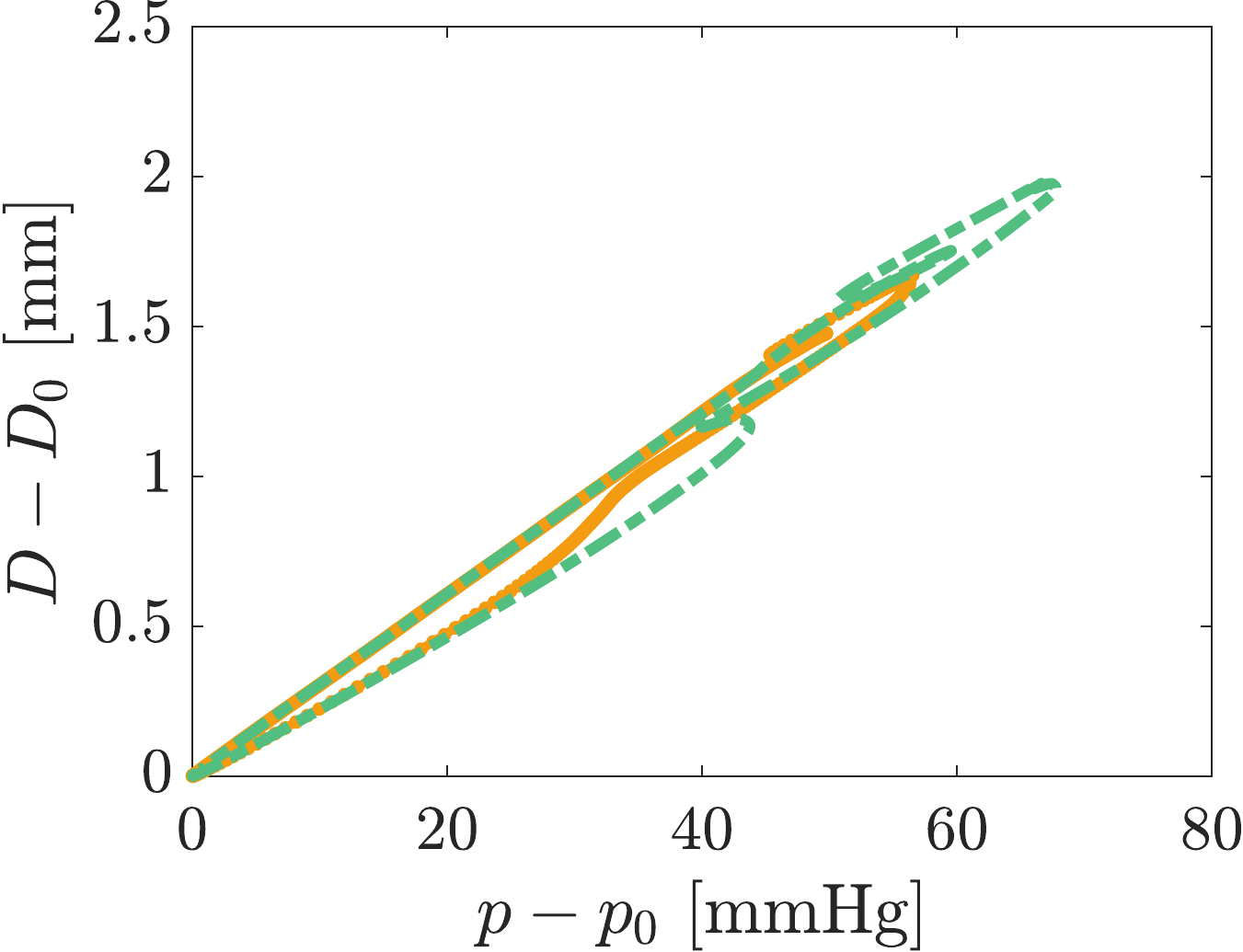}
\end{subfigure}
\hspace{0.1cm}
\begin{subfigure}{0.31\textwidth}
\centering
\includegraphics[height=3.7cm]{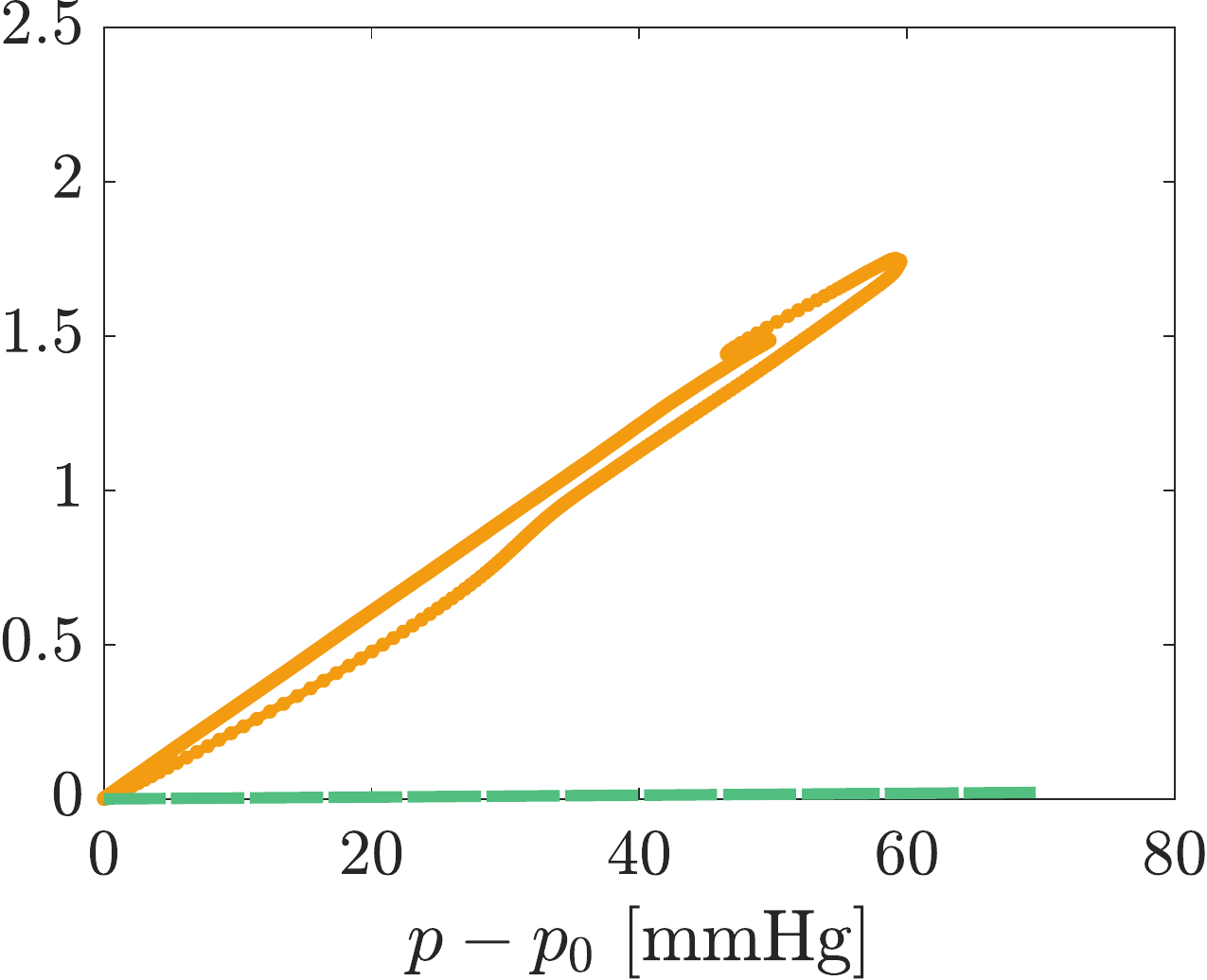}
\end{subfigure}
\hspace{0.1cm}
\begin{subfigure}{0.31\textwidth}
\centering
\includegraphics[height=3.7cm]{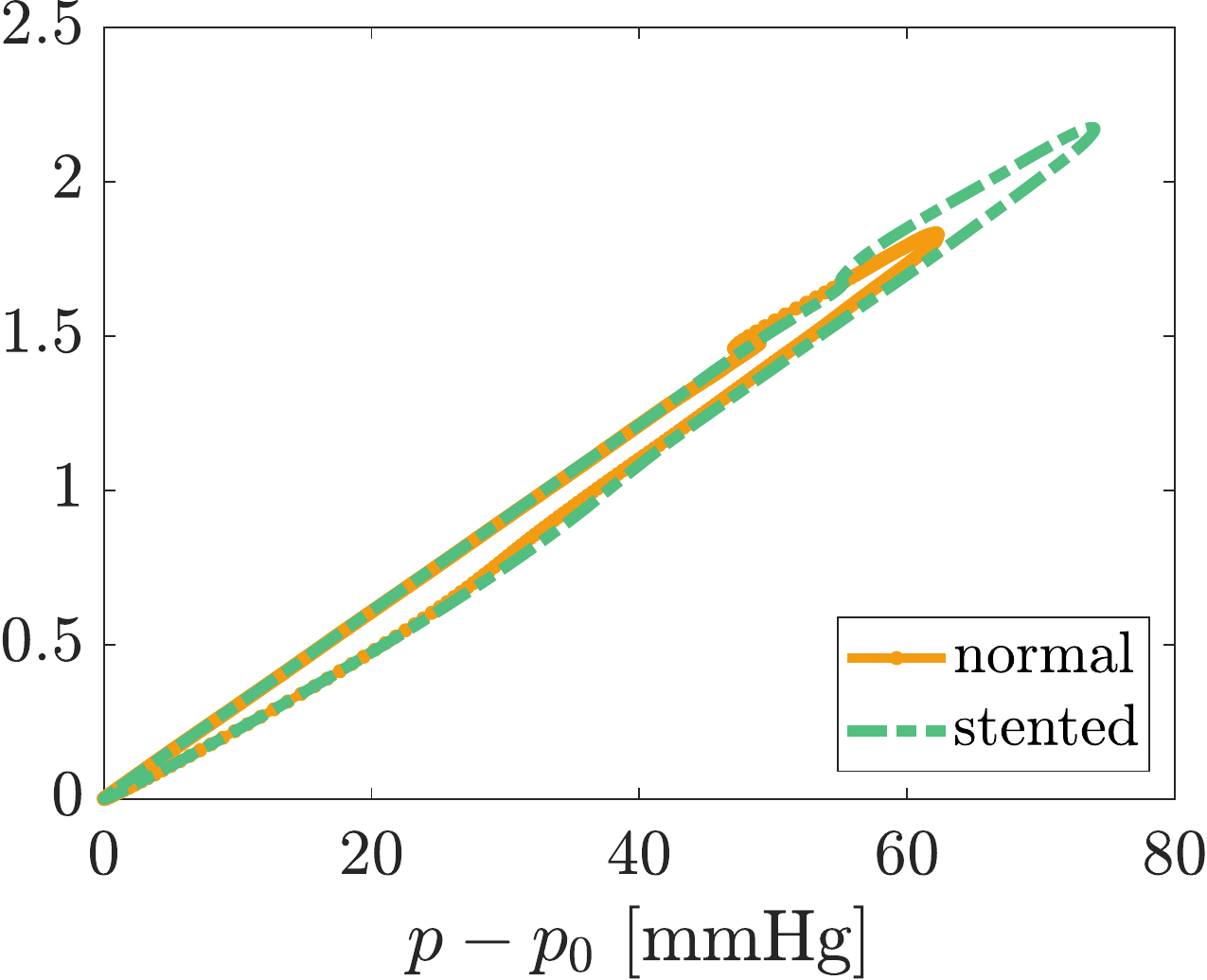}
\end{subfigure}
\caption{Multiscale test. Time evolution in one cardiac cycle of the flow rate, velocity, cross-sectional area ratio, pressure and hysteresis loops in a thoracic aorta under normal and healthy conditions, compared with those obtained in the presence of a stent in the middle of the vessel. Results are shown considering 3 control sections: upstream of the stent location (left column), in the middle of the vessel (middle column) and downstream of the stent (right column).} 
\label{fig.TCstent}
\end{figure}
\begin{figure}[!tbp]
\centering
\captionsetup[subfigure]{labelformat=empty}
\begin{subfigure}{0.31\textwidth}
\centering
\subcaption{$\qquad t = 9.55$ s}
\includegraphics[height=3.7cm]{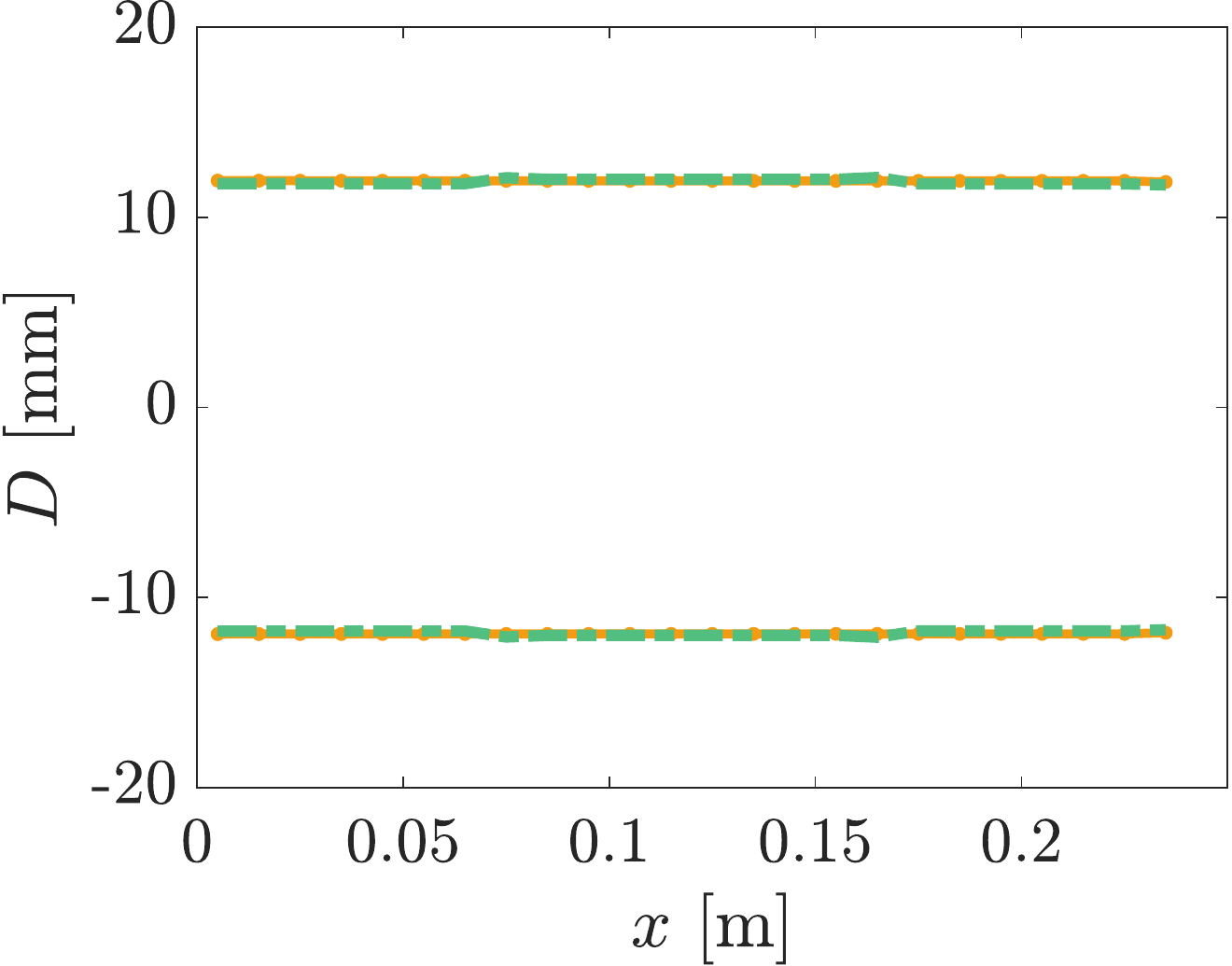}
\end{subfigure}
\hspace{0.1cm}
\begin{subfigure}{0.31\textwidth}
\centering
\subcaption{$\quad t = 9.78$ s}
\includegraphics[height=3.7cm]{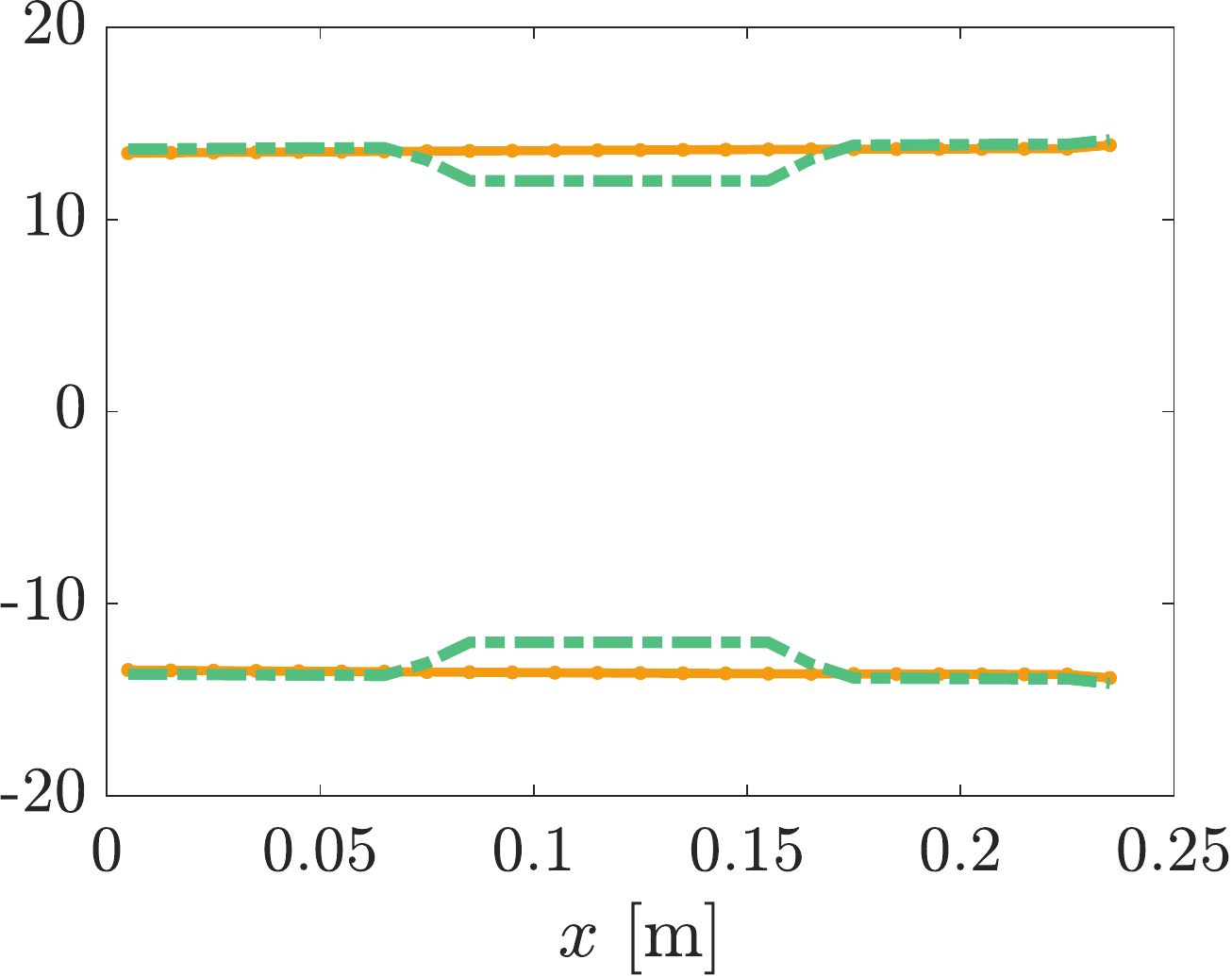}
\end{subfigure}
\hspace{0.1cm}
\begin{subfigure}{0.31\textwidth}
\centering
\subcaption{$\,\, t = 10.00$ s}
\includegraphics[height=3.7cm]{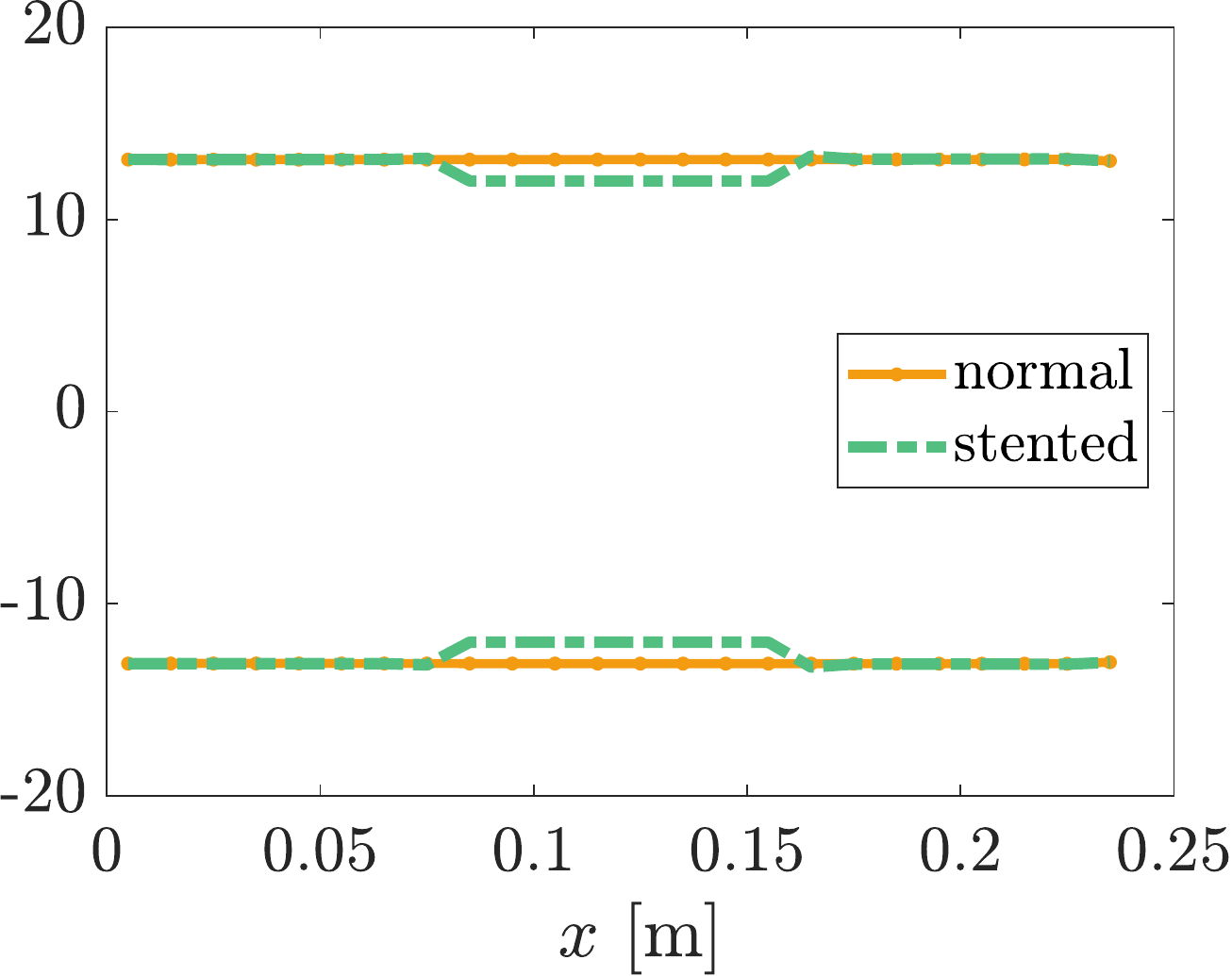}
\end{subfigure}
\caption{Multiscale test. Longitudinal section of a generic thoracic aorta at the beginning of the cardiac cycle (left), at the time of systolic peak (middle) and during the diastolic phase after the dicrotic notch (right) under normal, healthy conditions, compared with that in which a stent is placed in the centre of the artery.}
\label{fig.TCstent_D}
\end{figure}

The simulation has been run with $N_x=24$ computational cells for 10 cardiac cycles (duration $0.955\, \mathrm{s} = 63$ bpm each), both in presence of the stent and under normal, healthy conditions, to allow comparisons. In Figure \ref{fig.TCstent_xt}, the space-time solution in the stented configuration is shown for one cardiac cycle for the main variables. In the same figure the control sections are also indicated, corresponding to those depicted in the layout in Figure \ref{fig.TCstent_scheme}. In Figure \ref{fig.TCstent}, the time evolutions in one cardiac cycle of the main variables and hysteresis curves are presented for both the normal and stented configurations with respect to the 3 control sections. As expected, the presence of the stent significantly alters the blood flow propagation, creating, in particular, wave reflections due to the impact with the stiffer region, which cause an increase in pressure peaks \cite{formaggia2002}, especially visible in the upstream pressure plot. It can be observed that also the velocity undergoes an increase in the systolic phase, mostly in the stented tract. On the other hand, because of the greater stiffness of the stent, changes in cross-sectional area from equilibrium cannot be appreciated in the central region with respect to those produced by the original vessel wall (see $\alpha$ and hysteresis midpoint plots). In this regard, Figure \ref{fig.TCstent_D} compares the longitudinal section of the vessel without and with the stent at three different temporal instants: at the beginning of the cardiac cycle, at the time of the systolic peak, and during the diastolic phase (after the dichrotic notch), highlighting the almost zero stretchability of the stent compared with that of the regions without the stent \cite{sherwin2003}.

\section{Conclusions}
\label{section_conclusions}
In the present study, we introduce a multiscale constitutive framework for the purpose of modeling one-dimensional blood flow. 
We show that the proposed augmented model, which takes into account a linear viscoelastic constitutive characterization of the fluid-structure interaction occurring between the vessel wall and the blood flow, can describe different physical propagation phenomena ranging from hyperbolic transport to parabolic diffusion, recovering alternative rheological behaviors of blood vessels. 
This highly flexible, yet accurate, mathematical model is ideally suited for modeling the complex cardiovascular system, which is characterized by significant morphological and mechanical variability. We further derive a novel viscoelastic constitutive model by analyzing the perturbation of the local elastic equilibrium of the system allowing us to capture the second order small viscosity effects.

To solve the resulting multiscale hyperbolic system, we employ a state-of-the-art third-order asymptotic-preserving IMEX Runge-Kutta finite volume method that ensures consistency of the numerical scheme with the asymptotic limits of the mathematical model. Notably, our method enables us to choose a time step size that is not affected by restrictions related to the smallness of the scaling parameters and fulfills the well-balance property in time.
Several numerical tests confirm the validity of the approach, including a case study related to the hemodynamics of a thoracic aorta in the presence of a stent.

Further research will be directed toward a more in-depth treatment of the augmented blood flow model coupled with the new viscoelastic constitutive law derived from local elastic equilibrium perturbation. Additionally, we will explore the extension of our methodology to the main cardiovascular network \cite{piccioli2021}, as well as investigate the sensitivity of the model with respect to mechanical scaling parameters through uncertainty quantification approaches \cite{bertaglia2021a,bertaglia2022}.

\appendix
\section{Numerical implementation aspects}
\label{section_appendix}
\subsection{Dimensionless form of the model}
Due to the presence of variables in system \eqref{completesyst} that have very different orders of magnitude, a switch to the dimensionless form of the equations is necessary to avoid numerical fluctuations, especially when using high-order methods \cite{pimentel-garcia2023}. Fixing the characteristic values for length ($\bar L$), time ($\bar T$), blood density ($\bar \rho$), cross-sectional vessel area ($\bar A$), viscosity ($\bar \eta$), Young modulus ($\bar E = \bar \eta / \bar T$), and velocity $\bar U = \bar L/\bar T$, the following dimensionless variables are introduced:
\[ 
x^* = \frac{x}{\bar L}\,, \quad t^* = \frac{t}{\bar T}\,, \quad \rho^* = \frac{\rho}{\bar \rho}\,, \quad A^* = \frac{A}{\bar A}\,, \quad u^* = \frac{u}{\bar U}\,, \quad p^* = \frac{p}{\bar \rho \bar U^2}\,, 
\]
\[ 
A_0^* = \frac{A_0}{\bar A}\,, \quad E_0^* = \frac{E_0}{\bar E}\,, \quad E_{\infty}^* = \frac{E_{\infty}}{\bar E}\,, \quad p_0^* = \frac{p_0}{\bar \rho \bar U^2}\,, \quad \tau_r^* = \frac{\tau_r}{\bar T}\,.
\]
Thus, system \eqref{completesyst} can be written as:
\begin{subequations}
\begin{align}
	&\frac{\partial A^*}{\partial t^*} + \frac{\partial(A^*u^*)}{\partial x^*} = 0 \\
	&\frac{\partial(A^*u^*)}{\partial t^*}+ \frac{\partial((A^*u^*)^{2}/A^*)}{\partial x^*}  + \frac{A^*}{\rho^*} \, \frac{\partial p^*}{\partial x^*} = 0 \\
	&\frac{\partial p^*}{\partial t^*} + \frac1{\tilde{\mathsf{Re}}} \frac{E_0^*}{WA^*}\left( m\alpha^m - n\alpha^n\right) \,\frac{\partial(A^*u^*)}{\partial x^*} = -\frac{1}{\tau_r^*}\left( p^* - p_0^* - \frac1{\tilde{\mathsf{Re}}} \frac{E_{\infty}^*}{W}\left( \alpha ^m - \alpha ^n \right) \right) .
\end{align}
\label{completesyst_dimensionless}
\end{subequations}
It is worth to notice that here $\tilde{\mathsf{Re}} = \bar \rho \bar U \bar L / \bar \eta$ is the Reynolds number accounting for the viscosity $\bar \eta$ of the wall and not the viscosity of the fluid (in contrast with the classical definition).
Finally, in the above system we consider $\bar L = L$ m (length of the domain), $\bar T = 1$ s, $\bar \rho = 1050$ kg/m$^3$, $\bar A = \mathrm{mean}[A(x,0)]$ m$^2$, and $\bar E = \mathrm{mean}[E_0(x,0);E_{\infty}(x,0)]$ Pa.
\subsection{Third order IMEX method}
Following the third-order GSA BPR(3,4,3) scheme proposed in \cite{boscarino2017}, which is characterized by $s=4$ stages for the implicit part and 3 stages for the explicit part, the Butcher tableaux we employ (explicit on the left and implicit on the right) are:
\begin{equation}
\begin{tabular}{c | c c c c c}
0 & 0 & 0 & 0 & 0 & 0 \\
1 & 1 & 0 & 0 & 0 & 0 \\
2/3 & 4/9 & 2/9 & 0 & 0 & 0 \\
1 & 1/4 & 0 & 3/4 & 0 & 0 \\
1 & 1/4 & 0 & 3/4 & 0 & 0 \\ \hline
  & 1/4 & 0 & 3/4 & 0 & 0 
\end{tabular}
\hspace{0.5cm}
\begin{tabular}{c | c c c c c}
0 & 0 & 0 & 0 & 0 & 0 \\
1 & 1/2 & 1/2 & 0 & 0 & 0 \\
2/3 & 5/18 & -1/9 & 1/2 & 0 & 0 \\
1 & 1/2 & 0 & 0 & 1/2 & 0 \\
1 & 1/4 & 0 & 3/4 & -1/2 & 1/2 \\ \hline
  & 1/4 & 0 & 3/4 & -1/2 & 1/2  
\end{tabular}
\label{eq:tableaux}
\end{equation}

\section*{Acknowledgements}
This work was partially supported by MIUR (Ministero dell'Istruzione, dell'Universit\`a e della Ricerca) PRIN 2017 for the project \textit{``Innovative numerical methods for evolutionary partial differential equations and applications''}, code 2017KKJP4X. G.B. was also funded under ``Bando Giovani anno 2022 per progetti di ricerca finanziati con contributo 5x1000 anno 2020'' by the University of Ferrara, and acknowledges support from GNCS--INdAM under the Project E53C22001930001.\\

\bibliographystyle{abbrv}
\bibliography{BloodFlow_AsymptoticLimit_final}

\end{document}